\documentclass[a4paper,12pt]{article} 
\usepackage{amssymb}
\usepackage{amscd}
\usepackage{amsmath} 
\usepackage[dvipdfmx]{graphicx}
\usepackage{latexsym} 
\usepackage{enumerate}
\usepackage[usenames]{color}
\usepackage{lscape}

\usepackage{theorem}
\theoremstyle{plain}
\theorembodyfont{\itshape}
\newtheorem{theorem}{Theorem}[section]
\newtheorem{proposition}[theorem]{Proposition}
\newtheorem{lemma}[theorem]{Lemma}
\newtheorem{corollary}[theorem]{Corollary}

\theorembodyfont{\rmfamily}
\newtheorem{definition}[theorem]{Definition}
\newtheorem{remark}[theorem]{Remark}
\newtheorem{problem}[theorem]{Problem}

\newtheorem{algorithm}[theorem]{Algorithm}

\newtheorem{problemgpf}{Problem I}

\newtheorem{problemocf}{Problem I\!I}

\def\rA{\mathrm{A}}
\def\rB{\mathrm{B}}
\def\rI{\mathrm{I}}
\def\rII{\mathrm{I\!I}}
\def\rIm{\mathrm{Im}}
\def\rRe{\mathrm{Re}} 
\def\C{\mathbb{C}}
\def\N{\mathbb{N}}
\def\Q{\mathbb{Q}}
\def\R{\mathbb{R}}
\def\Z{\mathbb{Z}}
\def\cD{\mathcal{D}}
\def\cE{\mathcal{E}}
\def\cF{\mathcal{F}}

\def\cI{\mathcal{I}}
\def\ba{\mbox{\boldmath $a$}}

\def\bp{\mbox{\boldmath $p$}}
\def\bv{\mbox{\boldmath $v$}}
\def\bal{\mbox{\boldmath $\alpha$}}
\def\1{\mbox{\boldmath $1$}}
\def\hgF{{}_2F_1}
\def\hgG{{}_2G_1}
\def\hgH{{}_2H_1}
\def\hgK{{}_2K_1}
\def\ts{\textstyle}
\title{\bf Duality and Reciprocity for Hypergeometric \\ 
Series with Gamma Product Formula\thanks{ 
formula; closed-form expression; duality; reciprocity; contiguous 
relation; connection formula.}}  
\author{Katsunori Iwasaki\thanks{Department of Mathematics, 
Hokkaido University, Kita 10, Nishi 8, Kita-ku, Sapporo 060-0810 Japan. 
E-mail: {\tt iwasaki@math.sci.hokudai.ac.jp}}}
\date{February 9, 2018} 
\begin{document}
%%%%%%%%%%%%%%%%%%%%%%%%%%%%%%%%%%%%%%%%%%%%%%%%%%%%%%%%%%%%%%%%%%%%%%%%
\maketitle
\begin{abstract}
Following a previous article we continue our study on non-terminating  
hypergeometric series with one free parameter, which aims to find arithmetical 
constraints for a given hypergeometric series to admit a gamma product formula.   
In this article we exploit the concepts of duality and reciprocity not only to extend 
already obtained results to a larger region but also to strengthen themselves 
substantially.  
Among other things we are able to settle the rationality and 
finiteness conjectures posed in the previous article.    
\end{abstract}  
%%%%%%%%%%%%%%%%%%%%%%%%%%%%%%% sec:intro %%%%%%%%%%%%%%%%%%%%%%%%%%%%%%
\section{Introduction} \label{sec:intro}
%%%%%%%%%%%%%%%%%%%%%%%%%%%%%%%%%%%%%%%%%%%%%%%%%%%%%%%%%%%%%%%%%%%%%%%
Inspired by a series of works by Ebisu \cite{Ebisu1,Ebisu2,Ebisu3}, 
we set out to develop a theoretical study of non-terminating Gauss 
hypergeometric series with one free parameter in \cite{Iwasaki}. 
It aims to find arithmetical constraints for a given hypergeometric 
sum to admit a gamma product formula (GPF). 
We continue that study by discussing some of the problems arising there.   
In this article two symmetries which we call {\sl duality} and 
{\sl reciprocity} will play pivotal roles. 
We shall see how these symmetries can be used not only to extend our 
previous results to a larger region, but also to strengthen 
themselves substantially.  
In particular they are powerful enough to settle the 
{\sl rationality} and {\sl finite-cardinality} 
conjectures for the numbers $a$ and $b$ (see Theorem \ref{thm:ab-rf}).          
%%%%%
\par
%%%%% 
Given a data $\lambda = (p,q,r;a,b;x)$, we 
consider an entire meromorphic function       
%%%%%%%%%%%%%%%%%%%%%%%%%%%%%%% eqn:f %%%%%%%%%%%%%%%%%%%%%%%%%%%%%%%%%%
\begin{equation} \label{eqn:f}  
f(w;\lambda) := \hgF(p w+a, \, q w+b; \, r w; \, x),           
\end{equation} 
%%%%%%%%%%%%%%%%%%%%%%%%%%%%%%%%%%%%%%%%%%%%%%%%%%%%%%%%%%%%%%%%%%%%%%%%
where $\hgF(\alpha, \beta, \gamma; z)$ is the 
Gauss hypergeometric series. 
We are interested in the following.  
%%%%%%%%%%%%%%%%%%%%%%%%%%%%%%% prob:gpf %%%%%%%%%%%%%%%%%%%%%%%%%%%%%%%
\begin{problemgpf} \label{prob:gpf} 
Find a data $\lambda = (p,q,r;a,b;x)$ for which $f(w; \lambda)$ 
has a {\sl gamma product formula}: 
%%%%%%%%%%%%%%%%%%%%%%%%%%%%%%% eqn:gpf %%%%%%%%%%%%%%%%%%%%%%%%%%%
\begin{equation} \label{eqn:gpf}
f(w; \lambda) = C \cdot d^w \cdot 
\dfrac{\varGamma(w+u_1)\cdots\varGamma(w+u_m)}{\varGamma(w+v_1)\cdots\varGamma(w+v_n)},      
\end{equation}
%%%%%%%%%%%%%%%%%%%%%%%%%%%%%%%%%%%%%%%%%%%%%%%%%%%%%%%%%%%%%%%%%%%%%%%%
for some numbers $C, \, d \in \C^{\times}$; $m, \, n \in \Z_{\ge 0}$; 
$u_1, \dots, u_m \in \C$; and $v_1, \dots, v_n \in \C$.    
\end{problemgpf}
%%%%%%%%%%%%%%%%%%%%%%%%%%%%%%%%%%%%%%%%%%%%%%%%%%%%%%%%%%%%%%%%%%%%%%%%
\par
%%%%%% 
Problem $\rI$ is something like the peak of a high mountain, 
which is difficult to climb up directly.  
We need a base camp to attack it and this role is 
played by the following.  
%%%%%%%%%%%%%%%%%%%%%%%%%%%%%% prob:ocf %%%%%%%%%%%%%%%%%%%%%%%%%%%%%%%%
\begin{problemocf} \label{prob:ocf} 
Find a data $\lambda = (p,q,r;a,b;x)$ for which  
$f(w; \lambda)$ is {\sl of closed form}, that is,  
%%%%%%%%%%%%%%%%%%%%%%%%%%%%%% eqn:ocf %%%%%%%%%%%%%%%%%%%%%%%%%%%%%%%%%
\begin{equation} \label{eqn:ocf} 
\dfrac{f(w+1; \lambda)}{f(w; \lambda)} =: R(w; \lambda) 
\in \mathbb{C}(w) \,\,:\,\, \mbox{a rational function of $w$}.  
\end{equation}
%%%%%%%%%%%%%%%%%%%%%%%%%%%%%%%%%%%%%%%%%%%%%%%%%%%%%%%%%%%%%%%%%%%%%%%%
\end{problemocf}
%%%%%%%%%%%%%%%%%%%%%%%%%%%%%%%%%%%%%%%%%%%%%%%%%%%%%%%%%%%%%%%%%%%%%%%%
\par
%%%%% 
By the recursion formula for the gamma function 
$\varGamma(w+1) = w \, \varGamma(w)$, any solution to Problem $\rI$ 
is a solution to Problem $\rII$ with rational function  
%%%%%%%%%%%%%%%%%%%%%%%%%%%%% eqn:R %%%%%%%%%%%%%%%%%%%%%%%%%%%%%%%%%%%
\begin{equation} \label{eqn:R}
R(w; \lambda) = 
d \cdot \dfrac{(w+u_1) \cdots (w+u_m)}{(w+v_1) \cdots (w+v_n)}. 
\end{equation}  
%%%%%%%%%%%%%%%%%%%%%%%%%%%%%%%%%%%%%%%%%%%%%%%%%%%%%%%%%%%%%%%%%%%%%%%
%%%%%%%%%%%%%%%%%%%%%%%%%%%%% def:trichotomy %%%%%%%%%%%%%%%%%%%%%%%%%
\begin{definition} \label{def:trichtomy} 
The {\sl principal part} of a data $\lambda = (p,q,r;a,b;x)$ is 
the triple $\bp := (p,q,r)$. 
We say that $\lambda$ is {\sl integral} if $\bp \in \Z^3$; 
{\sl rational} if $\bp \in \Q^3$;  
and {\sl irrational} if $\bp \not\in \Q^3$.  
We can speak of an integral, rational or irrational 
solution to Problem $\rI$ or $\rII$. 
\end{definition}
%%%%%%%%%%%%%%%%%%%%%%%%%%%%%%%%%%%%%%%%%%%%%%%%%%%%%%%%%%%%%%%%%%%%%%%%%
\par
%%%%%
There is an efficient method to find {\sl integral} solutions to 
Problem $\rII$. 
A {\sl three-term relation} for $\hgF(\ba;z) := \hgF(\alpha,\beta;\gamma;z)$ 
with $\ba = (a_1,a_2;a_3) = (\alpha,\beta;\gamma)$ is an identity of 
the form: 
%%%%%%%%%%%%%%%%%%%%%%%%%%%%%%% eqn:3trF %%%%%%%%%%%%%%%%%%%%%%%%%%%%%%%
\begin{equation} \label{eqn:3trF}
\hgF(\ba+\bp;z) = r(\ba; z) \, \hgF(\ba;z) + q(\ba; z) \, \hgF(\ba+\1;z), 
\qquad \bp = (p,q;r) \in \mathbb{Z}^3, 
\end{equation}
%%%%%%%%%%%%%%%%%%%%%%%%%%%%%%%%%%%%%%%%%%%%%%%%%%%%%%%%%%%%%%%%%%%%%%%%
where $\1 := (1,1;1)$ and $q(\ba;z)$ and $r(\ba;z)$ are rational 
functions of $(\ba;z)$ depending uniquely on $\bp$.   
Relation \eqref{eqn:3trF} is obtained by composing a finite 
sequence of fifteen contiguous relations of Gauss.  
For an integral data $\lambda = (p,q,r;a,b;x)$, substituting 
$\ba = \bal(w) := (p w+a, q w +b; r w)$ and $z = x$ into equation  
\eqref{eqn:3trF}, we have a specialized three-term relation:   
%%%%%%%%%%%%%%%%%%%%%%%%%%%%%%% eqn:3trf %%%%%%%%%%%%%%%%%%%%%%%%%%%%%%
\begin{equation} \label{eqn:3trf}
f(w+1; \lambda) = R(w; \lambda) \, f(w; \lambda) + 
Q(w; \lambda) \, \tilde{f}(w; \lambda),  
\end{equation}
%%%%%%%%%%%%%%%%%%%%%%%%%%%%%%%%%%%%%%%%%%%%%%%%%%%%%%%%%%%%%%%%%%%%%%%
with $\tilde{f}(w; \lambda) := \hgF(\bal(w)+\1;x)$, where 
$Q(w;\lambda) := q(\bal(w);x)$ and $R(w;\lambda):= r(\bal(w);x)$ are 
rational functions of $w$ depending uniquely on $\lambda$. 
If $\lambda$ happens to be such a data that  
%%%%%%%%%%%%%%%%%%%%%%%%%%%%% eqn:cfcr %%%%%%%%%%%%%%%%%%%%%%%%%%%%%%%%
\begin{equation} \label{eqn:cfcr}
Q(w;\lambda) \equiv 0 \qquad \mbox{in} \quad \C(w), 
\end{equation}
%%%%%%%%%%%%%%%%%%%%%%%%%%%%%%%%%%%%%%%%%%%%%%%%%%%%%%%%%%%%%%%%%%%%%%% 
then three-term relation \eqref{eqn:3trf} reduces to a two-term one 
\eqref{eqn:ocf}, yielding a solution to Problem $\rII$. 
An integral solution to Problem $\rII$ so obtained is said to 
{\sl come from contiguous relations}.  
%%%%%%
\par
%%%%%%
Given a data $\lambda = (p,q,r; a,b; x)$ and a positive integer 
$k \in \N$, the new data $k \lambda := (k p, k q, kr; a,b; x)$ is 
referred to as the {\sl multiplication} of $\lambda$ by $k$.  
It is said to be {\it nontrivial} if $k \ge 2$, in particular  
$2 \lambda$ is the {\sl duplication} of $\lambda$.  
In view of definition \eqref{eqn:f} we have  
%%%%%%%%%%%%%%%%%%%%%%%%%%%%%%% eqn:f-mult %%%%%%%%%%%%%%%%%%%%%%%%%%%%%%
\begin{equation} \label{eqn:f-mult}
f(w; k \lambda) = f(k w; \lambda), \qquad k \in \N. 
\end{equation}
%%%%%%%%%%%%%%%%%%%%%%%%%%%%%%%%%%%%%%%%%%%%%%%%%%%%%%%%%%%%%%%%%%%%%%%%%%
Gauss's multiplication formula for the gamma function 
\cite[Theorem 1.5.2]{AAR}:  
%%%%%%%%%%%%%%%%%%%%%%%%%%%%%%% eqn:mult %%%%%%%%%%%%%%%%%%%%%%%%%%%%%%%%%
\begin{equation} \label{eqn:mult}
\varGamma(k w) = (2 \pi)^{(1-k)/2} \cdot k^{k w-1/2} \cdot 
\prod_{j=0}^{k-1} 
\varGamma\left(w + \textstyle \frac{j}{k}\right), \qquad k \in \N    
\end{equation}
%%%%%%%%%%%%%%%%%%%%%%%%%%%%%%%%%%%%%%%%%%%%%%%%%%%%%%%%%%%%%%%%%%%%%%%%%%
implies that if $\lambda$ is a solution to Problem $\rI$ with gamma 
product formula \eqref{eqn:gpf} then $k \lambda$ is also a solution to 
the same problem with the ``multiplied" gamma product formula 
%%%%%%%%%%%%%%%%%%%%%%%%%%%%% eqn:gpf-mult %%%%%%%%%%%%%%%%%%%%%%%%%%%%%%%
\begin{equation} \label{eqn:gpf-mult}
f(w; k \lambda) = C_k \cdot d_k^w \cdot 
\prod_{j=0}^{k-1} 
\dfrac{\varGamma\left(w+ \frac{u_1 + j}{k}\right) \cdots 
\varGamma\left(w+ \frac{u_m +j}{k}\right)}{\varGamma\left(w+ \frac{v_1 + j}{k}\right) 
\cdots \varGamma\left(w+ \frac{v_n +j}{k}\right)}, 
\end{equation}
%%%%%%%%%%%%%%%%%%%%%%%%%%%%%%%%%%%%%%%%%%%%%%%%%%%%%%%%%%%%%%%%%%%%%%%%%%
where $C_k := C \cdot (2 \pi)^{(k-1)(n-m)/2} \cdot 
k^{u-v+(n-m)/2}$ with $u := u_1+\cdots+u_m$, $v := v_1+\cdots+v_n$, and 
$d_k := d^k \cdot k^{k(m-n)}$.  
Similarly, if $\lambda$ is a solution to Problem $\rII$ with 
rational function $R(w; \lambda)$ in formula \eqref{eqn:R} then 
$k \lambda$ is also a solution to the same problem with  
%%%%%%%%%%%%%%%%%%%%%%%%%%%%% eqn:R-mult %%%%%%%%%%%%%%%%%%%%%%%%%%%%%%%%
\begin{equation} \label{eqn:R-mult}
R(w; k \lambda) = \prod_{j=0}^{k-1} R(k w + j; \lambda)
= d_k \cdot 
\prod_{j=0}^{k-1} 
\dfrac{\left(w+ \frac{u_1 + j}{k}\right) \cdots 
\left(w+ \frac{u_m +j}{k}\right)}{\left(w+ \frac{v_1 + j}{k}\right) 
\cdots \left(w+ \frac{v_n +j}{k}\right)}. 
\end{equation}
%%%%%%%%%%%%%%%%%%%%%%%%%%%%%%%%%%%%%%%%%%%%%%%%%%%%%%%%%%%%%%%%%%%%%%%%%%
\par 
%%%%% 
If $\lambda$ is a {\sl rational} solution to Problem $\rII$ and $k$ is the 
least common denominator of its principal part $\bp =(p, q; r) \in \Q^3$, 
then $k \lambda \in \Z^3$ is an integral solution to the same problem,  
so we can ask whether $k \lambda$ comes from contiguous relations.  
If the answer is ``yes", we say that the rational solution $\lambda$ 
{\sl essentially comes from contiguous relations}. 
It is easy to see that formula \eqref{eqn:gpf} can be recovered from 
formula \eqref{eqn:gpf-mult} via relation \eqref{eqn:f-mult}.   
Indeed, we have only to replace $w$ by $w/k$ in \eqref{eqn:gpf-mult} and 
use the multiplication formula \eqref{eqn:mult} in the other way round.   
However, it is totally unclear whether formula \eqref{eqn:R} can be 
recovered from formula \eqref{eqn:R-mult}, because $R(w; \lambda)$ does 
not have such a simple multiplicative structure as the function 
$f(w; \lambda)$ has in formula \eqref{eqn:f-mult}. 
%%%%%
\par
%%%%%
We now introduce two symmetries or transformations, which we call duality 
and reciprocity.     
%%%%%%%%%%%%%%%%%%%%%%%%%%%% def:duality %%%%%%%%%%%%%%%%%%%%%%%%%%%%%%%
\begin{definition} \label{def:duality} 
The {\sl duality} $\lambda = (p, q, r; a, b; x) \mapsto 
\lambda' := (p', q', r'; a', b'; x')$ is defined by   
%%%%%%%%%%%%%%%%%%%%%%%%%%% eqn:duality %%%%%%%%%%%%%%%%%%%%%%%%%%%%%%%%
\begin{equation} \label{eqn:duality}
p' := p, \quad q':= q, \quad r':=r; \quad 
a' := 1- \frac{2 p}{r} - a, \quad b' := 1 - \frac{2 q}{r} - b; 
\quad x':= x.        
\end{equation}
%%%%%%%%%%%%%%%%%%%%%%%%%%%%%%%%%%%%%%%%%%%%%%%%%%%%%%%%%%%%%%%%%%%%%%%%
It is a well-defined involution whenever $r$ does not vanish.    
\end{definition}
%%%%%%%%%%%%%%%%%%%%%%%%%%%%%%%%%%%%%%%%%%%%%%%%%%%%%%%%%%%%%%%%%%%%%%%%%
%%%%%%%%%%%%%%%%%%%%%%%% def:recip %%%%%%%%%%%%%%%%%%%%%%%%%%%%%%%%%%%%%%
\begin{definition} \label{def:recip} 
The {\sl reciprocity} $\lambda = (p, q, r; a, b; x) \mapsto 
\check{\lambda} := 
(\check{p}, \check{q}, \check{r}; \check{a}, \check{b}; \check{x})$ 
is defined by  
%%%%%%%%%%%%%%%%%%%%%%%%%% eqn:recip %%%%%%%%%%%%%%%%%%%%%%%%%%%%%%%%%%%%
\begin{equation} \label{eqn:recip}
\begin{split} 
\check{p} := -p, \qquad \check{q} :=-q, \qquad 
\check{r} :=r-p-q; \qquad \check{x} :=1-x,  \\[2mm] 
\check{a} := \dfrac{(r-q)(1-a)-p b}{r-p-q}, \qquad     
\check{b} := \dfrac{(r-p)(1-b)-q a}{r-p-q}.      
\end{split}  
\end{equation}
%%%%%%%%%%%%%%%%%%%%%%%%%%%%%%%%%%%%%%%%%%%%%%%%%%%%%%%%%%%%%%%%%%%%%%%%
This is a well-defined involution on the domain 
%%%%%%%%%%%%%%%%%%%%%%%%%% eqn:real1 %%%%%%%%%%%%%%%%%%%%%%%%%%%%%%%%%%%
\begin{equation} \label{eqn:real1}
p, \, q, \, r \in \R, \quad r(r-p-q) \neq 0; \qquad a, \, b \in \R; 
\qquad 0 < x < 1.     
\end{equation}
%%%%%%%%%%%%%%%%%%%%%%%%%%%%%%%%%%%%%%%%%%%%%%%%%%%%%%%%%%%%%%%%%%%%%%%% 
\end{definition}  
%%%%%%%%%%%%%%%%%%%%%%%%%%%%%%%%%%%%%%%%%%%%%%%%%%%%%%%%%%%%%%%%%%%%%%%%
\par
%%%%%%%
The origin of these transformations is very simple.   
The Gauss hypergeometric equation 
%%%%%%%%%%%%%%%%%%%%%%%%%% eqn:hge %%%%%%%%%%%%%%%%%%%%%%%%%%%%%%%%%%%%%
\begin{equation} \label{eqn:hge}
z(1-z) \dfrac{d u}{d z^2} + \{\gamma - (\alpha+\beta+1) z \} 
\dfrac{d u}{d z} - \alpha \beta \, u = 0 
\end{equation} 
%%%%%%%%%%%%%%%%%%%%%%%%%%%%%%%%%%%%%%%%%%%%%%%%%%%%%%%%%%%%%%%%%%%%%%%%
is a Fuchsian differential equation with the Riemann scheme: 
%%%%%%%%%%%%%%%%%%%%%%%%% eqn:rs %%%%%%%%%%%%%%%%%%%%%%%%%%%%%%%%%%%%%%%
\begin{equation} \label{eqn:rs}
\left\{\begin{matrix}
z = 0    & z = 1               & z = \infty \\[1mm]
0        & 0                   & \alpha     \\[1mm]
1-\gamma & \gamma-\alpha-\beta & \beta 
\end{matrix}\right\}, 
\end{equation}
%%%%%%%%%%%%%%%%%%%%%%%%%%%%%%%%%%%%%%%%%%%%%%%%%%%%%%%%%%%%%%%%%%%%%%%
in which the hypergeometric series $\hgF(\alpha, \beta; \gamma; z)$ is 
just the solution of local exponent $0$ at the origin $z = 0$. 
Exchanging it with the solution of exponent $1-\gamma$ at the same 
point $z = 0$ yields duality, while exchanging it with the solution 
of exponent $\gamma-\alpha-\beta$ at $z = 1$ gives reciprocity; 
we refer to \S\ref{sec:kummer} for their constructions.  
%%%%%
\par
\begin{figure}[t]
\begin{minipage}{0.45\hsize}
\begin{center}
\includegraphics[width=65mm,clip]{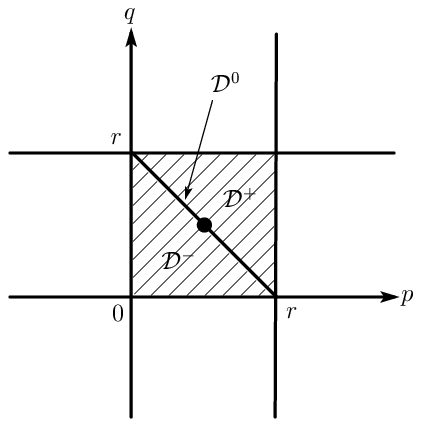}
\end{center} 
\caption{Central square $\cD$.}
\label{fig:D}
\end{minipage}
\begin{minipage}{0.45\hsize}
\begin{center}
\includegraphics[width=65mm,clip]{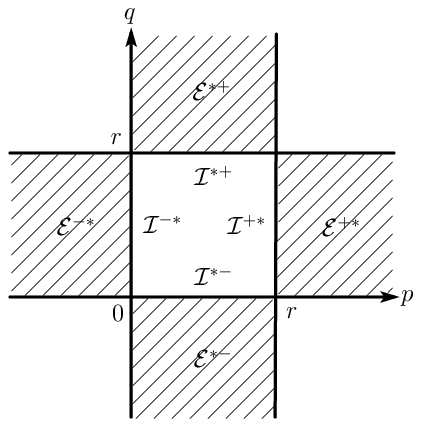}  
\end{center}
\caption{Wings $\cE$ and borders $\cI$.}
\label{fig:E}
\end{minipage} 
\end{figure}

%%%%%%
Suppose that the data $\lambda = (p,q,r;a,b;x)$ lies in the real domain 
%%%%%%%%%%%%%%%%%%%%%%%%%%% eqn:real2 %%%%%%%%%%%%%%%%%%%%%%%%%%%%%%%%%%
\begin{equation} \label{eqn:real2}
p, \, q \in \R; \qquad r > 0; \qquad a, \, b \in \R; \qquad 0 < x < 1.   
\end{equation} 
%%%%%%%%%%%%%%%%%%%%%%%%%%%%%%%%%%%%%%%%%%%%%%%%%%%%%%%%%%%%%%%%%%%%%%%%%
In \cite[\S1]{Iwasaki} we introduced its subregions 
$\cD = \cD^- \cup \cD^0 \cup \cD^+$ and    
$\cE = \cE^{*-} \cup \cE^{*+} \cup \cE^{-*} \cup \cE^{+*}$, the figures 
of which are depicted in Figures \ref{fig:D} and \ref{fig:E} upon 
projected down to the $(p, q)$-plane with a fixed value of $r > 0$.    
%%%%%%%%%%%%%%%%%%%%% fig:D %%%%% fig:E %%%%%%%%%%%%%%%%%%%%%%%%%%%%%%%%%
%%%%%%%%%%%%%%%%%%%%%%%%%%%%%%%%%%%%%%%%%%%%%%%%%%%%%%%%%%%%%%%%%%%%%%%%%
Note that the study on $\cD^0$ is finished by assertion (1) of 
\cite[Theorem 2.2]{Iwasaki}. 
We now introduce a new subregion $\cF = \cF^- \cup \cF^+$ with 
two components     
%%%%%%%
\begin{align*}  
\cF^- &:= \{\, \lambda \in \R^6 \,:\,\, 
p < 0, \quad q < 0, \quad r > 0; 
\quad a, \, b \in \R; \quad 0 < x < 1 \, \}, 
\\
\cF^+ &:= \{\, \lambda \in \R^6 \,:\,\, 
p > r, \quad q > r, \quad r > 0; 
\quad a, \, b \in \R; \quad 0 < x < 1 \, \}.          
\end{align*}
%%%%%%
Thanks to the classical symmetries in \cite[(8)]{Iwasaki} we have only 
to work on $\cD^-$, $\cE^{*-}$, $\cF^-$, because  
%%%%%%%%%%%%%%%%
\begin{alignat*}{5}
\cD^{-}  &\leftrightarrow \cD^{+},  \quad  &  
\cE^{*-} &\leftrightarrow \cE^{*+}, \quad  & 
\cE^{-*} &\leftrightarrow \cE^{+*},  \quad &
\cF^{-}  &\leftrightarrow \cF^{+}   \qquad & 
         & \mbox{by \cite[(8b)]{Iwasaki}}, \\    
         &                                 &  
\cE^{*-} &\leftrightarrow \cE^{-*}, \quad  & 
\cE^{*+} &\leftrightarrow \cE^{+*}   \quad &
         &                                 &   
         & \mbox{by \cite[(8a)]{Iwasaki}}. 
\end{alignat*}
%%%%%%%%%%%%%%%    
Duality maps each of $\cD^-$, $\cE^{*-}$ and $\cF^-$ 
bijectively onto itself, while reciprocity induces 
transpositions $\cD^{-} \leftrightarrow \cF^-$ and 
$\cE^{*-} \leftrightarrow \cE^{-*}$. 
A chief idea underlying these symmetries is the concept of  
Ebisu symmetries in \S \ref{sec:kummer}, based on 
which duality and reciprocity are constructed in 
\S \ref{sec:duality} and \S \ref{sec:r-cf}.   
%%%%%%%%%%%%%%%%%%%%%%%%%%%%%% fig:D-r %%%%%%%%%%%%%%%%%%%%%%%%%%%%%%
\begin{figure}[t]
\begin{center}
\includegraphics[width=160mm,clip]{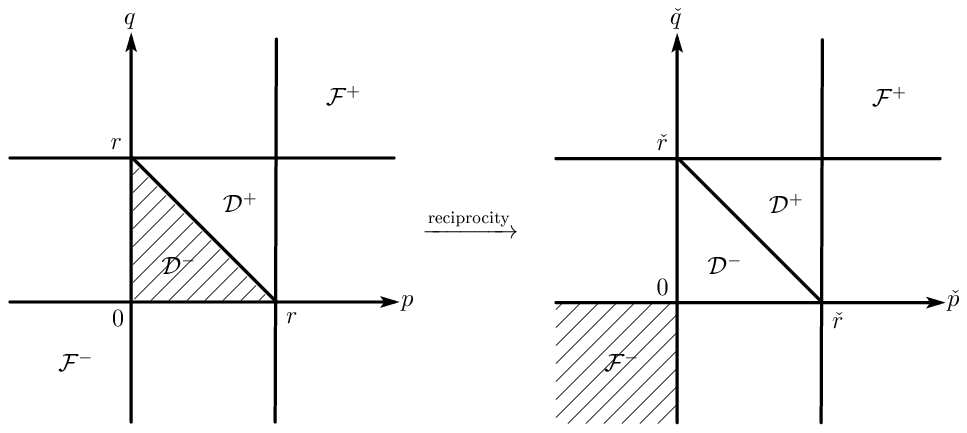}   
\end{center}
\caption{Reciprocity between $\cD^-$ and $\cF^-$.}
\label{fig:D-r} 
\end{figure}
%%%%%%%%%%%%%%%%%%%%%%%%%%%%%%%%%%%%%%%%%%%%%%%%%%%%%%%%%%%%%%%%%%%%%%
%%%%%%%%%%%%%%%%%%%%%%%%%%%%%% fig:E-r %%%%%%%%%%%%%%%%%%%%%%%%%%%%%%%
\begin{figure}[t]
\begin{center}
\includegraphics[width=160mm,clip]{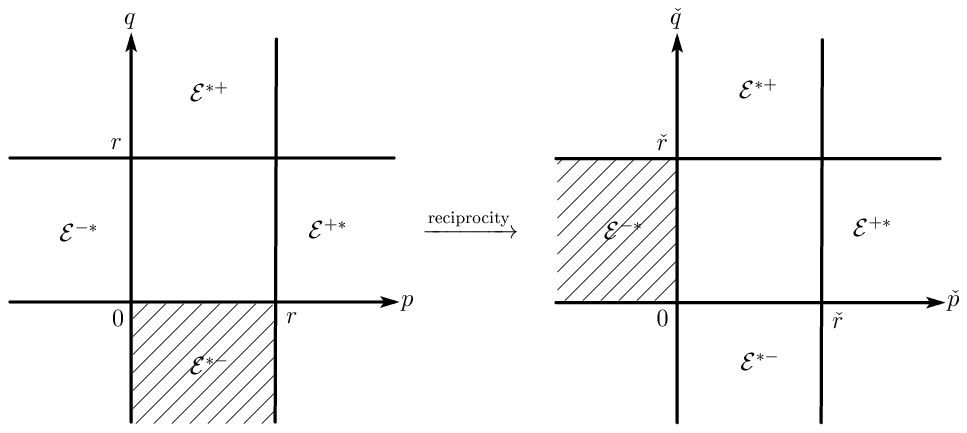}   
\end{center}
\caption{Reciprocity between $\cE^{*-}$ and $\cE^{-*}$.}  
\label{fig:E-r} 
\end{figure}
%%%%%%%%%%%%%%%%%%%%%%%%%%%%%%%%%%%%%%%%%%%%%%%%%%%%%%%%%%%%%%%%%%%%%%%%%         
%%%%%%%%%%%%%%%%%%%%%%%%%%%%% sec:results %%%%%%%%%%%%%%%%%%%%%%%%%%%%%%
\section{Main Results} \label{sec:results} 
%%%%%%%%%%%%%%%%%%%%%%%%%%%%%%%%%%%%%%%%%%%%%%%%%%%%%%%%%%%%%%%%%%%%%%%%
By \cite[Theorem 2.1]{Iwasaki} Problems $\rI$ and $\rII$ are 
equivalent in $\cD \cup \cE$, thus in this region we can speak of 
a solution without specifying to which problem it is a solution.   
By \cite[Theorem 2.2]{Iwasaki} any solution 
$\lambda = (p,q,r;a,b;x) \in \cD^-$ must satisfy either  
\begin{center}
$(\rA)$ \quad $p$, $q$, $r \in \Z$, \, $r-p-q \equiv 0 \mod 2$; 
\qquad or \qquad   
$(\rB)$ \quad $p$, $q \in \frac{1}{2} + \Z$, \, $r \in \Z$.  
\end{center}  
%%%%%%%
\par
%%%%%%  
By \cite[Theorem 2.3]{Iwasaki} any $(\rA)$-solution 
$\lambda = (p,q,r;a,b;x) \in \cD^-$ admits a GPF of the form  
%%%%%%%%%%%%%%%%%%%%%%%%% eqn:gpf2 %%%%%%%%%%%%%%%%%%%%%%%%%%%%%%%%%%
\begin{equation} \label{eqn:gpf2} 
f(w; \lambda) = C \cdot d^w \cdot 
\dfrac{\prod_{i=0}^{r-1}\varGamma\left(w+\frac{i}{r} \right)}{\prod_{i=1}^r 
\varGamma\left(w+v_i \right)},  
\end{equation} 
%%%%%%%%%%%%%%%%%%%%%%%%%%%%%%%%%%%%%%%%%%%%%%%%%%%%%%%%%%%%%%%%%%%%%
where $C$ is a positive constant and $d$ is given by 
%%%%%%%%%%%%%%%%%%%%%%%%%% eqn:d %%%%%%%%%%%%%%%%%%%%%%%%%%%%%%%%%%%%
\begin{equation} \label{eqn:d}
d = \dfrac{r^r}{\sqrt{p^p q^q (r-p)^{r-p} (r-q)^{r-q} x^r (1-x)^{p+q-r}}},    
\end{equation} 
%%%%%%%%%%%%%%%%%%%%%%%%%%%%%%%%%%%%%%%%%%%%%%%%%%%%%%%%%%%%%%%%%%%%%
while $v_1, \dots, v_r$ are such numbers that sum up to
%%%%%%%%%%%%%%%%%%%%%%%%% eqn:sum %%%%%%%%%%%%%%%%%%%%%%%%%%%%%%%%%%%
\begin{equation} \label{eqn:sum} 
v_1 + \cdots + v_r = \dfrac{r-1}{2}, 
\end{equation}
%%%%%%%%%%%%%%%%%%%%%%%%%%%%%%%%%%%%%%%%%%%%%%%%%%%%%%%%%%%%%%%%%%%%%
and that admit a division relation in $\C[w]$,     
%%%%%%%%%%%%%%%%%%%%%%%%% eqn:division %%%%%%%%%%%%%%%%%%%%%%%%%%%%%%
\begin{equation} \label{eqn:division}
\prod_{i=1}^r \left(w+v_i\right) \,\, \Big| \,\,    
\prod_{i=1}^{p-1}\left(w+{\ts\frac{i+a}{p}}\right) 
\prod_{i=1}^{q-1}\left(w+{\ts\frac{i+b}{q}}\right) 
\prod_{j=0}^{r-p-1}\left(w+{\ts\frac{j-a}{r-p}}\right) 
\prod_{j=0}^{r-q-1}\left(w+{\ts\frac{j-b}{r-q}}\right).      
\end{equation} 
%%%%%%%%%%%%%%%%%%%%%%%%%%%%%%%%%%%%%%%%%%%%%%%%%%%%%%%%%%%%%%%%%%%%%%
This allows us to introduce the numbers $v_1^*, \dots, v_r^*$ 
complementary to $v_1, \dots, v_r$ by       
%%%%%%%%%%%%%%%%%%%%%%%%%%%%% eqn:vi-s %%%%%%%%%%%%%%%%%%%%%%%%%%%%%%%
\begin{equation} \label{eqn:vi-s} 
\resizebox{0.9\hsize}{!}{$\displaystyle
\prod_{i=1}^r \left(w+v_i\right) \prod_{i=1}^r \left(w+v_i^*\right) 
= 
\prod_{i=0}^{p-1}\left(w+{\ts\frac{i+a}{p}}\right) 
\prod_{i=0}^{q-1}\left(w+{\ts\frac{i+b}{q}}\right) 
\prod_{j=0}^{r-p-1}\left(w+{\ts\frac{j-a}{r-p}}\right) 
\prod_{j=0}^{r-q-1}\left(w+{\ts\frac{j-b}{r-q}}\right). $}       
\end{equation}
%%%%%%%%%%%%%%%%%%%%%%%%%%%%%%%%%%%%%%%%%%%%%%%%%%%%%%%%%%%%%%%%%%%%%%%%%
\par
%%%%%
The first result is regarding a description of duality for 
$(\rA)$-solutions in $\cD^-$. 
%%%%%%%%%%%%%%%%%%%%%%%%%% thm:gpf-d %%%%%%%%%%%%%%%%%%%%%%%%%%%%%%%%%%%%
\begin{theorem}\label{thm:gpf-d} 
The duality $\lambda \mapsto \lambda'$ induces an involution on the set 
of all $(\rA)$-solutions in $\cD^-$.   
For those solutions it yields a transformation of gamma product formulas    
%%%%%%%%%%%%%%%%%%%%%%%%%% eqn:gpf-d %%%%%%%%%%%%%%%%%%%%%%%%%%%%%%%%%%%%
\begin{equation} \label{eqn:gpf-d}
f(w; \lambda) = C \cdot d^w \cdot 
\dfrac{\prod_{i=0}^{r-1} \varGamma\left(w+\frac{i}{r}\right)}{\prod_{i=1}^{r} \varGamma\left(w+v_i\right)} 
\quad \mapsto \quad 
f(w; \lambda') = C' \cdot d^w \cdot 
\dfrac{\prod_{i=0}^{r-1} \varGamma\left(w+\frac{i}{r}\right)}{\prod_{i=1}^{r} \varGamma\left(w+v_i'\right)},     
\end{equation}
%%%%%%%%%%%%%%%%%%%%%%%%%%%%%%%%%%%%%%%%%%%%%%%%%%%%%%%%%%%%%%%%%%%%%%%%
where $C'$ is a positive constant, $d$ is the number given by 
formula \eqref{eqn:d}, while 
%%%%%%%%%%%%%%%%%%%%%%%%%%%%%% eqn:vi-p %%%%%%%%%%%%%%%%%%%%%%%%%%%%%
\begin{equation} \label{eqn:vi-p}
v_i':= 1 - \dfrac{2}{r} - v_i^* \qquad (i = 1, \dots, r), 
\end{equation}
%%%%%%%%%%%%%%%%%%%%%%%%%%%%%%%%%%%%%%%%%%%%%%%%%%%%%%%%%%%%%%%%%%%%%%%%%
with $v_1^*, \dots, v_r^*$ being defined by equation \eqref{eqn:vi-s}.  
\end{theorem}
%%%%%%%%%%%%%%%%%%%%%%%%%%%%%%%%%%%%%%%%%%%%%%%%%%%%%%%%%%%%%%%%%%%%%%%%%
\par
%%%%%%
The second result is about the reciprocity $\lambda \mapsto \check{\lambda}$ 
of $(\rA)$-solutions $\lambda \in \cD^-$ into $\check{\lambda} \in \cF^-$. 
%%%%%%%%%%%%%%%%%%%%%%%%%% thm:D-r %%%%%%%%%%%%%%%%%%%%%%%%%%%%%%%%%%%%%%%
\begin{theorem} \label{thm:D-r}  
For any $(\rA)$-solution $\lambda = (p,q,r;a,b;x) \in \cD^-$  
there exists a division relation
%%%%%%%%%%%%%%%%%%%%%%%%%%%% eqn:division2 %%%%%%%%%%%%%%%%%%%%%%%%%% 
\begin{equation} \label{eqn:division2} 
\prod_{i=0}^{p-1}\left(w+{\ts\frac{i+a}{p}}\right) 
\prod_{i=0}^{q-1}\left(w+{\ts\frac{i+b}{q}}\right) 
\, \Big| \, \prod_{i=1}^r \left(w+v_i\right) 
\qquad \mbox{in \, $\C[w]$},    
\end{equation} 
%%%%%%%%%%%%%%%%%%%%%%%%%%%%%%%%%%%%%%%%%%%%%%%%%%%%%%%%%%%%%%%%%%%%%
which allows us to rearrange the numbers $v_1, \dots,v_r$ in    
\eqref{eqn:gpf2} so that 
%%%%%%%%%%%%%%%%%%%%%%%%%%% eqn:vi-pq %%%%%%%%%%%%%%%%%%%%%%%%%%%%%%%
\begin{equation} \label{eqn:vi-pq}
\prod_{i=0}^{p-1}\left(w+{\ts\frac{i+a}{p}}\right) 
\prod_{i=0}^{q-1}\left(w+{\ts\frac{i+b}{q}}\right) 
= \prod_{i=r-p-q+1}^{r} \left(w+v_i\right).  
\end{equation}
%%%%%%%%%%%%%%%%%%%%%%%%%%%%%%%%%%%%%%%%%%%%%%%%%%%%%%%%%%%%%%%%%%%%%%%
With this convention the reciprocity $\lambda \mapsto \check{\lambda}$ 
takes any $(\rA)$-solution $\lambda \in \cD^-$ to a solution 
$\check{\lambda} \in \cF^-$ to Problem $\rI$, inducing a  
transformation of gamma product formulas:   
%%%%%%%%%%%%%%%%%%%%%%%%%%%% eqn:gpf-c %%%%%%%%%%%%%%%%%%%%%%%%%%%%%%%%
\begin{equation} \label{eqn:gpf-c}
\resizebox{0.8\hsize}{!}{$\displaystyle
f(w; \lambda) = C \cdot d^w \cdot 
\dfrac{\prod_{i=0}^{r-1}\varGamma\left(w+\frac{i}{r} \right)}{\prod_{i=1}^r 
\varGamma\left(w+v_i \right)} \,\, \mapsto \,\,  
f(w; \check{\lambda}) = \check{C} \cdot \check{d}^w \cdot 
\dfrac{\prod_{i=0}^{r-p-q-1}\varGamma\left(w+\frac{i}{r-p-q} \right)}{
\prod_{i=1}^{r-p-q} \varGamma\left(w+ \check{v}_i \right)}, $}  
\end{equation}
%%%%%%%%%%%%%%%%%%%%%%%%%%%%%%%%%%%%%%%%%%%%%%%%%%%%%%%%%%%%%%%%%%%%%%%%%%%%%
where $\check{C}$ is a positive constant and 
%%%%%%%%%%%%%%%%%%%%%%%%% eqn:d-c %%%% eqn:vi-c %%%%%%%%%%%%%%%%%%%%%%%%%%%%%
\begin{align}
\check{d} &= (r-p-q)^{r-p-q} 
\sqrt{\frac{p^p \, q^q \, x^r}{(r-p)^{r-p} \, (r-q)^{r-q} \, (1-x)^{r-p-q}}}, 
\label{eqn:d-c} 
\\[3mm]
\check{v}_i &= v_i - \frac{1-a-b}{r-p-q} \qquad (i = 1, \dots, r-p-q). 
\label{eqn:vi-c} 
\end{align}
%%%%%%%%%%%%%%%%%%%%%%%%%%%%%%%%%%%%%%%%%%%%%%%%%%%%%%%%%%%%%%%%%%%%%%%%%%%%%
\end{theorem}
%%%%%%%%%%%%%%%%%%%%%%%%%%%%%%%%%%%%%%%%%%%%%%%%%%%%%%%%%%%%%%%%%%%%%%%%%%%%%
\par
%%%%%
For any $(\rA)$-solution in $\cD^{-}$ its principal part $\bp = (p,q,r)$ must belong 
to  
%%%%%%%%%%%%%%%%%%%%%%%%%%%%% eqn:D-m %%%%%%%%%%%%%%%%%%%%%%%%%%%%%%%%%%
\begin{equation} \label{eqn:D-m}
D^{-}_{\rA} := \{\, \bp = (p,q;r) \in \Z^3 \,:\, p > 0, \,\, q > 0, \,\,  
0 < r-p-q \equiv 0 \mod 2 \, \}.       
\end{equation} 
%%%%%%%%%%%%%%%%%%%%%%%%%%%%%%%%%%%%%%%%%%%%%%%%%%%%%%%%%%%%%%%%%%%%%%
The existence of duality and reciprocity puts further arithmetic constraints on 
the solutions.  
%%%%%%%%%%%%%%%%%%%%%%%%%%% thm:ab-rf %%%%%%%%%%%%%%%%%%%%%%%%%%%%%%%%%%%%%%
\begin{theorem} \label{thm:ab-rf} 
For any $(\rA)$-solution $\lambda = (p,q,r;a,b;x) \in \cD^-$ we must have 
%%%%%%%%%%%%%%%%%%% eqn:dr %%%% eqn:estimate %%%%%%%%%%%%%%%%%%%%%%%%%%%%%%%
\begin{gather}
(\, p|r \quad \mbox{or} \quad p|(r-p-q) \,) \quad \mbox{and} \quad 
(\, q|r \quad \mbox{or} \quad q|(r-p-q) \,),  \label{eqn:dr} \\[2mm] 
a, b, \, v_1, \dots, v_r \in \Q \cap [0, \, 1), \qquad 
x \in \overline{\Q} \cap (0, \, 1), \label{eqn:estimate} 
\end{gather} 
%%%%%%%%%%%%%%%%%%%%%%%%%%%%%%%%%%%%%%%%%%%%%%%%%%%%%%%%%%%%%%%%%%%%%%%%%%  
where $v_1, \dots, v_r$ are the numbers in GPF \eqref{eqn:gpf2}.   
Given any integer triple $\bp = (p,q;r)$, there are no or only a finite 
number of $(\rA)$-solutions $\lambda = (p,q,r;a,b;x) \in \cD^-$ with 
prescribed principal part $\bp$ and there is an algorithm to determine 
all of them in finite steps $($see Algorithm $\ref{algorithm})$.     
\end{theorem} 
%%%%%%%%%%%%%%%%%%%%%%%%%%%%%%%%%%%%%%%%%%%%%%%%%%%%%%%%%%%%%%%%%%%%%%%%
\par
%%%%%%
Since the duplication of a $(\rB)$-solution is an $(\rA)$-solution, 
Theorem \ref{thm:ab-rf} leads to the following.      
%%%%%%%%%%%%%%%%%%%%%%%%%%%%% thm:sol-B %%%%%%%%%%%%%%%%%%%%%%%%%%%%%%%%%%%%%
\begin{corollary} \label{thm:sol-B} 
If $\lambda = (p,q,r;a,b;x) \in \cD^-$ is a $(\rB)$-solution then $a$ and 
$b$ are rational numbers with $0 \le a, \, b < 1$, $x$ is an algebraic number, 
and $\lambda$ admits a GPF of the form \eqref{eqn:gpf2} where $d$ is given by 
formula \eqref{eqn:d} and $v_1, \dots, v_r$ are rational numbers such that  
$0 \le v_1, \dots, v_r < 1$.   
\end{corollary}
%%%%%%%%%%%%%%%%%%%%%%%%%%%%%%%%%%%%%%%%%%%%%%%%%%%%%%%%%%%%%%%%%%%%%%%%%%%%%
\par
%%%%%%
Theorem \ref{thm:D-r} concerns the reciprocity in the direction 
$\cD^- \to \cF^-$. 
Starting from $\cF^-$ take the reciprocity in the opposite direction  
$\cF^- \to \cD^-$ and then use Theorem \ref{thm:D-r} to get back in 
$\cF^-$.  
Then we have the following.         
%%%%%%%%%%%%%%%%%%%%%%%%%% thm:F %%%%%%%%%%%%%%%%%%%%%%%%%%%%%%%%%%%%%%%%%%%%
\begin{theorem} \label{thm:F}  
Let $\lambda = (p, q, r; a, b; x) \in \cF^-$ in what follows.  
\begin{enumerate} 
\item Any solution $\lambda \in \cF^-$ to Problem $\rI$ or $\rII$ 
is non-elementary, that is, $f(w; \lambda)$ has infinitely many poles in $\C_w$.     
\item If $\lambda \in \cF^-$ is an integral solution to Problem $\rII$, 
then $\lambda$ comes from contiguous relations with its reciprocal 
$\check{\lambda}$ being an $(\rA)$-solution in $\cD^-$,   
in other words, $\lambda$ is the reciprocal of an $(\rA)$-solution 
$\check{\lambda} \in \cD^-$ so that $r$ must be an even positive 
integer and $\lambda$ becomes a solution to Problem $\rI$ with  
$f(w; \lambda)$ having a GPF as in formula \eqref{eqn:gpf-F} below.  
\item If $\lambda \in \cF^-$ is a rational solution to Problem $\rII$, 
then $\lambda$ essentially comes from contiguous relations, 
$r$ must be a positive integer but not necessarily even; $a$, $b \in \Q$;   
$x$ algebraic, and $\lambda$ becomes a solution to Problem $\rI$ with 
$f(w; \lambda)$ having a GPF of the form  
%%%%%%%%%%%%%%%%%%%%%%%%%%%%%% eqn:gpf-F %%%%%%%%%%%%%%%%%%%%%%%%%%%%%%%%%%%%%
\begin{equation} \label{eqn:gpf-F} 
f(w; \lambda) = C \cdot d^w \cdot 
\dfrac{\prod_{i=0}^{r-1} \varGamma \left( w+\frac{i}{r}\right) }{\prod_{i=1}^{r} 
\varGamma\left(w+v_i\right)}, 
\end{equation}
%%%%%%%%%%%%%%%%%%%%%%%%%%%%%%%%%%%%%%%%%%%%%%%%%%%%%%%%%%%%%%%%%%%%%%%%%%%%% 
where $C$ is a positive constant, $d$ is a positive algebraic number defined 
by  
%%%%%%%%%%%%%%%%%%%%%%%%%%%%%% eqn:d-F %%%%%%%%%%%%%%%%%%%%%%%%%%%%%%%%%%%%%%%
\begin{equation} \label{eqn:d-F}
d = r^r \sqrt{\dfrac{|p|^{|p|} \, |q|^{|q|} \, 
(1-x)^{r-p-q}}{(r-p)^{r-p} \, (r-q)^{r-q} \, x^r}}, 
\end{equation}
%%%%%%%%%%%%%%%%%%%%%%%%%%%%%%%%%%%%%%%%%%%%%%%%%%%%%%%%%%%%%%%%%%%%%%%%%%%%%  
and $v_1, \dots, v_r$ are such numbers that satisfy the following conditions: 
%%%%%%%%%%%%%%%%%%%%%%%%%% eqn:vi-F %%%%%%%%%%%%%%%%%%%%%%%%%%%%%%%%%%%%%%%%% 
\begin{equation} \label{eqn:vi-F} 
v_1 + \cdots + v_r = \frac{r-1}{2}, \quad  
v_1, \dots, v_r \in \left(\Q \setminus \textstyle \frac{1}{r} \, \Z\right) 
\cap \left[c, \, c+1\right), \quad c:= \frac{1-a-b}{r-p-q}.  
\end{equation}
%%%%%%%%%%%%%%%%%%%%%%%%%%%%%%%%%%%%%%%%%%%%%%%%%%%%%%%%%%%%%%%%%%%%%%%%%%%%%
In particular, Problems $\rI$ and $\rII$ are equivalent for 
rational data in $\cF^-$. 
\end{enumerate}
\end{theorem}
%%%%%%%%%%%%%%%%%%%%%%%%%%%%%%%%%%%%%%%%%%%%%%%%%%%%%%%%%%%%%%%%%%%%%%%%%%%%%
%%%%%%%%%%%%%%%%%%%%%%%%%% rem:F %%%%%%%%%%%%%%%%%%%%%%%%%%%%%%%%%%%%%%%%%%%%
\begin{remark} \label{rem:F} 
We make three comments about Theorem \ref{thm:F}. 
\begin{enumerate} 
\item When $\lambda$ is integral, formula \eqref{eqn:gpf-F} is computable 
from the GPF for the $(\rA)$-solution $\check{\lambda} \in \cD^-$ through  
transformation \eqref{eqn:gpf-c}, where the roles of $\lambda$ and 
$\check{\lambda}$ are exchanged.  
Regarding duality it is possible to formulate a result similar to 
Theorem \ref{thm:gpf-d} for integral solutions in $\cF^-$.      
\item We mention how to calculate formula \eqref{eqn:gpf-F} when 
$\lambda$ is rational but not integral.    
If $k$ is the least common denominator of $p$, $q \in \Q$,  
then $k \lambda \in \cF^-$ is an integral solution to Problem $\rII$, 
so the GPF for $k \lambda$ is computable from the GPF for the 
$(\rA)$-solution $(k \lambda)^{\vee} \in \cD^-$ via 
transformation \eqref{eqn:gpf-c}.      
It turns out that the ensuing result is       
%%%%%%%%%%%%%%%%%%%%%%%%%%%%% eqn:gpf-F-k %%%%%%%%%%%%%%%%%%%%%%%%%%%%%%%%%%%
\begin{equation} \label{eqn:gpf-F-k} 
f(w; k \lambda) = C_k \cdot d^{k w} \cdot 
\dfrac{\prod_{i=0}^{k r-1} 
\varGamma\left(w+\frac{i}{k r}\right)}{\prod_{i=1}^{r} \prod_{j=0}^{k-1} 
\varGamma\left(w+\frac{v_i + j}{k}\right)},   
\end{equation}
%%%%%%%%%%%%%%%%%%%%%%%%%%%%%%%%%%%%%%%%%%%%%%%%%%%%%%%%%%%%%%%%%%%%%%%%%%%%%%
for some $C_k > 0$ and $v_1, \dots, v_r$ satisfying condition \eqref{eqn:vi-F}.    
Replacing $w$ with $w/k$ in \eqref{eqn:gpf-F-k} and using the 
multiplication formula \eqref{eqn:mult} for the gamma function we get the 
desired formula \eqref{eqn:gpf-F}. 
For details we refer to the proofs of Lemma \ref{lem:gpf-F} and 
Proposition \ref{prop:gpf-F}.    
\item Assertions (2) and (3) of Theorem \ref{thm:F} address only integral 
or rational solutions in $\cF^-$.     
It is an interesting open problem to know whether $\cF^-$ contains 
any irrational solution.   
\end{enumerate} 
\end{remark}
%%%%%%%%%%%%%%%%%%%%%%%%%%%%%%%%%%%%%%%%%%%%%%%%%%%%%%%%%%%%%%%%%%%%%%%%%%%%%
\par
%%%%%% 
Assertions (1) and (2) of Theorem \ref{thm:F} will be proved in 
\S \ref{sec:SW} as Propositions \ref{prop:cfcr-F}, while assertion (3) 
will be established as Proposition \ref{prop:gpf-F} by developing the idea 
in item (2) of Remark \ref{rem:F}.  
%%%%%%
\par
%%%%% 
Finally we present a small but important result on the domain 
$\cE^{*-} \cup \cE^{-*}$. 
%%%%%%%%%%%%%%%%%%%%%%%%%% thm:E %%%%%%%%%%%%%%%%%%%%%%%%%%%%%%%%%%%%%%%%
\begin{theorem} \label{thm:E}
Duality $\lambda \mapsto \lambda'$ maps all integral solutions in 
$\cE^{*-}$ onto themselves, while reciprocity 
$\lambda \mapsto \check{\lambda}$ maps all integral solutions in 
$\cE^{*-}$ onto those in $\cE^{-*}$, both bijectively.    
\end{theorem} 
%%%%%%%%%%%%%%%%%%%%%%%%%%%%%%%%%%%%%%%%%%%%%%%%%%%%%%%%%%%%%%%%%%%%%%%%%
\par
%%%%%%
The former and latter assertions of this theorem will be proved in 
\S \ref{sec:duality} and \S \ref{sec:r-cf} as Corollaries 
\ref{cor:duality1} and \ref{cor:r-cf1} respectively. 
Applications of the theorem will be discussed elsewhere.  
%%%%%%%%%%%%%%%%%%%%%%%%%% sec:kummner %%%%%%%%%%%%%%%%%%%%%%%%%%%%%%%%%
\section{Kummer's 24 Solutions and Ebisu Symmetries} \label{sec:kummer} 
%%%%%%%%%%%%%%%%%%%%%%%%%%%%%%%%%%%%%%%%%%%%%%%%%%%%%%%%%%%%%%%%%%%%%%%%
Kummer \cite{Kummer} constructed twenty-four solutions (or more precisely 
power series representations of solutions) to the Gauss hypergeometric 
equation \eqref{eqn:hge}.  
They are known as {\sl Kummer's twenty-four solutions}, among which 
the hypergeometric series $\hgF(\ba;z)$ is the most 
representative member.   
A complete list of them can be found in Erd\'elyi 
\cite[Chap. I\!I, \S 2.9, formulas (1)--(24)]{Erdelyi}.  
Ebisu \cite[Lemma 2.2]{Ebisu3} showed that each of Kummer's 
solutions, say $\hgK(\ba;z)$, admits a three-term relation of the 
following form: for every integer vector $\bp = (p,q;r) \in \Z^3$,    
%%%%%%%%%%%%%%%%%%%%%%%%%%%%%%% eqn:3trK %%%%%%%%%%%%%%%%%%%%%%%%%%%%%%%
\begin{equation} \label{eqn:3trK}
\hgK(\ba+\bp;z) = \psi(\ba;\bp) \, r(\ba;z) \, \hgK(\ba;z) + 
\phi(\ba;\bp) \, q(\ba;z) \, \hgK(\ba+\1;z),    
\end{equation}
%%%%%%%%%%%%%%%%%%%%%%%%%%%%%%%%%%%%%%%%%%%%%%%%%%%%%%%%%%%%%%%%%%%%%%%%
where $q(\ba;z)$ and $r(\ba;z)$ are the functions appearing in the 
original three-term relation \eqref{eqn:3trF} for $\hgF(\ba;z)$, 
whereas $\phi(\ba;\bp)$ and $\psi(\ba;\bp)$ are nontrivial rational 
functions of $\ba$ depending uniquely on $\hgK(\ba;z)$ and $\bp$, 
explicit formulas for which can be found in Ebisu \cite[\S 2.2]{Ebisu3}.   
%%%%%%%
\par
%%%%%%%
Problem $\rII$ makes sense not only for $\hgF(\ba;z)$ but 
also for any other member $\hgK(\ba;z)$ of Kummer's solutions.    
Given an integral data $\lambda = (p,q,r;a,b;x)$, if we put 
$k(w;\lambda) := \hgK(\bal(w);x)$, $\tilde{k}(w; \lambda) 
:= \hgK(\bal(w)+\1;x)$, $\Phi(w; \lambda) := \phi(\bal(w); x)$ and 
$\Psi(w; \lambda) := \psi(\bal(w); x)$ with 
$\bal(w) := (pw+a, \, qw+b; \, rw)$, then the three-term relation 
\eqref{eqn:3trK} leads to 
%%%%%%%%%%%%%%%%%%%%%%%%%%%%%%% eqn:3trk %%%%%%%%%%%%%%%%%%%%%%%%%%%%%%
\begin{equation} \label{eqn:3trk}
k(w+1; \lambda) = \Psi(w; \lambda) R(w; \lambda) \cdot k(w; \lambda) + 
\Phi(w;\lambda) Q(w; \lambda) \cdot \tilde{k}(w; \lambda),  
\end{equation}
%%%%%%%%%%%%%%%%%%%%%%%%%%%%%%%%%%%%%%%%%%%%%%%%%%%%%%%%%%%%%%%%%%%%%%%
just as relation \eqref{eqn:3trF} leads to \eqref{eqn:3trf}. 
Notice that $\Phi(w; \lambda)$ and $\Psi(w; \lambda)$ are nontrivial 
rational functions of $w$. 
This observation gives the following important lemma.  
%%%%%%%%%%%%%%%%%%%%%%%%%%%%%% lem:kummer %%%%%%%%%%%%%%%%%%%%%%%%%%%%%%
\begin{lemma} \label{lem:kummer} 
Let $\hgK(\ba; z)$ be any member of Kummer's twenty-four solutions.  
An integral data $\lambda$ is a solution to Problem $\rII$ for 
the function $\hgF(\ba; z)$ that comes from contiguous relations, 
if and only if the same is true for the function $\hgK(\ba; z)$.     
If this is the case, then  
%%%%%%%%%%%%%%%%%%%%%%%%%%%%%% eqn:ocf-k %%%%%%%%%%%%%%%%%%%%%%%%%%%%%%%
\begin{equation} \label{eqn:ocf-k}
\dfrac{k(w+1; \lambda)}{k(w; \lambda)} = \Psi(w; \lambda) \cdot 
R(w; \lambda) \in \C(w),  
\end{equation}
%%%%%%%%%%%%%%%%%%%%%%%%%%%%%%%%%%%%%%%%%%%%%%%%%%%%%%%%%%%%%%%%%%%%%%%%%
which is corresponding to condition \eqref{eqn:ocf} for the original 
function $f(w; \lambda)$. 
\end{lemma}   
%%%%%%%%%%%%%%%%%%%%%%%%%%%%%%%%%%%%%%%%%%%%%%%%%%%%%%%%%%%%%%%%%%%%%%%%
%%%%%%%%%%%%%%%%%%%%%%%%%%%%%% begin proof %%%%%%%%%%%%%%%%%%%%%%%%%%%%%
{\it Proof}. 
Recall that an integral data $\lambda$ is a solution to Problem $\rII$ 
for $\hgF(\ba; z)$ that comes from contiguous relations if and only if 
condition \eqref{eqn:cfcr} is satisfied. 
In view of formula \eqref{eqn:3trk} the corresponding condition for 
$\hgK(\ba; z)$ is $\Phi(w; \lambda) Q(w; \lambda) \equiv 0$ in $\C(w)$.  
But this is just equivalent to condition \eqref{eqn:cfcr}, because 
$\Phi(w; \lambda)$ is nontrivial. 
Now formula \eqref{eqn:3trk} leads to condition \eqref{eqn:ocf-k}. 
\hfill $\Box$
%%%%%%%%%%%%%%%%%%%%%%%%%%%%%% end proof %%%%%%%%%%%%%%%%%%%%%%%%%%%%%%%
%%%%%%%%%%%%%%%%%%%%%%%%%%%%%% rem:kummer %%%%%%%%%%%%%%%%%%%%%%%%%%%%%%
\begin{remark} \label{rem:kummer} 
Any member $\hgK(\ba; z)$ of Kummer's $24$ solutions can be written 
$\hgK(\ba; z) = \mbox{(an elementary factor)} \times  
\hgF(\tilde{\ba}; \tilde{z})$ in terms of the original hypergeometric 
function $\hgF(\ba; z)$ and a certain transformation of variables 
$(\ba; z) \mapsto (\tilde{\ba}; \tilde{z})$. 
So Lemma \ref{lem:kummer} suggests that each $\hgK(\ba; z)$ brings a 
symmetry to Problem $\rII$ for the original function $\hgF(\ba; z)$. 
It may be referred to as an {\sl Ebisu symmetry} because it 
originates from Ebisu's observation \eqref{eqn:3trK}. 
The existence of Ebisu symmetries is an advantage of dealing with  
Problem $\rII$, whereas such a helpful structure cannot be expected 
for Problem $\rI$, although we must keep it in mind that {\sl Ebisu 
symmetries make sense only for those solutions which come from 
contiguous relations}.  
\end{remark}
%%%%%%%%%%%%%%%%%%%%%%%%%%%%%%%%%%%%%%%%%%%%%%%%%%%%%%%%%%%%%%%%%%%%%%%%
\par
%%%%%
In this article we exhibit two special choices of Kummer's solutions 
other than the original one $\hgF(\ba;z)$. 
The resulting Ebisu symmetries will be the main players in this article, 
that is, duality and reciprocity.  
Here one choice of $\hgK(\ba; z)$ is to take 
%%%%%%%%%%%%%%%%%%%%%%%%%%%%%% eqn:hgG %%%%%%%%%%%%%%%%%%%%%%%%%%%%%%%%%
\begin{equation} \label{eqn:hgG}
\hgG(\ba;z) := z^{1-\gamma} \hgF(\alpha-\gamma+1,\beta-\gamma+1;2-\gamma;z),    
\end{equation}
%%%%%%%%%%%%%%%%%%%%%%%%%%%%%%%%%%%%%%%%%%%%%%%%%%%%%%%%%%%%%%%%%%%%%%%%
which is the solution of local exponent $1-\gamma$ at $z=0$ in Riemann 
scheme \eqref{eqn:rs}.  
Note that solutions $\hgF(\ba; z)$ and $\hgG(\ba; z)$ form a linear 
basis of local solutions to the hypergeometric equation \eqref{eqn:hge} 
at $z = 0$, unless $\gamma$ is an integer.  
The other choice is  
%%%%%%%%%%%%%%%%%%%%%%%%%%%%%% eqn:hgH %%%%%%%%%%%%%%%%%%%%%%%%%%%%%%%%%%%%%
\begin{equation} \label{eqn:hgH}
\hgH(\ba;z) := z^{1-\gamma}(1-z)^{\gamma-\alpha-\beta} 
\hgF(1-\alpha,1-\beta;\gamma-\alpha-\beta+1;1-z), 
\end{equation}
%%%%%%%%%%%%%%%%%%%%%%%%%%%%%%%%%%%%%%%%%%%%%%%%%%%%%%%%%%%%%%%%%%%%%%%%%%%%
which is an expression for the local solution of exponent  
$\gamma-\alpha-\beta$ at $z = 1$ in the scheme \eqref{eqn:rs}. 
%%%%%%
\par
%%%%%% 
When $\hgK(\ba; z)$ is either $\hgG(\ba; z)$ or $\hgH(\ba; z)$, 
we use Lemma \ref{lem:kummer} to construct duality or reciprocity, 
where we employ the following notation. 
For $\hgK(\ba; z) = \hgG(\ba; z)$ the functions $k(w;\lambda)$, $\psi(\ba;\bp)$ 
and $\Psi(w; \lambda)$ are denoted by $g(w;\lambda)$, $\psi_g(\ba;\bp)$ and 
$\Psi_g(w; \lambda)$, while for $\hgK(\ba; z) = \hgH(\ba; z)$ they are denoted 
by $h(w;\lambda)$, $\psi_h(\ba;\bp)$ and $\Psi_h(w; \lambda)$.  
Note that 
%%%%%%%%%%%%%%%%%%%%%%%%%%%%%%% eqn:g %%%%% eqn:h %%%%%%%%%%%%%%%%%%%%%%%%%%
\begin{align}
g(w; \lambda) &= x^{1-r w} \hgF((p-r)w+a+1, (q-r)w+b+1; 2 -r w; x), 
\label{eqn:g} \\[1mm]  
h(w;\lambda) &= x^{1-r w}(1-x)^{(r-p-q)w-a-b} \nonumber \\
&\phantom{==}  \times \hgF(1-a-p w, \, 1-b-q w; \, (r-p-q)w+1-a-b; \, 1-x). 
\label{eqn:h} 
\end{align}
%%%%%%%%%%%%%%%%%%%%%%%%%%%%%%%%%%%%%%%%%%%%%%%%%%%%%%%%%%%%%%%%%%%%%%%%%%%%
\par
%%%%%%%  
From a result of Ebisu \cite[Lemma 2.2]{Ebisu3} we have  
%%%%%%%%%%%%%%%%%%%
\begin{align*}
\psi_g(\ba;\bp) &= (-1)^{r-p-q} 
\dfrac{(\alpha)_p(\beta)_q(\gamma-\alpha)_{r-p}(\gamma-\beta)_{r-q}}{(\gamma-1)_r (\gamma)_r},  
\\[2mm] 
\psi_h(\ba;\bp) &= (-1)^{r-p-q} 
\dfrac{(\alpha)_p(\beta)_q(\gamma-\alpha-\beta+1)_{r-p-q}}{(\gamma)_r},     
\end{align*} 
%%%%%%%%%%%%%%%%%%% 
in formula \eqref{eqn:3trK}, from which we find     
%%%%%%%%%%%%%%%%%%%%%% eqn:Psi-g %%%% eqn:Psi-h %%%%%%%%%%%%%%%%%%%%%%%%%%%%%
\begin{align}
\Psi_g(w; \lambda) &= (-1)^{r-p-q} 
\dfrac{(p w+a)_p(q w+b)_q((r-p)w-a)_{r-p}((r-q)w-b)_{r-q}}{(r w-1)_r (r w)_r}, 
\label{eqn:Psi-g} \\[2mm]
\Psi_h(w; \lambda) &= (-1)^{r-p-q} 
\dfrac{(p w+a)_p(q w+b)_q((r-p-q)w-a-b+1)_{r-p-q}}{(r w)_r}, \label{eqn:Psi-h}  
\end{align} 
%%%%%%%%%%%%%%%%%%%%%%%%%%%%%%%%%%%%%%%%%%%%%%%%%%%%%%%%%%%%%%%%%%%%%%%%%%%%%
in formula \eqref{eqn:ocf-k}, where $(s)_n := 
\varGamma(s+n)/\varGamma(s)$ is Pochhammer's symbol or the rising 
factorial.     
Solution \eqref{eqn:hgG} and formula \eqref{eqn:Psi-g} will be used to 
construct duality in \S\ref{sec:duality}, while solution \eqref{eqn:hgH} 
and formula \eqref{eqn:Psi-h} will be employed to construct reciprocity 
in \S\ref{sec:r-cf}, respectively.  
%%%%%%%%%%%%%%%%%%%%%%%%%%%%% sec:duality %%%%%%%%%%%%%%%%%%%%%%%%%%%%%%%%%%% 
\section{Duality} \label{sec:duality} 
%%%%%%%%%%%%%%%%%%%%%%%%%%%%%%%%%%%%%%%%%%%%%%%%%%%%%%%%%%%%%%%%%%%%%%%%%%%%
Applying Lemma \ref{lem:kummer} to $\hgK(\ba;z) = \hgG(\ba;z)$ leads to the 
duality \eqref{eqn:duality} in Definition \ref{def:duality}. 
%%%%%%%%%%%%%%%%%%%%%%%%%% lem:duality1 %%%%%%%%%%%%%%%%%%%%%%%%%%%%%%%%%%%%
\begin{lemma} \label{lem:duality1} 
Let $\lambda = (p,q,r;a,b;x)$ be an integral data in domain \eqref{eqn:real2}.   
If $\lambda$ is a solution to Problem $\rII$ that comes from contiguous 
relations, with rational function $R(w; \lambda)$ in condition 
\eqref{eqn:ocf}, then its dual $\lambda' = (p',q',r';a',b';x')$ is also 
a solution to Problem $\rII$ that comes from contiguous relations,  
with the corresponding rational function  
%%%%%%%%%%%%%%%%%%%%%%%%%% eqn:R-p %%%%%%%%%%%%%%%%%%%%%%%%%%%%%%%%%%%%%
\begin{equation} \label{eqn:R-p}
R(w;\lambda') = 
\dfrac{x^{-r}(1-x)^{r-p-q}}{\Psi_g(w'; \lambda) \, R(w';\lambda)},      
\end{equation}
%%%%%%%%%%%%%%%%%%%%%%%%%%%%%%%%%%%%%%%%%%%%%%%%%%%%%%%%%%%%%%%%%%%%%%%%%%%%
where the function $\Psi_g(w; \lambda)$ is given by formula 
\eqref{eqn:Psi-g} and $w \mapsto w'$ is the reflection   
%%%%%%%%%%%%%%%%%%%%%%%%%%% eqn:w-p %%%%%%%%%%%%%%%%%%%%%%%%%%%%%%%%%%%%
\begin{equation} \label{eqn:w-p}
w' := \frac{2}{r}-1-w. 
\end{equation}
%%%%%%%%%%%%%%%%%%%%%%%%%%%%%%%%%%%%%%%%%%%%%%%%%%%%%%%%%%%%%%%%%%%%%%%%%%%% 
\end{lemma} 
%%%%%%%%%%%%%%%%%%%%%%%%%% begin proof %%%%%%%%%%%%%%%%%%%%%%%%%%%%%%%%%%%%%
{\it Proof}. 
Replacing $w$ by $w+1 = \frac{2}{r} - w'$ in formula \eqref{eqn:g}, 
we observe that   
%%%%%
\[
\begin{split}
g(w+1;\lambda) 
&= \resizebox{0.8\hsize}{!}{$\displaystyle
x^{1-r\left(\frac{2}{r} - w'\right)} 
\hgF\left((p-r)\left(\ts\frac{2}{r} - w'\right)+a+1, \, 
(q-r)\left(\ts\frac{2}{r} - w'\right)+b+1; 
2 - r\left(\ts\frac{2}{r}-w'\right); x\right)$} \\
&= x^{r w'-1} \hgF\left((r-p) w' - a', \, (r-q) w' - b'; r w'; x\right) \\
&= x^{r w'-1} (1-x)^{(p+q-r)w'+a'+b'} f(w'; \lambda'), 
\end{split}
\]
%%%%%
where definition \eqref{eqn:duality} and Euler's transformation 
\cite[(7b)]{Iwasaki} are used in the second and third 
equalities respectively. 
Since the shift $w \mapsto w-1$ is equivalent to 
$w' \mapsto w'+1$, we have  
%%%%%%%%%%%%%%%%%%%%%%%%%%%%%%%%%%% eqn:f-p %%%%%%%%%%%%%%%%%%%%%%%%%%%%%%%%%%
\begin{align} 
f(w'; \lambda') 
&= x^{1-r w'} (1-x)^{(r-p-q) w'-a'-b'} g(w+1;\lambda),  \nonumber \\
f(w'+1; \lambda') 
&= x^{1-r (w'+1)} (1-x)^{(r-p-q) (w'+1)-a'-b'} g(w;\lambda), 
\label{eqn:f-p} 
\end{align}   
%%%%%%%%%%%%%%%%%%%%%%%%%%%%%%%%%%%%%%%%%%%%%%%%%%%%%%%%%%%%%%%%%%%%%%%%%%%%%%
which together with formula \eqref{eqn:ocf-k} in Lemma \ref{lem:kummer} 
yields 
%%%%%%%%%%
\[
R(w'; \lambda') := \dfrac{f(w'+1; \lambda')}{f(w'; \lambda')} 
= x^{-r}(1-x)^{r-p-q} \cdot \dfrac{g(w;\lambda)}{g(w+1;\lambda)} = 
\dfrac{ x^{-r}(1-x)^{r-p-q}}{\Psi_g(w; \lambda) \, R(w; \lambda)}.   
\]
%%%%%%%%%%
Replacing $w$ by $w'$ in the above and noting $w'' = w$, 
we obtain formula \eqref{eqn:R-p}. \hfill $\Box$ \par\medskip
%%%%%%%%%%%%%%%%%%%%%%%%%% end proof %%%%%%%%%%%%%%%%%%%%%%%%%%%%%%%%%%%%
Using Lemma \ref{lem:duality1} in domain $\cE^{*-}$ leads to an  
immediate consequence. 
%%%%%%%%%%%%%%%%%%%%%%%%%% cor:duality1 %%%%%%%%%%%%%%%%%%%%%%%%%%%%%%%%%
\begin{corollary} \label{cor:duality1}  
Duality $\lambda \mapsto \lambda'$ in \eqref{eqn:duality} induces a 
self-bijection on the set of all solutions $\lambda = 
(p,q,r;a,b;x) \in \cE^{*-}$ with $q \in \Z$.    
\end{corollary}   
%%%%%%%%%%%%%%%%%%%%%%%%%%%%%%%%%%%%%%%%%%%%%%%%%%%%%%%%%%%%%%%%%%%%%%%%%
%%%%%%%%%%%%%%%%%%%%%%%%%% begin proof %%%%%%%%%%%%%%%%%%%%%%%%%%%%%%%%%%
{\it Proof}. 
In view of definition \eqref{eqn:duality} the duality 
$\lambda \mapsto \lambda'$ is a bijection $\cE^{*-} \to \cE^{*-}$ in 
the data level. 
By assertion (3) of \cite[Theorem 2.5]{Iwasaki} and Lemma \ref{lem:duality1}, 
it induces a bijection in the solution level among all solutions 
$\lambda \in \cE^{*-}$ to Problem $\rII$ (and to Problem $\rI$ 
by \cite[Theorem 2.1]{Iwasaki}) with $q \in \Z$.  \hfill $\Box$ \par\medskip
%%%%%%%%%%%%%%%%%%%%%%%%%% end proof %%%%%%%%%%%%%%%%%%%%%%%%%%%%%%%%%%%%
The same results holds true for $(\rA)$-solutions in $\cD^-$ since 
the duality is also a bijection $\cD^- \to \cD^-$ in the data level, 
but in fact we are able to obtain more detailed results on $\cD^-$. 
%%%%%%%%%%%%%%%%%%%%%%%%%% lem:duality2 %%%%%%%%%%%%%%%%%%%%%%%%%%%%%%%%%
\begin{lemma} \label{lem:duality2} 
If $\lambda = (p,q,r;a,b;x) \in \cD^-$ is an $(\rA)$-solution to 
Problem $\rII$, then its dual $\lambda' = (p,q,r;a',b';x) \in \cD^-$ is 
also an $(\rA)$-solution to the same problem with  
%%%%%%%%%%%%%%%%%%%%%%%%%% eqn:R-p2 %%%%%%%%%%%%%%%%%%%%%%%%%%%%%%%%%
\begin{equation} \label{eqn:R-p2}
R(w;\lambda') = d \cdot  
\frac{\prod_{i=0}^{r-1} \left(w+\frac{i}{r}
\right)}{\prod_{i=1}^{r} \left(w+v_i' \right)}, 
\end{equation}
%%%%%%%%%%%%%%%%%%%%%%%%%%%%%%%%%%%%%%%%%%%%%%%%%%%%%%%%%%%%%%%%%%%%%%%%%
where $d$ and $v_1', \dots, v_r'$ are defined by formulas 
\eqref{eqn:d} and \eqref{eqn:vi-p} respectively.  
\end{lemma} 
%%%%%%%%%%%%%%%%%%%%%%%%%%%%%%%%%%%%%%%%%%%%%%%%%%%%%%%%%%%%%%%%%%%%%%%%%
%%%%%%%%%%%%%%%%%%%%%%%%%% begin proof %%%%%%%%%%%%%%%%%%%%%%%%%%%%%%%%%% 
{\it Proof}. 
By assertion (2) of \cite[Theorem 2.3]{Iwasaki} we have the gamma product 
formula \eqref{eqn:gpf2}, so that the rational function 
$R(w; \lambda)$ in formula \eqref{eqn:R} is given by  
%%%%%%%%%%%%%%%%%%%%%%%%%%%%% eqn:R2 %%%%%%%%%%%%%%%%%%%%%%%%%%%%%%%%%%%
\begin{equation} \label{eqn:R2}
R(w; \lambda) = d \cdot \dfrac{\prod_{i=0}^{r-1} \left(w+\frac{i}{r}
\right)}{\prod_{i=1}^{r} \left(w+v_i\right)},   
\end{equation}
%%%%%%%%%%%%%%%%%%%%%%%%%%%%%%%%%%%%%%%%%%%%%%%%%%%%%%%%%%%%%%%%%%%%%%%%%
where $d$ is defined in formula \eqref{eqn:d}. 
On the other hand, formula \eqref{eqn:Psi-g} can be rewritten 
%%%%%
\[
\begin{split}
\Psi_g(w; \lambda) &= \dfrac{p^p q^q (r-p)^{r-p} (r-q)^{r-q}}{r^{2r}} \\ 
&\phantom{=} \times 
\dfrac{\prod_{i=0}^{p-1}\left(w+{\ts\frac{i+a}{p}}\right) 
\prod_{i=0}^{q-1}\left(w+{\ts\frac{i+b}{q}}\right) 
\prod_{i=0}^{r-p-1}\left(w+{\ts\frac{i-a}{r-p}}\right) 
\prod_{i=0}^{r-q-1}\left(w+{\ts\frac{i-b}{r-q}}\right)}{
\prod_{i=0}^{r-1} \left(w+\frac{i-1}{r}\right) 
\prod_{i=0}^{r-1} \left(w+\frac{i}{r}\right)},     
\end{split}  
\]
%%%%%
where $(-1)^{r-p-q} = 1$ is taken into account, which follows from 
assertion (2) of \cite[Theorem 2.2]{Iwasaki}. 
Thus taking the product of equations \eqref{eqn:Psi-g} and \eqref{eqn:R2} 
and using definition \eqref{eqn:vi-s}, we have 
%%%%%%
\[
\dfrac{1}{\Psi_g(w; \lambda) \, R(w; \lambda)} = 
d \cdot x^r\cdot (1-x)^{p+q-r} \cdot 
\dfrac{\prod_{i=0}^{r-1} \left(w+\frac{i-1}{r}
\right)}{\prod_{i=1}^r \left(w+v_i^*\right)}. 
\]
%%%%%%
Replacing $w$ by $w'$ in the above, where $w'$ is defined by 
\eqref{eqn:w-p}, and using   
$\prod_{i=0}^{r-1}\left(w'+ \frac{i-1}{r}\right) = (-1)^r 
\prod_{i=0}^{r-1}\left(w + \frac{i}{r}\right)$ and 
$\prod_{i=0}^{r-1}\left(w'+ v_i^*\right) = (-1)^r 
\prod_{i=0}^{r-1}\left(w + v_i'\right)$, we have  
%%%%%%%
\[
\dfrac{1}{\Psi_g(w'; \lambda) \, R(w'; \lambda)} = 
d \cdot x^r\cdot (1-x)^{p+q-r} \cdot 
\dfrac{\prod_{i=0}^{r-1} \left(w+\frac{i}{r}
\right)}{\prod_{i=1}^r \left(w+v_i'\right)}. 
\]
%%%%%%%
Substituting this into formula \eqref{eqn:R-p} leads to 
the desired formula \eqref{eqn:R-p2}. 
\hfill $\Box$ 
%%%%%%%%%%%%%%%%%%%%%%%%%%%%%% end proof %%%%%%%%%%%%%%%%%%%%%%%%%%%%%%%%
%%%%%%%%%%%%%%%%%%%%%%%%%%%%%% prop:duality %%%%%%%%%%%%%%%%%%%%%%%%%%%%%
\begin{proposition} \label{prop:duality} 
If $\lambda = (p,q,r;a,b;x) \in \cD^-$ is an $(\rA)$-solution to  
Problem $\rI$ with gamma product formula \eqref{eqn:gpf2}, then 
there exists a positive constant $C' > 0$ such that 
%%%%%%%%%%%%%%%%%%%% eqn:gpf3 %%%%% eqn:gpf-g %%%%%%%%%%%%%%%%%%%%%%%%%%%
\begin{align}
f(w; \lambda') &= C' \cdot d^w \cdot 
\dfrac{\prod_{i=0}^{r-1} \varGamma\left(w+\frac{i}{r}\right)}{\prod_{i=1}^{r} \varGamma\left(w+v_i'\right)}, 
\label{eqn:gpf3} \\[2mm]
g(w; \lambda) &= D \cdot \delta^w \cdot
\dfrac{\prod_{i=1}^{r} \sin \pi \left(w+v_i^*\right)}{\prod_{i=0}^{r-1} \sin \pi \left(w+\frac{i-1}{r}\right)}
\cdot  
\dfrac{\prod_{i=1}^{r} \varGamma\left(w+v_i^*\right)}{\prod_{i=0}^{r-1} \varGamma\left(w+\frac{i-1}{r}\right)}, 
\label{eqn:gpf-g}  
\end{align} 
%%%%%%%%%%%%%%%%%%%%%%%%%%%%%%%%%%%%%%%%%%%%%%%%%%%%%%%%%%%%%%%%%%%%%%%%%
where $d$, $v_i'$ and $v_i^*$ are given by formulas 
\eqref{eqn:d}, \eqref{eqn:vi-p} and \eqref{eqn:vi-s}, whereas 
$D$ and $\delta$ are defined by  
%%%%%%%%%%%%%%%%%%%%%%%%%%%%%% eqn:D-d %%%%%%%%%%%%%%%%%%%%%%%%%%%%%%%%%%
\begin{equation} \label{eqn:D-d}
D := \dfrac{x \cdot d^{2/r} \cdot C'}{(1-x)^{a+b}}, \qquad 
\delta := \dfrac{(1-x)^{r-p-q}}{d \cdot x^r}. 
\end{equation} 
%%%%%%%%%%%%%%%%%%%%%%%%%%%%%%%%%%%%%%%%%%%%%%%%%%%%%%%%%%%%%%%%%%%%%%%%%
\end{proposition}
%%%%%%%%%%%%%%%%%%%%%%%%%%%%%%%%%%%%%%%%%%%%%%%%%%%%%%%%%%%%%%%%%%%%%%%%% 
%%%%%%%%%%%%%%%%%%%%%%%%%%%%%%% begin proof %%%%%%%%%%%%%%%%%%%%%%%%%%%%% 
{\it Proof}. 
By Lemma \ref{lem:duality2}, $\lambda'$ is a solution to Problem $\rII$ 
in $\cD^-$ with rational function $R(w;\lambda')$ in formula 
\eqref{eqn:R-p2}.   
Theorem 2.1 of \cite{Iwasaki} then implies that it leads back to a 
solution to Problem $\rI$, which is exactly the gamma product formula 
in \eqref{eqn:gpf3}. 
Next we have 
%%%%%
\[
\begin{split}
g(w;\lambda) 
&= x^{r(w'+1)-1} (1-x)^{a'+b'-(r-p-q)(w'+1)} \, f(w'+1; \lambda') \\
&= \dfrac{x}{(1-x)^{a+b}} \cdot \left\{\frac{(1-x)^{r-p-q}}{x^r} \right\}^w 
\cdot f\left(\textstyle \frac{2}{r}-w; \lambda' \right) \\
&= \dfrac{x}{(1-x)^{a+b}} \cdot \left\{\frac{(1-x)^{r-p-q}}{x^r} \right\}^w 
\cdot C' \cdot d^{\frac{2}{r} -w} \cdot 
\dfrac{\prod_{i=0}^{r-1} 
\varGamma\left(\frac{2}{r}-w+\frac{i}{r} \right)}{\prod_{i=1}^r 
\varGamma\left(\frac{2}{r}-w + v_i'\right)} \\
&= \dfrac{C' \cdot x \cdot d^{\frac{2}{r}}}{(1-x)^{a+b}} \cdot 
\left\{\frac{(1-x)^{r-p-q}}{d \cdot x^r} \right\}^w  
\dfrac{\prod_{i=0}^{r-1} 
\varGamma\left(1-(w+\frac{(r-1-i)-1}{r}) \right)}{\prod_{i=1}^r 
\varGamma\left(1-(w + v_i^*)\right)} \\
&= D \cdot \delta^w \cdot \dfrac{\prod_{i=0}^{r-1} 
\varGamma\left(1-(w+\frac{i-1}{r}) \right)}{\prod_{i=1}^r 
\varGamma\left(1-(w + v_i^*)\right)}
\end{split}
\]
%%%%%
where the first equality follows from \eqref{eqn:f-p}, the second 
from \eqref{eqn:duality} and \eqref{eqn:w-p}, the third from 
\eqref{eqn:gpf3}, the fourth from \eqref{eqn:vi-p}, the fifth 
from \eqref{eqn:D-d} and the replacement of indices  
$i \leftrightarrow r-1-i$, respectively.    
Finally we use Euler's reflection formula for the gamma 
function \cite[Theorem 1.2.1]{AAR}:  
%%%%%%%%%%%%%%%%%%%%%%%%%%%%%%% eqn:refl %%%%%%%%%%%%%%%%%%%%%%%%%%%%%%%%
\begin{equation} \label{eqn:refl}
\varGamma(w) \varGamma(1-w) = \dfrac{\pi}{\sin \pi w},  
\end{equation} 
%%%%%%%%%%%%%%%%%%%%%%%%%%%%%%%%%%%%%%%%%%%%%%%%%%%%%%%%%%%%%%%%%%%%%%%%% 
to establish the desired formula \eqref{eqn:gpf-g} 
\hfill $\Box$ \par\medskip
%%%%%%%%%%%%%%%%%%%%%%%%%%%% end proof %%%%%%%%%%%%%%%%%%%%%%%%%%%%%%%%%%
With the proof of formula \eqref{eqn:gpf3} in Proposition 
\ref{prop:duality}, Theorem \ref{thm:gpf-d} has been established. 
%%%%%%%%%%%%%%%%%%%%%%%%%%%%% sec:r-cf %%%%%%%%%%%%%%%%%%%%%%%%%%%%%%%%%% 
\section{Reciprocity and Connection Formula} \label{sec:r-cf} 
%%%%%%%%%%%%%%%%%%%%%%%%%%%%%%%%%%%%%%%%%%%%%%%%%%%%%%%%%%%%%%%%%%%%%%%%
We shall show that applying Lemma \ref{lem:kummer} to $\hgK(\ba;z) = 
\hgH(\ba;z)$ yields the reciprocity \eqref{eqn:recip} in 
Definition \ref{def:recip}. 
First we observe that under the involution 
$\lambda \mapsto \check{\lambda}$, the transformation   
%%%%%%%%%%%%%%%%%%%%%%%%%%% eqn:w-c %%%%%%%%%%%%%%%%%%%%%%%%%%%%%%%%%%%%
\begin{equation} \label{eqn:w-c}
w \mapsto \check{w} := w + c, \qquad 
c = c(\lambda) := \frac{1-a-b}{r-p-q}  
\end{equation}
%%%%%%%%%%%%%%%%%%%%%%%%%%%%%%%%%%%%%%%%%%%%%%%%%%%%%%%%%%%%%%%%%%%%%%%%
also yields an involution, since definition \eqref{eqn:recip} implies 
$\check{c}:= c(\check{\lambda}) = - c(\lambda) = -c$.     
%%%%%%%%%%%%%%%%%%%%%%%%%% lem:r-cf1 %%%%%%%%%%%%%%%%%%%%%%%%%%%%%%%%%%%%
\begin{lemma} \label{lem:r-cf1} 
Let $\lambda = (p,q,r;a,b;x)$ be an integral data in domain \eqref{eqn:real1}. 
If $\lambda$ is a solution to Problem $\rII$ that comes from contiguous 
relations, with rational function $R(w;\lambda)$ in condition 
\eqref{eqn:ocf}, then its reciprocal $\check{\lambda} 
= (\check{p},\check{q},\check{r};\check{a},\check{b}; \check{x})$ 
is also a solution to Problem $\rII$ that comes from contiguous 
relations, with the corresponding rational function 
%%%%%%%%%%%%%%%%%%%%%%%%%% eqn:R-c %%%%%%%%%%%%%%%%%%%%%%%%%%%%%%%%%%%%%
\begin{equation} \label{eqn:R-c}
R(w; \check{\lambda}) = x^r (1-x)^{p+q-r} \cdot 
\Psi_h(\hat{w}; \lambda) R(\hat{w}; \lambda), 
\end{equation}
%%%%%%%%%%%%%%%%%%%%%%%%%%%%%%%%%%%%%%%%%%%%%%%%%%%%%%%%%%%%%%%%%%%%%%%% 
where $\Psi_h(w; \lambda)$ is given by formula \eqref{eqn:Psi-h} and 
$\hat{w}$ is defined by
%%%%%%%%%%%%%%%%%%%%%%%%% eqn:w-h %%%%%%%%%%%%%%%%%%%%%%%%%%%%%%%%%%%%%%
\begin{equation} \label{eqn:w-h}
\hat{w} := w - c, \qquad c = c(\lambda) := \frac{1-a-b}{r-p-q}.  
\end{equation}   
\end{lemma}
%%%%%%%%%%%%%%%%%%%%%%%%% begin proof %%%%%%%%%%%%%%%%%%%%%%%%%%%%%%%%%
{\it Proof}. 
Definitions \eqref{eqn:recip} and \eqref{eqn:w-c} imply  
$(r-p-q) w+1-a-b = \check{r} \check{w}$ and   
%%%%%%%%%%%
\[
1-a - p w = 1-a - p (\check{w} - c) = \check{p} \check{w} + \check{a},    
\]
%%%%%%%%%%% 
and similarly $1-b-q w = \check{q} \check{w} + \check{b}$.  
Substituting these into formula \eqref{eqn:h} we find  
$h(w; \lambda) = x^{1-r w}(1-x)^{(r-p-q)w-a-b} \,   
f(\check{w}; \check{\lambda})$, or in other words,  
%%%%%%%%%%%%%%%%%%%%%%%%%%% eqn:f-h %%%%%%%%%%%%%%%%%%%%%%%%%%%%%%%%%%
\begin{equation} \label{eqn:f-h}
f(\check{w}; \check{\lambda}) = x^{r w-1}(1-x)^{(p+q-r)w+a+b} 
\, h(w; \lambda).  
\end{equation}
%%%%%%%%%%%%%%%%%%%%%%%%%%%%%%%%%%%%%%%%%%%%%%%%%%%%%%%%%%%%%%%%%%%%%%
Since increasing $\check{w}$ by $1$ is equivalent to increasing $w$ 
by $1$, we also have  
%%%%%%%%%%%
\[
f(\check{w}+1; \check{\lambda}) = 
x^{r(w+1)-1}(1-x)^{(p+q-r)(w+1)+a+b} \, h(w+1; \lambda),    
\]
%%%%%%%%%%%
so that formula \eqref{eqn:ocf-k} in Lemma \ref{lem:kummer} and 
definition \eqref{eqn:w-c} yield  
%%%%%%%%%%%
\[
\begin{split}
R(\check{w}; \check{\lambda}) &= 
\dfrac{f(\check{w}+1; \check{\lambda})}{f(\check{w}; \check{\lambda})} 
= x^r (1-x)^{p+q-r} \dfrac{h(w+1; \lambda)}{h(w; \lambda)} = 
x^r (1-x)^{p+q-r} \cdot \Psi_h(w; \lambda) R(w;\lambda) \\
&= x^r (1-x)^{p+q-r} \cdot \Psi_h(\check{w} - c; \lambda) \, 
R(\check{w} - c; \lambda).  
\end{split} 
\]
%%%%%%%%%%%
Since $\check{w}$ is an indeterminate variable, we can replace 
$\check{w}$ by $w$ in the above to obtain 
%%%%%%%%%%%
\[
\begin{split}
R(w; \check{\lambda}) &= 
x^r (1-x)^{p+q-r} \cdot \Psi_h(w - c; \lambda) R(w - c; \lambda) \\
&= x^r (1-x)^{p+q-r} \cdot \Psi_h(\hat{w}; \lambda) \, R(\hat{w}; \lambda), 
\end{split} 
\]
%%%%%%%%%%% 
where definition \eqref{eqn:w-h} is used in the last equality. 
\hfill $\Box$ \par\medskip 
%%%%%%%%%%%%%%%%%%%%%%%%%% end proof %%%%%%%%%%%%%%%%%%%%%%%%%%%%%%%%%%%%
Using Lemma \ref{lem:r-cf1} in domains $\cE^{*-}$ and $\cE^{-*}$ 
yields the following direct consequence. 
%%%%%%%%%%%%%%%%%%%%%%%%%% cor:r-cf1 %%%%%%%%%%%%%%%%%%%%%%%%%%%%%%%%%%%%
\begin{corollary} \label{cor:r-cf1}  
Reciprocity $\lambda \mapsto \check{\lambda}$ in \eqref{eqn:recip} 
induces a bijection between the set of all solutions $\lambda = 
(p,q,r;a,b;x) \in \cE^{*-}$ with $q \in \Z$ and the set of all solutions 
$\check{\lambda} = 
(\check{p}, \check{q}, \check{r}; \check{a}, \check{b}; \check{x}) 
\in \cE^{-*}$ with $\check{p} \in \Z$.   
\end{corollary}   
%%%%%%%%%%%%%%%%%%%%%%%%%%%%%%%%%%%%%%%%%%%%%%%%%%%%%%%%%%%%%%%%%%%%%%%%%
%%%%%%%%%%%%%%%%%%%%%%%%%% begin proof %%%%%%%%%%%%%%%%%%%%%%%%%%%%%%%%%%
{\it Proof}.  
In view of definition \eqref{eqn:recip} the reciprocity  
$\lambda \mapsto \check{\lambda}$ is a bijection $\cE^{*-} \to \cE^{-*}$ 
in the data level. 
By \cite[Theorem 2.5]{Iwasaki} and Lemma \ref{lem:r-cf1} it induces  
a bijection in the solution level between the solutions 
$\lambda \in \cE^{*-}$ with $q \in \Z$ and the solutions 
$\check{\lambda} \in \cE^{-*}$ with $\check{p} \in \Z$.  
\hfill $\Box$ \par\medskip 
%%%%%%%%%%%%%%%%%%%%%%%%%% end proof %%%%%%%%%%%%%%%%%%%%%%%%%%%%%%%%%%%%
Although the reciprocity is also a bijection $\cD^- \to \cF^-$ in the 
data level, it does not immediately induce a bijection in the 
integral solution level as in Corollary \ref{cor:r-cf1}. 
This is because we have not yet known whether every integral solution 
in $\cF^-$ to Problem $\rII$ comes from contiguous relations, 
so that the backward reciprocity $\cF^- \to \cD^-$ is not established 
yet in the solution level (see Remark \ref{rem:kummer}). 
This issue is postponed until it is settled in Proposition 
\ref{prop:cfcr-F}. 
In the rest of this section we develop a detailed study of the forward 
reciprocity $\cD^- \to \cF^-$ in the integral solution level by using 
the connection formula for hypergeometric functions.    
%%%%%%%%%%%%%%%%%%%%%%%%% lem:r-cf2 %%%%%%%%%%%%%%%%%%%%%%%%%%%%%%%%%%%
\begin{lemma} \label{lem:r-cf2} 
If $\lambda = (p,q,r;a,b;x) \in \cD^-$ is an $(\rA)$-solution  
with GPF \eqref{eqn:gpf2}, then  
%%%%%%%%%%%%%%%%%%%%%%%%% eqn:h2 %%%%%%%%%%%%%%%%%%%%%%%%%%%%%%%%%%%%%%
\begin{equation} \label{eqn:h2}
\resizebox{0.85\hsize}{!}{$\displaystyle 
h(w;\lambda) = \varGamma((r-p-q)w+1-a-b) \cdot \tilde{d}^w \cdot \chi(w)   
\cdot \dfrac{\prod_{i=0}^{p-1}\varGamma\left(w+\frac{i+a}{p}\right) 
\prod_{i=0}^{q-1}\varGamma\left(w+\frac{i+b}{q}\right)}{\prod_{i=1}^{r} 
\varGamma\left(w+v_i\right)}$},  
\end{equation}
%%%%%%%%%%%%%%%%%%%%%%%%%%%%%%%%%%%%%%%%%%%%%%%%%%%%%%%%%%%%%%%%%%%%%%%
where the constant $\tilde{d}$ and the function $\chi(w)$ are given by  
%%%%%%%%%%% eqn:d-t %%% eqn:chi%% eqn:C1 %% eqn:C2 %%%%%%%%%%%%%%%%%%%%
\begin{align}  
\tilde{d} &= \sqrt{p^p q^q (r-p)^{p-r} (r-q)^{q-r} x^{-r} (1-x)^{r-p-q}}, 
\label{eqn:d-t} \\[1mm] 
\chi(w) &= C_1 \cdot \dfrac{\sin \pi(p w+a) \sin \pi(q w+b)}{\sin(\pi r w)} 
+ C_2 \cdot \dfrac{\prod_{i=1}^r \sin \pi(w+v_i^*)}{\prod_{i=0}^{r-1} 
\sin \pi\left(w+\frac{i-1}{r}\right)}, \label{eqn:chi} \\[1mm]
C_1 &:= 2 \, (2\pi)^{(r-p-q-1)/2} \cdot C \cdot p^{a-1/2} q^{b-1/2} r^{1/2}, 
\label{eqn:C1} \\    
C_2 &:= (2\pi)^{(r-p-q-1)/2} \cdot D \cdot (r-p)^{a+1/2} (r-q)^{b+1/2} r^{-3/2}, 
\label{eqn:C2} 
\end{align} 
%%%%%%%%%%%%%%%%%%%%%%%%%%%%%%%%%%%%%%%%%%%%%%%%%%%%%%%%%%%%%%%%%%%%%%%
with $C$ and $D$ being the constants in formulas \eqref{eqn:gpf2} and 
\eqref{eqn:D-d} respectively.  
\end{lemma}
%%%%%%%%%%%%%%%%%%%%%%%%%%%%%%%%%%%%%%%%%%%%%%%%%%%%%%%%%%%%%%%%%%%%%%%
%%%%%%%%%%%%%%%%%%%%%%%%% begin proof %%%%%%%%%%%%%%%%%%%%%%%%%%%%%%%%%
{\it Proof}. 
A connection formula in Erd\'elyi 
\cite[Chap.\! I\!I, \S 2.9, formula (43)]{Erdelyi}) reads    
%%%%%%%%%%%%%%%%
\[
\hgH(\ba;z) = 
\dfrac{\varGamma(\gamma-\alpha-\beta+1) \varGamma(1-\gamma)}{\varGamma(1-\alpha) \varGamma(1-\beta)} \, \hgF(\ba;z) +
\dfrac{\varGamma(\gamma-\alpha-\beta+1) \varGamma(\gamma-1)}{\varGamma(\gamma-\alpha) \varGamma(\gamma-\beta)} \, \hgG(\ba;z).   
\]
%%%%%%%%%%%%%%%%
Substituting $\ba = \bal(w) := (p w+a, \, q w+b; \, r w)$ and $z = x$ 
into the connection formula and using the definitions of 
$g(w; \lambda)$ and $h(w; \lambda)$ (just before formulas 
\eqref{eqn:g} and \eqref{eqn:h}), we have 
%%%%%%%%%%%%%%%%%%%%%%%%%%% eqn:cf %%%%%%%%%%%%%%%%%%%%%%%%%%%%%%%%%%%%
\begin{equation} \label{eqn:cf} 
h(w; \lambda) = C_f(w) \, f(w; \lambda) + C_g(w) \, g(w; \lambda), 
\end{equation}
%%%%%%%%%%%%%%%%%%%%%%%%%%%%%%%%%%%%%%%%%%%%%%%%%%%%%%%%%%%%%%%%%%%%%%%
where the connection coefficients $C_f(w)$ and $C_g(w)$ are given by 
%%%%%%%%%%%%%%
\begin{align*} 
C_f(w) 
&= \frac{\varGamma((r-p-q)w+1-a-b)\, \varGamma(1-r w)}{\varGamma(1-a-p w)\, \varGamma(1-b-q w)}, 
\\[1mm] 
C_g(w)
&= \frac{\varGamma((r-p-q)w+1-a-b) \, \varGamma(r w-1)}{\varGamma((r-p)w-a) \varGamma((r-q)w-b)}.     
\end{align*}
%%%%%
\par
%%%%%
The connection coefficient $C_f(w)$ can be written  
%%%%%%%% 
\[
\begin{split}
C_f(w) 
&= \varGamma((r-p-q)w+1-a-b) \cdot 
\dfrac{\sin \pi(p w+a) \, \sin \pi(q w+b)}{\pi \, \sin(\pi r w)} \cdot 
\dfrac{\varGamma(p w+a) \, \varGamma(q w+b)}{\varGamma(r w)} \\
&= \varGamma((r-p-q)w+1-a-b) \cdot 
\dfrac{\sin \pi(p w+a) \, \sin \pi(q w+b)}{\pi \, \sin(\pi r w)} \\
&\phantom{=} \times (2 \pi)^{(r-p-q+1)/2} p^{a-1/2} q^{b-1/2} r^{1/2} 
\left(\frac{p^p q^q}{r^r}\right)^w 
\dfrac{\prod_{i=0}^{p-1} \varGamma\left(w+\frac{i+a}{p}\right)
\prod_{i=0}^{q-1} \varGamma\left(w+\frac{i+b}{q}\right)}{\prod_{i=0}^{r-1} 
\varGamma\left(w+\frac{i}{r}\right)}, 
\end{split}
\]
%%%%%%%%
where the first equality follows from the reflection formula 
\eqref{eqn:refl}, while the second equality from Gauss's multiplication 
formula \eqref{eqn:mult} for the gamma function.  
The above expression for $C_f(w)$ is multiplied by formula 
\eqref{eqn:gpf2} to yield 
%%%%%%%%%%%%%%%%%%%%%%%%%%%%%% eqn:Cf %%%%%%%%%%%%%%%%%%%%%%%%%%%%%%%%%%%%
\begin{equation} \label{eqn:Cf}
\begin{split}
C_f(w) \, f(w;\lambda) 
&= \varGamma((r-p-q)w+1-a-b) \cdot 
\dfrac{\sin \pi(p w+a) \, \sin \pi(q w+b)}{\sin(\pi r w)} \\
&\phantom{=} \times C_1 \cdot 
\tilde{d}^w \cdot 
\dfrac{\prod_{i=0}^{p-1} \varGamma\left(w+\frac{i+a}{p}\right)
\prod_{i=0}^{q-1} \varGamma\left(w+\frac{i+b}{q}\right)}{\prod_{i=0}^{r-1} 
\varGamma(w+v_i)}, 
\end{split}
\end{equation}
%%%%%%%%%%%%%%%%%%%%%%%%%%%%%%%%%%%%%%%%%%%%%%%%%%%%%%%%%%%%%%%%%%%%%%%%%%
where one uses definition \eqref{eqn:C1} and relation 
$d \cdot p^p q^q r^{-r} = \tilde{d}$ that follows from \eqref{eqn:d} 
and \eqref{eqn:d-t}. 
%%%%%
\par
%%%%%
In a similar manner the multiplication formula \eqref{eqn:mult} allows us 
to write   
%%%%%
\[
\begin{split}
C_g(w) &= \varGamma((r-p-q)w+1-a-b) \cdot (2\pi)^{(r-p-q-1)/2} 
(r-p)^{a+1/2} (r-q)^{b+1/2} r^{-3/2} \\
&\phantom{=} \times 
\left\{ \frac{r^r}{(r-p)^{r-p}(r-q)^{r-q}} \right\}^w \cdot 
\dfrac{\prod_{i=0}^{r-1} \varGamma\left(w+\frac{i-1}{r}\right)}{\prod_{i=0}^{r-p-1} 
\varGamma\left(w+\frac{i-a}{r-p}\right)
\prod_{i=0}^{r-q-1} \varGamma\left(w+\frac{i-b}{r-q}\right)}, 
\end{split}
\]
%%%%%
which is multiplied by formula \eqref{eqn:gpf-g} to yield 
%%%%%%%%%%%%%%
\begin{align*}
C_g(w) \, g(w;\lambda) 
&= \varGamma((r-p-q)w+1-a-b) \cdot 
\dfrac{\prod_{i=1}^r \sin \pi(w+v_i^*)}{\prod_{i=0}^{r-1} 
\sin \pi\left(w+\frac{i-1}{r}\right)} \\
&\phantom{=} \times C_2 \cdot \tilde{d}^w \cdot 
\dfrac{\prod_{i=1}^r \varGamma(w+v_i^*)}{\prod_{i=0}^{r-p-1} 
\varGamma\left(w+\frac{i-a}{r-p}\right)
\prod_{i=0}^{r-q-1} \varGamma\left(w+\frac{i-b}{r-q}\right)}, 
\end{align*}
%%%%%%%%%%%%%%
where one uses definition \eqref{eqn:C2} and the relation 
$\delta \cdot r^r (r-p)^{p-r} (r-q)^{q-r} = \tilde{d}$ that 
follows from \eqref{eqn:D-d} and \eqref{eqn:d-t}. 
Taking equation \eqref{eqn:vi-s} into account we have          
%%%%%%%%%%%%%%%%%%%%%%%%%%%%%%% eqn:Cg %%%%%%%%%%%%%%%%%%%%%%%%%%%%%%%%%%%
\begin{equation} \label{eqn:Cg}
\begin{split}
C_g(w) \, g(w;\lambda) &= \varGamma((r-p-q)w+1-a-b) \cdot 
\dfrac{\prod_{i=1}^r \sin \pi(w+v_i^*)}{\prod_{i=0}^{r-1} 
\sin \pi\left(w+\frac{i-1}{r}\right)} \\
&\phantom{=} \times C_2 \cdot \tilde{d}^w \cdot 
\dfrac{\prod_{i=0}^{p-1} \varGamma\left(w+\frac{i+a}{p}\right)
\prod_{i=0}^{q-1} \varGamma\left(w+\frac{i+b}{q}\right)}{\prod_{i=0}^{r-1} 
\varGamma(w+v_i)}. 
\end{split}
\end{equation}
%%%%%%%%%%%%%%%%%%%%%%%%%%%%%%%%%%%%%%%%%%%%%%%%%%%%%%%%%%%%%%%%%%%%%%%%%% 
Substituting \eqref{eqn:Cf} and \eqref{eqn:Cg} into \eqref{eqn:cf} yields 
the desired formula \eqref{eqn:h2} with \eqref{eqn:d-t} and 
\eqref{eqn:chi}. \hfill $\Box$ 
%%%%%%%%%%%%%%%%%%%%%%%%%%%%%%% end proof %%%%%%%%%%%%%%%%%%%%%%%%%%%%%%%%%
%%%%%%%%%%%%%%%%%%%%%%%%%%%%%%% lem:chi %%%%%%%%%%%%%%%%%%%%%%%%%%%%%%%%%%%
\begin{lemma} \label{lem:chi} 
In Lemma $\ref{lem:r-cf2}$ the function $\chi(w)$ in formula \eqref{eqn:chi} 
must be constant so that there exists a real constant $C_3$ such that 
formula \eqref{eqn:h2} becomes 
%%%%%%%%%%%%%%%%%%%%%%%%% eqn:h3 %%%%%%%%%%%%%%%%%%%%%%%%%%%%%%%%%%%%%%
\begin{equation} \label{eqn:h3}
\dfrac{h(w;\lambda)}{\varGamma((r-p-q)w+1-a-b)} = C_3 \cdot \tilde{d}^w    
\cdot \dfrac{\prod_{i=0}^{p-1}\varGamma\left(w+\frac{i+a}{p}\right) 
\prod_{i=0}^{q-1}\varGamma\left(w+\frac{i+b}{q}\right)}{\prod_{i=1}^{r} 
\varGamma\left(w+v_i\right)}.   
\end{equation}
%%%%%%%%%%%%%%%%%%%%%%%%%%%%%%%%%%%%%%%%%%%%%%%%%%%%%%%%%%%%%%%%%%%%%%% 
\end{lemma}
%%%%%%%%%%%%%%%%%%%%%%%%%%%%%%%%%%%%%%%%%%%%%%%%%%%%%%%%%%%%%%%%%%%%%%%%%%
%%%%%%%%%%%%%%%%%%%%%%%%%%%%%%% begin proof %%%%%%%%%%%%%%%%%%%%%%%%%%%%%%
{\it Proof}. 
The left-hand side of equation \eqref{eqn:h3}, which is denoted by  
$\tilde{h}(w;\lambda)$, is an entire holomorphic function of $w$, 
since by formula \eqref{eqn:h} any pole of $h(w;\lambda)$ must be 
simple and located where $(r-p-q)w+1-a-b$ becomes a nonpositive 
integer, so that it is canceled by a zero of 
$1/\varGamma((r-p-q)w+1-a-b)$. 
It follows from formula \eqref{eqn:h2} that 
%%%%%%%%%%%%%%%%%
\begin{equation*}
\chi(w) = \tilde{d}^{-w} \cdot \tilde{h}(w;\lambda)      
\cdot \dfrac{\prod_{i=1}^{r} \varGamma\left(w+v_i\right)}{\prod_{i=0}^{p-1}
\varGamma\left(w+\frac{i+a}{p}\right) 
\prod_{i=0}^{q-1}\varGamma\left(w+\frac{i+b}{q}\right)},    
\end{equation*}
%%%%%%%%%%%%%%%%
and so $\chi(w)$ is holomorphic in the half-plane $\rRe \, w > 
- \min\{ v_i \,:\, i = 1, \dots, n \}$. 
On the other hand formula \eqref{eqn:chi} implies that 
$\chi(w)$ is a periodic function of period one, since  
$p$, $q$, and $r$ are positive integers with $r-p-q$ even. 
Thus $\chi(w)$ must be an entire holomorphic and periodic function   
of period one. 
In particular $\chi(w)$ is bounded on the horizontal strip 
$|\rIm \, w| \le 1$.  
%%%%%
\par
%%%%%
Now notice that $\sin z$ admits a two-sided bound 
$\frac{1}{4} \, e^{|\rIm \,z|} \le |\sin z| \le e^{|\rIm \, z|}$ on 
the outer region $|\rIm \, z| \ge 1$.  
Applying it to formula \eqref{eqn:chi} yields an estimate:  
%%%%%%
\[
|\chi(w)| \le |C_1| \, 
\dfrac{e^{\pi p |\rIm \, w|} \cdot e^{\pi q |\rIm \, w|}}{\frac{1}{4} \, 
e^{\pi r |\rIm \, w|}} + |C_2| \, \dfrac{ \prod_{i=1}^r 
e^{\pi |\rIm \, w|}}{\prod_{i=0}^{r-1} \frac{1}{4} \, e^{\pi |\rIm \, w|}} 
= 4|C_1| \, e^{-\pi(r-p-q)|\rIm \, w|} + 4^r |C_2|,  
\]
%%%%%% 
and hence $|\chi(w)| \le 4|C_1| \, e^{-\pi(r-p-q)} + 4^r |C_2|$ on 
$|\rIm \, w| \ge 1$. 
Thus the entire function $\chi(w)$ is bounded on $\C_w$. 
By Liouville's theorem $\chi(w)$ must be a constant $C_3$, 
which is clearly real. \hfill $\Box$ \par\medskip
%%%%%%%%%%%%%%%%%%%%%%%%%%%%%%% end proof %%%%%%%%%%%%%%%%%%%%%%%%%%%%%%%%
Given a positive integer $k$, we put $[k] := \{0,1,\dots, k-1\}$. 
By division relation \eqref{eqn:division} there exist subsets 
$I_p \subset [p]$, $I_q \subset [q]$, $J_p \subset [r-p]$, 
$J_q \subset [r-q]$, with $0 \not \in I_p$ and $0 \not \in I_q$,  
such that 
%%%%%%%%%%%%%%%%%%%%%%%%%%%%%%% eqn:IJ %%%%%%%%%%%%%%%%%%%%%%%%%%%%%%%%%%%
\begin{equation} \label{eqn:IJ}
\prod_{i=1}^r (w+v_i) = 
\prod_{i \in I_p} \left(w+{\ts\frac{i+a}{p}}\right) 
\prod_{i \in I_q} \left(w+{\ts\frac{i+b}{q}}\right) 
\prod_{j \in J_p} \left(w+{\ts\frac{j-a}{r-p}}\right) 
\prod_{j \in J_q} \left(w+{\ts\frac{j-b}{r-q}}\right).   
\end{equation}
%%%%%%%%%%%%%%%%%%%%%%%%%%%%%%%%%%%%%%%%%%%%%%%%%%%%%%%%%%%%%%%%%%%%%%%%%%%
If we put $\bar{I}_p := [p] \setminus I_p$ and $\bar{I}_q := [q] \setminus I_q$, 
then equation \eqref{eqn:h3} becomes  
%%%%%%%%%%%%%%%%%%%%%%%%%%%% eqn:h4 %%%%%%%%%%%%%%%%%%%%%%%%%%%%%%%%%%%%%%%
\begin{equation} \label{eqn:h4}
\dfrac{h(w;\lambda)}{\varGamma((r-p-q)w+1-a-b)} = C_3 \cdot \tilde{d}^w    
\cdot \dfrac{\prod_{i\in \bar{I}_p} \varGamma\left(w+\frac{i+a}{p}\right) 
\prod_{i \in \bar{I}_q} \varGamma\left(w+\frac{i+b}{q}\right)}{
\prod_{j \in J_p} \varGamma\left(w+\frac{j-a}{r-p}\right) 
\prod_{j \in J_q} \varGamma\left(w+\frac{j-b}{r-q}\right)}.   
\end{equation}
%%%%%%%%%%%%%%%%%%%%%%%%%%%%%%%%%%%%%%%%%%%%%%%%%%%%%%%%%%%%%%%%%%%%%%%
\par
%%%%
To exploit formula \eqref{eqn:h4} we need a preliminary lemma. 
A {\sl multi-set} is a set allowing repeated elements. 
For multi-sets $S = \{s_1, \dots, s_m\}$ and $T = \{t_1, \dots, t_n\}$, 
we write $S \succ T$ if $m \le n$ and there exists a re-indexing of 
$t_1, \dots, t_n$ such that $s_i - t_i \in \Z_{\ge 0}$ for every  
$i = 1, \dots, m$. 
%%%%%%%%%%%%%%%%%%%%%%%%%%%%%%% lem:nni %%%%%%%%%%%%%%%%%%%%%%%%%%%%%%%%%%
\begin{lemma} \label{lem:nni}
Let $S = \{s_1, \dots, s_m\}$ and $T = \{t_1, \dots, t_n\}$ be 
multi-sets of real numbers. 
If  
%%%%%
\[ 
h(w) = \dfrac{\varGamma(w+s_1) \cdots \varGamma(w+s_m)}{\varGamma(w+t_1) \cdots \varGamma(w+t_n)}
\]
%%%%%
is an entire function of $w$, then $S \succ T$.  
\end{lemma} 
%%%%%%%%%%%%%%%%%%%%%%%% begin proof %%%%%%%%%%%%%%%%%%%%%%%%%%%%%%%%%%%
{\it Proof}. 
The proof is by induction on $m$. 
When $m = 0$ the assertion is obvious as the numerator of $h(w)$ is $1$. 
Let $m \ge 1$. 
We may assume $s_1 \le \cdots \le s_m$.   
An upper factor $\varGamma(w+s_1)$ of $h(w)$ has a pole at $w = -s_1$. 
In order for $h(w)$ to be holomorphic, a lower factor 
$\varGamma(w+t_j)$ must have a pole at the same point.    
After a transposition of $t_j$ and $t_1$ we may put $j = 1$. 
Then $-s_1+t_1=-r_1 \in \mathbb{Z}_{\le 0}$, that is, 
$r_1 \in \mathbb{Z}_{\ge0}$. 
Since $\varGamma(w+s_1)/\varGamma(w+t_1) = (w+t_1)_{r_1}$, we have 
%%%%%%
\[
h(w) = (w+t_1)_{r_1} \cdot h_1(w), \qquad \mbox{where} \qquad 
h_1(w) := \dfrac{\varGamma(w+s_2) \cdots \varGamma(w+s_m)}{\varGamma(w+t_2) \cdots \varGamma(w+t_n)}. 
\]
%%%%%% 
We claim that $h_1(w)$ is entire holomorphic. 
If $r_1 = 0$ this is obvious since $h(w) = h_1(w)$. 
Let $r_1 \ge 1$. 
Any pole of $\varGamma(w+s_2) \cdots \varGamma(w+s_m)$ is separated 
from all the roots of $(w+t_1)_{r_1}$, because we have 
$-s_m \le \cdots \le -s_2 \le -s_1$ and the roots of $(w+t_1)_{r_1}$ 
are located at $-s_1+1 < -s_1+2 < \dots < -t_1$ in an increasing order.    
Thus any pole of $h_1(w)$, if it exists, cannot be 
canceled by a root of $(w+t_1)_{r_1}$.  
Accordingly, $h_1(w)$ has no poles, since $h(w)$ has no poles.  
Now we can apply the induction hypothesis to $h_1(w)$ to conclude the proof.      
\hfill $\Box$ \par\medskip
%%%%%%%%%%%%%%%%%%%%%%%% end proof %%%%%%%%%%%%%%%%%%%%%%%%%%%%%%%%%%%%%
Since the left-hand side of equation \eqref{eqn:h4} is entire holomorphic 
in $w$, so must be the gamma products on the right. 
Thus Lemma \ref{lem:nni} yields   
%%%%%%%%%%%%%%%%%%%%%%%% eqn:succ1 %%%%%%%%%%%%%%%%%%%%%%%%%%%%%%%%%%%%%
\begin{equation} \label{eqn:succ1}
\left\{ \frac{i+a}{p}\right\}_{i\in \bar{I}_p} \bigcup \,\,  
\left\{\frac{i+b}{q}\right\}_{i \in \bar{I}_q}
\succ \, 
\left\{\frac{j-a}{r-p}\right\}_{j \in J_p} \bigcup \,\, 
\left\{\frac{j-b}{r-q}\right\}_{j \in J_q},   
\end{equation}
%%%%%%%%%%%%%%%%%%%%%%%%%%%%%%%%%%%%%%%%%%%%%%%%%%%%%%%%%%%%%%%%%%%%%%%%
where the both sides above are thought of as multi-sets. 
Note that $0 \in \bar{I}_p$ and $0 \in \bar{I}_q$.  
%%%%%%%%%%%%%%%%%%%%%%%%%%%%% tab:lem-ab %%%%%%%%%%%%%%%%%%%%%%%%%%%%%%%
\begin{table}[t]
\begin{center}
\begin{tabular}{|c|c|c|c|} 
\hline 
     &      &     &     \\[-4mm]
case & type & $a$ & $b$ \\[1mm]
\hline 
     &      &     &     \\[-3mm]
$1$ & $(\rI, \, \rI)$   & $\dfrac{p i}{r}$                         & $\dfrac{q j}{r}$ \\[3mm]
\hline
     &      &     &     \\[-3mm]
$2$ & $(\rII, \, \rII)$ & \hspace{2mm} $\dfrac{p \{ (r-p)i - q j \}}{r(r-p-q)}$ \hspace{2mm} 
& \hspace{2mm} $\dfrac{q \{ (r-q)j - p i \}}{r(r-p-q)}$ \hspace{2mm} \\[5mm] 
\hline 
     &      &     &     \\[-3mm]
$3$ & $(\rI, \, \rII)$  & $\dfrac{p i}{r}$                         & $\dfrac{q (r j - p i)}{r(r-p)}$ \\[4mm]
\hline
     &      &     &     \\[-3mm]
$4$ & $(\rII, \, \rI)$  & $\dfrac{p (r i -q j)}{r(r-q)}$           & $\dfrac{q j}{r}$ \\[4mm] 
\hline
\end{tabular}   
\end{center}
\caption{Formula for $(a, b)$ in terms of $(p,q,r)$ and $(i, j) \in \Z_{\ge 0}^2$.}
\label{tab:lem-ab}
\end{table}
%%%%%%%%%%%%%%%%%%%%%%%%%%%%%%%%%%%%%%%%%%%%%%%%%%%%%%%%%%%%%%%%%%%%%%%%%
%%%%%%%%%%%%%%%%%%%%%%%%%%%%% lem:ab %%%%%%%%%%%%%%%%%%%%%%%%%%%%%%%%%%%%
\begin{lemma} \label{lem:ab}
For any $(\rA)$-solution $\lambda = (p,q,r;a,b;x) \in \cD^-$ there exists 
a pair $(i, j)$ of nonnegative integers such that one of the four cases 
in Table $\ref{tab:lem-ab}$ occurs.   
In particular, in any case $a$ and $b$ must be rational numbers.       
\end{lemma}
%%%%%%%%%%%%%%%%%%%%%%%%%%%%%%%%%%%%%%%%%%%%%%%%%%%%%%%%%%%%%%%%%%%%%%%%
%%%%%%%%%%%%%%%%%%%%%%%%%%%%% begin proof %%%%%%%%%%%%%%%%%%%%%%%%%%%%%%
{\it Proof}.  
Since $0 \in \bar{I}_p$, condition \eqref{eqn:succ1} implies that either 
$(\rI)_p$ or $(\rII)_p$ below holds:   
%%%%%%%%
\begin{enumerate}
\item[$(\rI_p)$] there exists an integer $i \in J_p + (r-p) \, \Z_{\ge0}$ such that  
$\dfrac{a}{p} = \dfrac{i-a}{r-p}$, 
\item[$(\rII_p)$] there exists an integer $i \in J_q + (r-q) \, \Z_{\ge0}$ such that 
$\dfrac{a}{p} = \dfrac{i-b}{r-q}$.   
\end{enumerate}
%%%%%%%% 
In a similar manner, since $0 \in \bar{I}_q$, condition \eqref{eqn:succ1} 
implies that either $(\rI_q)$ or $(\rII_q)$ below holds:    
%%%%%%%%
\begin{enumerate}
\item[$(\rI_q)$] there exists an integer $j \in J_q + (r-q) \, \Z_{\ge0}$ such that  
$\dfrac{b}{q} = \dfrac{j-b}{r-q}$, 
\item[$(\rII_q)$] there exists an integer $j \in J_p + (r-p) \, \Z_{\ge0}$ such that 
$\dfrac{b}{q} = \dfrac{j-a}{r-p}$.   
\end{enumerate}
%%%%%%%% 
There are a total of four types $(\rI_p, \, \rI_q)$, $(\rII_p, \, \rII_q)$, 
$(\rI_p, \, \rII_q)$, $(\rII_p, \, \rI_q)$, which are exactly the Cases 
$1$--$4$ of Table \ref{tab:lem-ab} respectively, where suffixes 
$* \in \{p, \, q\}$ of $\rI_*$ and $\rII_*$ are omitted. 
In each case we have a pair of linear equations for $(a, b)$, 
which can be solved as indicated in Table \ref{tab:lem-ab}.  
\hfill $\Box$
%%%%%%%%%%%%%%%%%%%%%%%%%%%%% end proof %%%%%%%%%%%%%%%%%%%%%%%%%%%%%%%%
%%%%%%%%%%%%%%%%%%%%%%%%%%%%% lem:ab-J %%%%%%%%%%%%%%%%%%%%%%%%%%%%%%%%%
\begin{lemma} \label{lem:ab-J} 
If $a > 0$ then $0 \not\in J_p$. 
Similarly if $b > 0$ then $0 \not\in J_q$.  
\end{lemma}
%%%%%%%%%%%%%%%%%%%%%%%%%%%%%%%%%%%%%%%%%%%%%%%%%%%%%%%%%%%%%%%%%%%%%%%%
%%%%%%%%%%%%%%%%%%%%%%%%%%%%% begin proof %%%%%%%%%%%%%%%%%%%%%%%%%%%%%%
{\it Proof}. 
We use assertion (1) of Iwasaki \cite[Proposition 11.3]{Iwasaki} with  
$k = 0$.  
In the notation there it states that for each $k = 0, \dots, r-p-1$,  
one has $(w-w_k^*)|\prod_{i=1}^r(w+v_i)$ if and only if 
%%%%%%%%%%%%%%%%%%%%%%%%%%%%% eqn:ab-J1 %%%%%%%%%%%%%%%%%%%%%%%%%%%%%%%%
\begin{equation} \label{eqn:ab-J1}
(\gamma_k^*+k)_{p+1} \cdot \cF_k(\beta_k^*; \gamma_k^*;x) = 0, 
\end{equation}
%%%%%%%%%%%%%%%%%%%%%%%%%%%%%%%%%%%%%%%%%%%%%%%%%%%%%%%%%%%%%%%%%%%%%%%%
unless \cite[condition (101a)]{Iwasaki} is satisfied, that is, unless  
%%%%%%%%%%%%%%%%%%%%%%%%%%%% eqn:ab-J2 %%%%%%%%%%%%%%%%%%%%%%%%%%%%%%%%%
\begin{equation} \label{eqn:ab-J2}
\tilde{\beta}_k^*, \,\, \tilde{\gamma}_k^* \in \Z, \qquad 
0 \le - \tilde{\beta}_k^* \le - \tilde{\gamma}_k^* \le r-p-k-2.  
\end{equation}
%%%%%%%%%%%%%%%%%%%%%%%%%%%%%%%%%%%%%%%%%%%%%%%%%%%%%%%%%%%%%%%%%%%%%%%
When $k = 0$, definitions in 
\cite[formula (96) and Proposition 11.1]{Iwasaki} read 
%%%%%
\[
w_0^* = \dfrac{a}{r-p}, \qquad  
\gamma_0^* = r w_0^* =  \dfrac{r a}{r-p}, \qquad  
\tilde{\gamma}_0^* = 2-r(w_0^* +1) = 2-r-\dfrac{r a}{r-p},  
\]
%%%%% 
so condition $- \tilde{\gamma}_k^* \le r-p-k-2$ in \eqref{eqn:ab-J2} 
is equivalent to $a \le -p(r-p)/r \,\, (< 0)$. 
Thus if $a > 0$ then condition \eqref{eqn:ab-J2} with $k = 0$ is not 
fulfilled.  
When $k = 0$, condition \eqref{eqn:ab-J1} becomes 
$(\gamma_0^*)_{p+1} = 0$, since $\cF_0(\beta_0^*; \gamma_0^*;x) = 1$ 
by \cite[definition (95)]{Iwasaki}. 
But if $a > 0$ then $\gamma_0^* > 0$ and so $(\gamma_0^*)_{p+1} > 0$, 
which means that $w-w_0^* = w + \frac{0-a}{r-p}$ is {\sl not} a factor 
of $\prod_{i=1}^r(w+v_i)$. 
This in turn implies $0 \not\in J_p$ by formula $\eqref{eqn:IJ}$. 
The implication $b > 0 \Rightarrow 0 \not\in J_q$ is proved in 
the same way. \hfill $\Box$ 
%%%%%%%%%%%%%%%%%%%%%%%%%%%%% end proof %%%%%%%%%%%%%%%%%%%%%%%%%%%%%%%%
%%%%%%%%%%%%%%%%%%%%%%%%%%%%% sec:dr-fc %%%%%%%%%%%%%%%%%%%%%%%%%%%%%%%%
\section{Duality Revisited and Finite Cardinality} \label{sec:dr-fc}  
%%%%%%%%%%%%%%%%%%%%%%%%%%%%%%%%%%%%%%%%%%%%%%%%%%%%%%%%%%%%%%%%%%%%%%%%
Consider an $(\rA)$-solution $\lambda = (p,q,r;a,b;x) \in \cD^-$  
and its dual solution $\lambda' = (p,q,r;a',b';x) \in \cD^-$.  
To the pair $(\lambda, \lambda')$ we associate the three-by-two matrix 
%%%%%%%%
\[
\begin{pmatrix}
p  & q \\ 
a  & b \\
a' & b' 
\end{pmatrix}.  
\]
%%%%%%%%
Exchange of its columns represents the trivial symmetry \cite[(8a)]{Iwasaki} 
applied to $(\lambda, \lambda')$, whereas exchange of its middle and 
bottom rows represents duality \eqref{eqn:duality} itself. 
With the omission of its top row $(p, q)$, the matrix above 
is abbreviated to the square matrix  
%%%%%%%%%%%%%%%%%%%%%%%%%%%%% eqn:matrix %%%%%%%%%%%%%%%%%%%%%%%%%%%%%%%%
\begin{equation} \label{eqn:matrix}
\begin{pmatrix} 
a  & b \\
a' & b' 
\end{pmatrix}.    
\end{equation}
%%%%%%%%%%%%%%%%%%%%%%%%%%%%%%%%%%%%%%%%%%%%%%%%%%%%%%%%%%%%%%%%%%%%%%%%%
Two dual pairs are said to be {\sl equivalent} if their matrices are 
transitive via the column or row exchange or the composition of them.      
The matrix \eqref{eqn:matrix} is said to be of type 
%%%%%%%%%%%%%%%%%%%%%%%%%%%%% eqn:type %%%%%%%%%%%%%%%%%%%%%%%%%%%%%%%%
\begin{equation} \label{eqn:type}
\begin{pmatrix} 
J_1 & J_2 \\
J_3 & J_4
\end{pmatrix}, \qquad J_1, J_2, J_3, J_4 \in \{\rI, \, \rII\}, 
\end{equation}
%%%%%%%%%%%%%%%%%%%%%%%%%%%%%%%%%%%%%%%%%%%%%%%%%%%%%%%%%%%%%%%%%%%%%%%
if the top row $(a, b)$ of matrix \eqref{eqn:matrix} is of type 
$(J_1, J_2)$, while its bottom row $(a', b')$ is of type $(J_3, J_4)$ 
in the sense of Table \ref{tab:lem-ab}. 
%%%%%%%%%%%%%%%%%%%%%%%%%% tab:type %%%%%%%%%%%%%%%%%%%%%%%%%%%%%%%%%%%
\begin{table}[t]
\begin{center}
\begin{tabular}{|c|c|c|c|c|c|c|c|}
\hline  
     &   &   &   &   &   &   &   \\[-3mm]
case & 1 & 2 & 3 & 4 & 5 & 6 & 7 \\[1mm]
\hline 
     &   &   &   &   &   &   &   \\[-3mm] 
type & 
$\begin{matrix} \rI & \rI \\[1mm] \rI & \rI \end{matrix}$ &  
$\begin{matrix} \rII & \rII \\[2mm] \rII & \rII \end{matrix}$ &
$\begin{matrix} \rI & \rII \\[2mm] \rI & \rII \end{matrix}$ &   
$\begin{matrix} \rI & \rI \\[2mm] \rI & \rII \end{matrix}$ & 
$\begin{matrix} \rII & \rII \\[2mm] \rII & \rI \end{matrix}$ & 
$\begin{matrix} \rI & \rII \\[2mm] \rII & \rI \end{matrix}$ &
$\begin{matrix} \rI & \rI \\[2mm] \rII & \rII \end{matrix}$ \\[4mm]
\hline 
\end{tabular}
\end{center}
\caption{Types of matrix \eqref{eqn:matrix} up to equivalence.}
\label{tab:type}
\end{table}
%%%%%%%%%%%%%%%%%%%%%%%%%%%%%%%%%%%%%%%%%%%%%%%%%%%%%%%%%%%%%%%%%%%%%%%%
There are a total of seven types up to equivalence, which are tabulated 
in Table \ref{tab:type}, where parentheses in \eqref{eqn:type} are omitted. 
%%%%% 
\par
%%%%%
The essential parts of definition \eqref{eqn:duality} for duality can be 
rewritten as  
%%%%%%%%%%%%%%%%%%%%%%%%%%%%% eqn:duality2 %%%%%%%%%%%%%%%%%%%%%%%%%%%%%
\begin{equation} \label{eqn:duality2}
a + a' = 1- \dfrac{2 p}{r}, \qquad b + b' = 1- \dfrac{2 q}{r}. 
\end{equation}
%%%%%%%%%%%%%%%%%%%%%%%%%%%%%%%%%%%%%%%%%%%%%%%%%%%%%%%%%%%%%%%%%%%%%%%%
In each case of Table \ref{tab:type} we shall see what kinds of 
consequences are derived from equations \eqref{eqn:duality2}. 
In this section we employ the following notation. 
%%%%%%%%%%%%%%%%%%%%%%%%%% tab:ab %%%%%%%%%%%%%%%%%%%%%%%%%%%%%%%%%%%%%%%%
\begin{table}[p]
\begin{center}
\rotatebox{90}{
\begin{tabular}{|c|l|l|l|l|c|}
\hline
     &     &     &     &     &     \\[-3mm]
     & (1) & (2) & (3) & (4) & (5) \\[1mm]
\hline 
     &     &     &     &     &     \\[-3mm]
case & {\small division relations} & {\small formula for $(a, b)$} & {\small formula for $(a', b')$} & 
$\Z$-linear equations for $(i, j; i', j')$ & type \\[2mm]
\hline 
     &     &     &     &     &     \\[-3mm]
$1$ & 
$\begin{array}{l} p|r \\[1mm] q|r \end{array}$ & 
$\begin{array}{l} a = \frac{i}{r_p} \\[2mm] b = \frac{j}{r_q} \end{array}$ & 
$\begin{array}{l} a'= \frac{i'}{r_p} \\[2mm] b'= \frac{j'}{r_q} \end{array}$ & 
$\begin{array}{l} i+i'= r_p-2 \\[1mm] j+j' = r_q-2 \end{array}$ & 
$\begin{matrix} \rI & \rI \\[1mm] \rI & \rI \end{matrix}$ \\[5mm]
\hline  
     &     &     &     &     &     \\[-3mm] 
$2$ & 
$\begin{array}{l} p|(r-p-q) \\[2mm] q|(r-p-q) \end{array}$ & 
$\begin{array}{l} a  = \frac{(r-p)i-q j}{r(r-p-q)_p} \\[2mm] b  = \frac{(r-q)j-p i}{r(r-p-q)_q} \end{array}$ & 
$\begin{array}{l} a' = \frac{(r-p)i'-q j'}{r(r-p-q)_p} \\[2mm] b' = \frac{(r-q)j'-p i'}{r(r-p-q)_q} \end{array}$ & 
$\begin{array}{l} i+i' = (r-p-q)_p \\[2mm] j+j' = (r-p-q)_q \end{array}$ & 
$\begin{matrix} \rII & \rII \\[2mm] \rII & \rII \end{matrix}$ \\[6mm]
\hline 
     &     &     &     &     &     \\[-3mm]
$3$ & 
$\begin{array}{l} p|r \\[2mm] q|(r-p-q) \end{array}$ & 
$\begin{array}{l} a  = \frac{i}{r_p} \\[2mm] b  = \frac{r_p j-i}{r_p (r-p)_q} \end{array}$ &
$\begin{array}{l} a' = \frac{i'}{r_p} \\[2mm] b' = \frac{r_p j'-i'}{r_p (r-p)_q} \end{array}$ &
$\begin{array}{l} i+i' = r_p-2 \\[2mm] j+j' = (r-p-q)_q \end{array}$ & 
$\begin{matrix} \rI & \rII \\[2mm] \rI & \rII \end{matrix}$ \\[6mm]
\hline 
     &     &     &     &     &     \\[-3mm]    
$4$ & 
$\begin{array}{l} p|r \\[2mm] q|r_p(r-p-q) \end{array}$  & 
$\begin{array}{l} a = \frac{i}{r_p} \\[2mm] b  = \frac{q j}{r} \end{array}$ & 
$\begin{array}{l} a' = \frac{i'}{r_p} \\[2mm] b' = \frac{r_p j'-i'}{(r_p(r-p))_q} \end{array}$ & 
$\begin{array}{l} i+i' = r_p-2 \\[2mm] {\scriptstyle i+(r_p-1)j+r_p j' = (r_p(r-p-q))_q} 
\end{array}$ & 
$\begin{matrix} \rI & \rI \\[2mm] \rI & \rII \end{matrix}$ \\[6mm]
\hline   
     &     &     &     &     &     \\[-3mm]    
$5$ & 
$\begin{array}{l} p|(r-p-q) \\[2mm] q|r(r-p-q)_p \end{array}$  & 
$\begin{array}{l} a = \frac{(r-p)i-q j}{r(r-p-q)_p} \\[2mm] b  = \frac{(r-q)_p j-i}{(r(r-p-q)_p)_q} \end{array}$ & 
$\begin{array}{l} a' = \frac{r i'-q j'}{r(r-q)_p} \\[2mm] b' = \frac{q j'}{r} \end{array}$ & 
$\begin{array}{l} i+i' = (r-p-q)_p \\[2mm] {\scriptstyle i'+(r-q)_p j+ (r-p-q)_p j' = ((r-q)(r-p-q)_p)_q} \end{array}$ & 
$\begin{matrix} \rII & \rII \\[2mm] \rII & \rI \end{matrix}$ \\[6mm] 
\hline   
     &     &     &     &     &     \\[-3mm] 
$6$ & 
$\begin{array}{l} p|r(r-p-q) \\[2mm] q|r(r-p-q) \end{array}$ & 
$\begin{array}{l} a  = \frac{p i}{r} \\[2mm] b  = \frac{r j-p i}{(r(r-p))_q} \end{array}$ &
$\begin{array}{l} a' = \frac{r i'-q j'}{(r(r-q))_p} \\[2mm] b' = \frac{q j'}{r} \end{array}$ &
$\begin{array}{l} {\scriptstyle (r-p-q)i+(r-p)i'+q j = ((r-p)(r-p-q))_p} \\[2mm] 
{\scriptstyle (r-p-q)j'+p i' + (r-q)j = ((r-q)(r-p-q))_q} \end{array}$ & 
$\begin{matrix} \rI & \rII \\[2mm] \rII & \rI \end{matrix}$ \\[6mm]
\hline 
\end{tabular}
} 
\end{center}
\caption{Candidates for dual pairs of $(\rA)$-solutions in $\cD^-$.}
\label{tab:ab}
\end{table}
%%%%%%%%%%%%%%%%%%%%%%%%%%%%%%%%%%%%%%%%%%%%%%%%%%%%%%%%%%%%%%%%%%%%%%%%   
For positive integers $s$ and $t$ we write 
%%%%%%%%%%%%%%%%%%%%%%%%%%%%% eqn:s_t %%%%%%%%%%%%%%%%%%%%%%%%%%%%%%%%%%
\begin{equation} \label{eqn:s_t}
s_t := \frac{s}{t} \qquad \mbox{when and only when} \quad t|s.
\end{equation}
%%%%%%%%%%%%%%%%%%%%%%%%%%%%%%%%%%%%%%%%%%%%%%%%%%%%%%%%%%%%%%%%%%%%%%%% 
It is a convenient notation which indicates at once that $s_t$ stands 
for a positive integer; it also saves space when $s$ is a large 
expression, but for example $(r-p-q)_p$ should not be confused with 
a rising  factorial number. 
In this section $(\cdots)_p$ never represents a factorial number.    
%%%%%%%%%%%%%%%%%%%%%%%%%%%%% lem:case1-6 %%%%%%%%%%%%%%%%%%%%%%%%%%%%%%
\begin{lemma} \label{lem:case1-6} 
In cases $1$--$6$ of Table $\ref{tab:type}$, equations 
\eqref{eqn:duality2} lead to the conditions in Table $\ref{tab:ab}$, 
which consist of the following four items:  
\begin{enumerate} 
\item two division relations for the integer triple $\bp = (p,q;r)$,   
\item a formula for $(a,b)$ in terms of $\bp$ and a pair  
of nonnegative integers $(i,j) \in \Z_{\ge 0}^2$,    
\item a formula for $(a',b')$ in terms of $\bp$ and a pair 
of nonnegative integers  $(i',j') \in \Z_{\ge 0}^2$,      
\item two $\Z$-linear equations for quadruple $(i,j;i',j') \in \Z_{\ge 0}^4$, 
\end{enumerate}
where for each $\nu \in \{1,2,3,4\}$ item $(\nu)$ is exhibited in 
column $(\nu)$ of Table $\ref{tab:ab}$. 
\end{lemma}
%%%%%%%%%%%%%%%%%%%%%%%%%%%%%%%%%%%%%%%%%%%%%%%%%%%%%%%%%%%%%%%%%%%%%%%%
%%%%%%%%%%%%%%%%%%%%%%%%%%%%% begin proof %%%%%%%%%%%%%%%%%%%%%%%%%%%%%%
{\it Proof}. This lemma is proved by case-by-case treatments presented 
below. \\[1mm] 
%%%%%%%%%%%%%%%%%%%%%%%%%%%%% case 1 %%%%%%%%%%%%%%%%%%%%%%%%%%%%%%%%%%%
{\bf Case 1}. By Lemma \ref{lem:ab} there exists a quadruple of 
nonnegative integers $(i,j;i',j')$ such that 
%%%%%%%%%%%%%%%%%%%%%%%%%%%%% eqn:c1-1 %%%%%%%%%%%%%%%%%%%%%%%%%%%%%%%%%
\begin{equation} \label{eqn:c1-1} 
a = \dfrac{p i}{r}, \qquad b = \dfrac{q j}{r}, \qquad 
a' = \dfrac{p i'}{r}, \qquad b' =  \dfrac{q j'}{r}. 
\end{equation}
%%%%%%%%%%%%%%%%%%%%%%%%%%%%%%%%%%%%%%%%%%%%%%%%%%%%%%%%%%%%%%%%%%%%%%%%
Substituting formula \eqref{eqn:c1-1} into equations \eqref{eqn:duality2} 
yields $i + i' = \frac{r}{p} -2 \in \Z$ and $j + j' = \frac{r}{q} -2 \in \Z$, 
which implies $p|r$ and $q|r$, so that we have 
%%%%%%%%%%%%%%%%%%%%%%%%%%%%% eqn:ls-c1 %%%%%%%%%%%%%%%%%%%%%%%%%%%%%%%%
\begin{equation} \label{eqn:ls-c1}
i + i' = r_p -2, \qquad  j + j' = r_q -2, \qquad 
i, j, i', j' \in \Z_{\ge 0}. 
\end{equation}
%%%%%%%%%%%%%%%%%%%%%%%%%%%%%%%%%%%%%%%%%%%%%%%%%%%%%%%%%%%%%%%%%%%%%%%% 
Formula \eqref{eqn:c1-1} then becomes 
%%%%%%%%%%%%%%%%%%%%%%%%%%%%% eqn:c1-2 %%%%%%%%%%%%%%%%%%%%%%%%%%%%%%%%%
\begin{equation} \label{eqn:c1-2} 
a = \dfrac{i}{r_p}, \qquad b = \dfrac{j}{r_q}, \qquad 
a' = \dfrac{i'}{r_p}, \qquad b' =  \dfrac{j'}{r_q}. 
\end{equation}
%%%%%%%%%%%%%%%%%%%%%%%%%%%%%%%%%%%%%%%%%%%%%%%%%%%%%%%%%%%%%%%%%%%%%%%%
Thus all the conditions in case 1 of Table \ref{tab:ab} have been 
obtained. \\[1mm]
%%%%%%%%%%%%%%%%%%%%%%%%%%%%% case 2 %%%%%%%%%%%%%%%%%%%%%%%%%%%%%%%%%%% 
{\bf Case 2}. By Lemma \ref{lem:ab} there exists a quadruple of 
nonnegative integers $(i,j;i',j')$ such that 
%%%%%%%%%%%%%%%%%%%%%%%%%%%% eqn:c2-1 %%%%%%%%%%%%%%%%%%%%%%%%%%%%%%%%%%
\begin{equation} \label{eqn:c2-1}
\begin{array}{rlrl}
a &= \dfrac{p \{(r-p)i-q j\}}{r(r-p-q)}, \qquad &  
b &= \dfrac{q \{(r-q)j-p i\}}{r(r-p-q)}, \\[4mm] 
a' &= \dfrac{p \{(r-p)i'-q j'\}}{r(r-p-q)}, \qquad &  
b' &= \dfrac{q \{(r-q)j'-p i'\}}{r(r-p-q)}. 
\end{array}
\end{equation}  
%%%%%%%%%%%%%%%%%%%%%%%%%%%%%%%%%%%%%%%%%%%%%%%%%%%%%%%%%%%%%%%%%%%%%%%%
Substituting formula \eqref{eqn:c2-1} into equations \eqref{eqn:duality2} 
yields 
%%%%%
\[
\resizebox{0.9\hsize}{!}{$\displaystyle
(r-p) (i+i') - q(j+j') = (r-p-q) \left(\textstyle \frac{r}{p} - 2\right), 
\qquad 
-p (i+i') +(r-q)(j+j') = (r-p-q) \left(\textstyle \frac{r}{q} - 2\right). 
$} 
\]
%%%%% 
They are solved with respect to $i+i'$ and $j+j'$ to obtain 
$i+i' = (r-p-q)/p \in \Z$ and $j+j' = (r-p-q)/q \in \Z$, which imply 
$p|(r-p-q)$ and $q|(r-p-q)$, so that we have 
%%%%%%%%%%%%%%%%%%%%%%%%%%%% eqn:ls-c2 %%%%%%%%%%%%%%%%%%%%%%%%%%%%%%%%%
\begin{equation} \label{eqn:ls-c2}
i+i' = (r-p-q)_p, \qquad j+j' = (r-p-q)_q, \qquad 
i, j, i', j' \in \Z_{\ge 0}. 
\end{equation}
%%%%%%%%%%%%%%%%%%%%%%%%%%%%%%%%%%%%%%%%%%%%%%%%%%%%%%%%%%%%%%%%%%%%%%%% 
Formula \eqref{eqn:c2-1} then becomes 
%%%%%%%%%%%%%%%%%%%%%%%%%%%% eqn:c2-2 %%%%%%%%%%%%%%%%%%%%%%%%%%%%%%%%%%
\begin{equation} \label{eqn:c2-2}  
a = \dfrac{(r-p)i-q j}{r(r-p-q)_p}, \quad   
b = \dfrac{(r-q)j-p i}{r(r-p-q)_q}, \quad 
a' = \dfrac{(r-p)i'-q j'}{r(r-p-q)_p}, \quad   
b' = \dfrac{(r-q)j'-p i'}{r(r-p-q)_q}. 
\end{equation}
%%%%%%%%%%%%%%%%%%%%%%%%%%%%%%%%%%%%%%%%%%%%%%%%%%%%%%%%%%%%%%%%%%%%%%%%
Thus all the conditions in case 2 of Table \ref{tab:ab} have been 
obtained. \\[1mm]
%%%%%%%%%%%%%%%%%%%%%%%%%%%%% case 3 %%%%%%%%%%%%%%%%%%%%%%%%%%%%%%%%%%%
{\bf Case 3}. By Lemma \ref{lem:ab} there exists a quadruple of 
nonnegative integers $(i,j;i',j')$ such that 
%%%%%%%%%%%%%%%%%%%%%%%%%%%%% eqn:c3-1 %%%%%%%%%%%%%%%%%%%%%%%%%%%%%%%%%
\begin{equation} \label{eqn:c3-1}
a = \dfrac{p i}{r}, \qquad  
b = \dfrac{q (r j- p i)}{r(r-p)}, \qquad  
a'= \dfrac{p i'}{r}, \qquad  
b'= \dfrac{q(r j'- p i')}{r(r-p)}. 
\end{equation}
%%%%%%%%%%%%%%%%%%%%%%%%%%%%%%%%%%%%%%%%%%%%%%%%%%%%%%%%%%%%%%%%%%%%%%%
Substituting formula \eqref{eqn:c3-1} into equations \eqref{eqn:duality2} 
yields $i+i' = \frac{r}{p}-2 \in \Z$ and $r(j+j')-p(i+i') = (r-p) 
\left(\frac{r}{q}-2\right)$, the former of which is put into the 
latter to yields $j+j' = (r-p-q)/q \in \Z$. 
Thus we have $p|r$ and $q|(r-p-q)$, and hence  
%%%%%%%%%%%%%%%%%%%%%%%%%%%%% eqn:ls-c3 %%%%%%%%%%%%%%%%%%%%%%%%%%%%%%%
\begin{equation} \label{eqn:ls-c3}
i+i' = r_p -2, \qquad j+j' = (r-p-q)_q, \qquad 
i, j, i', j' \in \Z_{\ge 0}. 
\end{equation}
%%%%%%%%%%%%%%%%%%%%%%%%%%%%%%%%%%%%%%%%%%%%%%%%%%%%%%%%%%%%%%%%%%%%%%%
Formula \eqref{eqn:c3-1} then becomes 
%%%%%%%%%%%%%%%%%%%%%%%%%%%%% eqn:c3-2 %%%%%%%%%%%%%%%%%%%%%%%%%%%%%%%%
\begin{equation} \label{eqn:c3-2} 
a = \dfrac{i}{r_p}, \qquad  
b = \dfrac{r_p \, j - i}{r_p (r-p)_q}, \qquad  
a'= \dfrac{i'}{r_p}, \qquad  
b'= \dfrac{r_p \, j' - i'}{r_p (r-p)_q}.  
\end{equation}
%%%%%%%%%%%%%%%%%%%%%%%%%%%%%%%%%%%%%%%%%%%%%%%%%%%%%%%%%%%%%%%%%%%%%%
Thus all the conditions in case 3 of Table \ref{tab:ab} have been 
obtained. \\[1mm]
%%%%%%%%%%%%%%%%%%%%%%%%%%%%% case 4 %%%%%%%%%%%%%%%%%%%%%%%%%%%%%%%%%%% 
{\bf Case 4}. By Lemma \ref{lem:ab} there exists a quadruple of 
nonnegative integers $(i,j;i',j')$ such that 
%%%%%%%%%%%%%%%%%%%%%%%%%%%%% eqn:c4-1 %%%%%%%%%%%%%%%%%%%%%%%%%%%%%%%%%
\begin{equation} \label{eqn:c4-1}
a = \dfrac{p i}{r}, \qquad b = \dfrac{q j}{r}, \qquad 
a' = \dfrac{p i'}{r}, \qquad b' = \dfrac{q (r j'-p i')}{r(r-p)}. 
\end{equation}
%%%%%%%%%%%%%%%%%%%%%%%%%%%%%%%%%%%%%%%%%%%%%%%%%%%%%%%%%%%%%%%%%%%%%%%%
Substituting formula \eqref{eqn:c4-1} into equations \eqref{eqn:duality2} 
yields 
%%%%%
\[
i+i' = \textstyle \frac{r}{p} - 2 \in \Z, \qquad 
(r-p)j - p i' + r j' = (r-p)\left(\textstyle \frac{r}{q}-2\right).  
\]
%%%%% 
The former equation implies $p|r$ and $i+i' = r_p - 2$, while the 
division of the latter by $p$ makes $(r_p -1)j -i' + r_p \, j' = (r_p -1) 
\left(\frac{r}{q} - 2\right)$, to which $i+i'=r_p - 2$ is added to give 
$i + (r_p-1)j + r_p \, j' = r_p(r-p-q)/q \in \Z$. 
Summing up we have $q|r_p(r-p-q)$ and 
%%%%%%%%%%%%%%%%%%%%%%%%%%%%% eqn:ls-c4 %%%%%%%%%%%%%%%%%%%%%%%%%%%%%%%
\begin{equation} \label{eqn:ls-c4}
i+i' = r_p - 2, \qquad i + (r_p-1)j + r_p \, j' = (r_p(r-p-q))_q, 
\qquad i, j, i', j' \in \Z_{\ge 0}. 
\end{equation}
%%%%%%%%%%%%%%%%%%%%%%%%%%%%%%%%%%%%%%%%%%%%%%%%%%%%%%%%%%%%%%%%%%%%%%%
Formula \eqref{eqn:c4-1} then becomes 
%%%%%%%%%%%%%%%%%%%%%%%%%%%%% eqn:c4-2 %%%%%%%%%%%%%%%%%%%%%%%%%%%%%%%%
\begin{equation} \label{eqn:c4-2}  
a = \dfrac{i}{r_p}, \qquad b = \dfrac{q j}{r}, \qquad 
a' = \dfrac{i'}{r_p}, \qquad b' = 
\dfrac{r_p \, j'- i'}{(r_p(r -p))_q}. 
\end{equation}
%%%%%%%%%%%%%%%%%%%%%%%%%%%%%%%%%%%%%%%%%%%%%%%%%%%%%%%%%%%%%%%%%%%%%%%
Thus all the conditions in case 4 of Table \ref{tab:ab} have been 
obtained. \\[1mm]
%%%%%%%%%%%%%%%%%%%%%%%%%%%%% case 5 %%%%%%%%%%%%%%%%%%%%%%%%%%%%%%%%%%%
{\bf Case 5}. By Lemma \ref{lem:ab} there exists a quadruple of 
nonnegative integers $(i,j;i',j')$ such that 
%%%%%%%%%%%%%%%%%%%%%%%%%%%%% eqn:c5-1 %%%%%%%%%%%%%%%%%%%%%%%%%%%%%%%%%
\begin{equation} \label{eqn:c5-1}
a = \dfrac{p \{(r-p)i-q j\}}{r(r-p-q)}, \qquad  
b = \dfrac{q \{(r-q) j- p i\}}{r(r-p-q)}, \qquad  
a'= \dfrac{p(r i'-q j')}{r(r-q)}, \qquad  
b'= \dfrac{q j'}{r}. 
\end{equation}
%%%%%%%%%%%%%%%%%%%%%%%%%%%%%%%%%%%%%%%%%%%%%%%%%%%%%%%%%%%%%%%%%%%%%%%
Substituting formula \eqref{eqn:c5-1} into equations \eqref{eqn:duality2} 
yields 
%%%%%%%%%%%%%%%%%%%%%%%%%%% eqn:c5-2 %%%%%%%%%%%%%%%%%%%%%%%%%%%%%%%%%%
\begin{subequations} \label{eqn:c5-2}
\begin{align}
\resizebox{0.6\hsize}{!}{$\displaystyle (r-p)(r-q) i -q(r-q) j + r(r-p-q)i' - q(r-p-q) j'$} 
&= \resizebox{0.2\hsize}{!}{$(r-p-q)(r-q)\left(\frac{r}{p} - 2\right)$}, 
\label{eqn:c5-2a} \\ 
-p i + (r-q)j + (r-p-q) j' &= 
(r-p-q)\left(\textstyle \frac{r}{q}-2\right). \label{eqn:c5-2b}
\end{align} 
\end{subequations}
%%%%%%%%%%%%%%%%%%%%%%%%%%%%%%%%%%%%%%%%%%%%%%%%%%%%%%%%%%%%%%%%%%%%%%%  
Calculating $\eqref{eqn:c5-2a} + q \times \eqref{eqn:c5-2b}$ we have  
$i+i' = (r-p-q)/p \in \Z$, which implies $p|(r-p-q)$ and so  
$i+i' = (r-p-q)_p$.  
Division of \eqref{eqn:c5-2b} by $p$ makes $-i + (r-q)_p \, j + 
(r-p-q)_p \, j' = (r-p-q)_p \left(\frac{r}{q}-2\right)$, to which 
$i+i' = (r-p-q)_p$ is added to yield $i' + (r-q)_p \, j + 
(r-p-q)_p \, j' = (r-q)(r-p-q)_p/q \in \Z$. 
Summing up we have $q|r(r-p-q)_p$ and  
%%%%%%%%%%%%%%%%%%%%%%%%%%%%% eqn:ls-c5 %%%%%%%%%%%%%%%%%%%%%%%%%%%%%%%
\begin{subequations} \label{eqn:ls-c5}
\begin{align}
i+i' &= (r-p-q)_p,  & i, i' &\in \Z_{\ge 0}, \label{eqn:ls-c5a} \\
i' + (r-q)_p \, j + (r-p-q)_p \, j' &= ((r-q)(r-p-q)_p)_q, & 
i', j, j' &\in \Z_{\ge 0}. \label{eqn:ls-c5b}
\end{align} 
\end{subequations}
%%%%%%%%%%%%%%%%%%%%%%%%%%%%%%%%%%%%%%%%%%%%%%%%%%%%%%%%%%%%%%%%%%%%%%% 
Formula \eqref{eqn:c5-1} then becomes 
%%%%%%%%%%%%%%%%%%%%%%%%%%%%% eqn:c5-3 %%%%%%%%%%%%%%%%%%%%%%%%%%%%%%%%
\begin{equation} \label{eqn:c5-3} 
a = \dfrac{(r-p)i-q j}{r(r-p-q)_p}, \qquad  
b = \dfrac{(r-q)_p \, j- i}{(r(r-p-q)_p)_q}, \qquad  
a'= \dfrac{r i'-q j'}{r(r-q)_p}, \qquad  
b'= \dfrac{q j'}{r}. 
\end{equation}
%%%%%%%%%%%%%%%%%%%%%%%%%%%%%%%%%%%%%%%%%%%%%%%%%%%%%%%%%%%%%%%%%%%%%%%
Thus all the conditions in case 5 of Table \ref{tab:ab} have been 
obtained. \\[1mm]
%%%%%%%%%%%%%%%%%%%%%%%%%%%%% case 6 %%%%%%%%%%%%%%%%%%%%%%%%%%%%%%%%%%%
{\bf Case 6}. By Lemma \ref{lem:ab} there exists a quadruple of 
nonnegative integers $(i,j;i',j')$ such that 
%%%%%%%%%%%%%%%%%%%%%%%%%%%%% eqn:c6-1 %%%%%%%%%%%%%%%%%%%%%%%%%%%%%%%%%
\begin{equation} \label{eqn:c6-1}
a = \dfrac{p i}{r}, \qquad  
b = \dfrac{q (r j- p i)}{r(r-p)}, \qquad 
a'= \dfrac{p(r i'- q j')}{r(r-q)}, \qquad  
b'= \dfrac{q j'}{r}. 
\end{equation}
%%%%%%%%%%%%%%%%%%%%%%%%%%%%%%%%%%%%%%%%%%%%%%%%%%%%%%%%%%%%%%%%%%%%%%%
Substituting formula \eqref{eqn:c6-1} into equations \eqref{eqn:duality2} 
yields
%%%%%%%%%%%%%%
\[ 
(r-q)i - q j' + r i' = (r-q)\left(\textstyle \frac{r}{p}-2\right), \qquad 
-p i + (r-p)j'+ r j  = (r-p)\left(\textstyle \frac{r}{q}-2\right).   
\]
%%%%%%%%%%%%%%
They are recast to $(r-p-q)i +(r-p) i' + q j = (r-p)(r-p-q)/p \in \Z$ 
and $(r-p-q)j'+ p i' + (r-q)j = (r-q)(r-p-q)/q \in \Z$. 
Thus $p|r(r-p-q)$ and $q|r(r-p-q)$, so that 
%%%%%%%%%%%%%%%%%%%%%%%%%%%%% eqn:ls-c6 %%%%%%%%%%%%%%%%%%%%%%%%%%%%%%%%%
\begin{subequations} \label{eqn:ls-c6}
\begin{align}
(r-p-q)i +(r-p) i' + q j &= ((r-p)(r-p-q))_p, & i, i', j &\in \Z_{\ge0},  
\label{eqn:ls-c6a} \\ 
(r-p-q)j'+ p i' + (r-q)j &= ((r-q)(r-p-q))_q, & j', i', j &\in \Z_{\ge0}.  
\label{eqn:ls-c6b}  
\end{align}   
\end{subequations}
%%%%%%%%%%%%%%%%%%%%%%%%%%%%%%%%%%%%%%%%%%%%%%%%%%%%%%%%%%%%%%%%%%%%%%%%
Formula \eqref{eqn:c6-1} then becomes 
%%%%%%%%%%%%%%%%%%%%%%%%%%%%%%%% eqn:c6-2 %%%%%%%%%%%%%%%%%%%%%%%%%%%%%
\begin{equation} \label{eqn:c6-2}
a = \dfrac{p i}{r}, \qquad  
b = \dfrac{r j- p i}{(r(r-p))_q}, \qquad 
a'= \dfrac{r i'- q j'}{(r(r-q))_p}, \qquad  
b'= \dfrac{q j'}{r}. 
\end{equation}
%%%%%%%%%%%%%%%%%%%%%%%%%%%%%%%%%%%%%%%%%%%%%%%%%%%%%%%%%%%%%%%%%%%%%%
Thus all the conditions in case 6 of Table \ref{tab:ab} have been 
obtained. \hfill $\Box$
%%%%%%%%%%%%%%%%%%%%%%%%%%%%% end proof %%%%%%%%%%%%%%%%%%%%%%%%%%%%%%%% 
%%%%%%%%%%%%%%%%%%%%%%%%%%%%% lem:case7 %%%%%%%%%%%%%%%%%%%%%%%%%%%%%%%%
\begin{lemma} \label{lem:case7} 
Case $7$ of Table $\ref{tab:type}$ cannot occur. 
\end{lemma}
%%%%%%%%%%%%%%%%%%%%%%%%%%%%%%%%%%%%%%%%%%%%%%%%%%%%%%%%%%%%%%%%%%%%%%%
%%%%%%%%%%%%%%%%%%%%%%%%%%%%% begin proof %%%%%%%%%%%%%%%%%%%%%%%%%%%%%
{\it Proof}. 
By Lemma \ref{lem:ab} there exists a quadruple of nonnegative 
integers $(i,j;i',j')$ such that 
%%%%%
\[
a = \dfrac{p i}{r}, \qquad 
b = \dfrac{q j}{r}, \qquad 
a' = \dfrac{p \{(r-p) i'-q j'\}}{r(r-p-q)}, \qquad 
b' =  \dfrac{q \{(r-q) j'-p i'\}}{r(r-p-q)}. 
\]
%%%%%
Substituting these into relations \eqref{eqn:duality2} yields 
%%%%%%%%%%%%%%%
\begin{align*}
(r-p-q) i + (r-p) i' -q j' &= (r-p-q)\left(\textstyle \frac{r}{p} -2\right), \\
(r-p-q) j + (r-q) j' -p j' &= (r-p-q)\left(\textstyle \frac{r}{q} -2\right), 
\end{align*}
%%%%%%%%%%%%%%
which can readily be converted into
%%%%%%%%%%%%%%%%%%%%%%%%%%%%% eqn:c7-1 %%%%%%%%%%%%%%%%%%%%%%%%%%%%%%%
\begin{subequations} \label{eqn:c7-1}
\begin{align}
(r-q) i + q j + r i' &= \textstyle \frac{r}{p} -2 \in \Z, \label{eqn:c7-1a} \\
p i + (r-p) j + r j' &= \textstyle \frac{r}{q} -2 \in \Z, \label{eqn:c7-1b},  
\end{align}
\end{subequations}
%%%%%%%%%%%%%%%%%%%%%%%%%%%%%%%%%%%%%%%%%%%%%%%%%%%%%%%%%%%%%%%%%%%%%%% 
which imply $p|r$ and $q|r$. 
Now $r$ must be even, for otherwise $r$ is odd and hence so are 
$p$ and $q$ by $p|r$ and $q|r$, but then $r-p-q \equiv 1-1-1 \equiv 
1 \mod 2$, that is, $r-p-q$ is odd, which is absurd since it must be 
even by \cite[assertion (2) of Theorem 2.2]{Iwasaki}.     
If $i' \ge 1$ then \eqref{eqn:c7-1a} gives an absurd estimate  
$r = (r-q) \cdot 0 + q \cdot 0 + r \cdot 1 \le (r-q) i + q j + r i' = 
r_p -2 \le r-2$. 
Thus we must have $i' = 0$. 
In a similar manner equation \eqref{eqn:c7-1b} forces $j' = 0$. 
Therefore \eqref{eqn:c7-1} reduces to 
%%%%%%%%%%%%%%%%%%%%%%%%%%%%% eqn:c7-2 %%%%%%%%%%%%%%%%%%%%%%%%%%%%%%%
\begin{subequations} \label{eqn:c7-2}
\begin{align}
(r-q) i + q j &= r_p -2, \label{eqn:c7-2a} \\
p i + (r-p) j &= r_q -2, \label{eqn:c7-2b},  
\end{align}
\end{subequations}
%%%%%%%%%%%%%%%%%%%%%%%%%%%%%%%%%%%%%%%%%%%%%%%%%%%%%%%%%%%%%%%%%%%%%%% 
If $i \ge 1$ and $j \ge 1$ then \eqref{eqn:c7-2a} gives an absurd 
estimate $r = (r-q) \cdot 1 + q \cdot 1 \le (r-q) i + q j = r_p -2 
\le r-2$. 
Thus we must have either $i = 0$ or $j = 0$. 
Here we may and shall assume $j = 0$ due to the symmetry of conditions 
\eqref{eqn:c7-2a} and \eqref{eqn:c7-2b}. 
Then \eqref{eqn:c7-2a} and \eqref{eqn:c7-2b} yield 
%%%%%%%%%%%%%%%%%%%%%%%%%%%%% eqn:c7-3 %%%%%%%%%%%%%%%%%%%%%%%%%%%%%%%%
\begin{equation} \label{eqn:c7-3}
i = \frac{r_p-2}{r-q} = \frac{r_q-2}{p} \in \Z,  
\end{equation}
%%%%%%%%%%%%%%%%%%%%%%%%%%%%%%%%%%%%%%%%%%%%%%%%%%%%%%%%%%%%%%%%%%%%%%%
which in particular implies $p|(r_q-2)$. 
Now $r_q$ must be even, for otherwise $r_q$ is odd and hence so is 
$p$ by $p|(r_q-2)$, while $q$ is even since $r = q r_q$ is even with  
$r_q$ odd, but then $r-p-q \equiv 0-1-0 \equiv 1 \mod 2$, that is, 
$r-p-q$ is odd, which is again absurd. 
By the second equality in \eqref{eqn:c7-3} we have 
$2 q = 4 r - 2 p - r \cdot r_q = 4 r - 2 p - 2 r \cdot r_{2q}$, that is, 
$q/p = 2 r_p - 1 - r_p \cdot r_{2q} \in \Z$, where $r_{2 q}$ makes 
sense as $r_q$ is even. 
Thus $p|q$ and so we have integer equations 
$1 = 2 r_p - q_p - r_p \cdot r_{2 q} = 4 r_{2 q} 
\cdot q_p - q_p - 2 (r_{2 q})^2 q_p = 
q_p \cdot \{4 r_{2 q} - 2 (r_{2 q})^2 - 1\}$, where   
$r_p = 2 r_{2 q} \cdot q_p$ is used in the second equality. 
Thus $q_p = 1$ and $4 r_{2 q}- 2 (r_{2 q})^2 -1 = 1$, 
the former of which means $p = q$ while the latter yields  
$r_{2 q} (2-r_{2 q}) = 1$, that is, $r_{2 q} = 1$ and so 
$r = 2 q$. 
Then $r-p-q = r-2 q = 0$, which is absurd since we have 
$r-p-q > 0$ in $\cD^-$. 
This last contradiction shows that the occurrence of case $7$ in 
Table \ref{tab:type} is impossible. \hfill $\Box$
%%%%%%%%%%%%%%%%%%%%%%%%%%%%% end proof %%%%%%%%%%%%%%%%%%%%%%%%%%%%%%%
%%%%%%%%%%%%%%%%%%%%%%%%%%%%% lem:ab-bound %%%%%%%%%%%%%%%%%%%%%%%%%%%%
\begin{lemma} \label{lem:ab-bound} 
We must have $-1 < a, \, b, \, a', \, b' < 1$. 
\end{lemma} 
%%%%%%%%%%%%%%%%%%%%%%%%%%%%%%%%%%%%%%%%%%%%%%%%%%%%%%%%%%%%%%%%%%%%%%%
%%%%%%%%%%%%%%%%%%%%%%%%%%%%% begin proof %%%%%%%%%%%%%%%%%%%%%%%%%%%%%
{\it Proof}. The lemma is proved by a case-by-case check. \\[1mm]
%%%%%%%%%%%%%%%%%%%%%%%%%%%%% case 1 %%%%%%%%%%%%%%%%%%%%%%%%%%%%%%%%%%% 
{\bf Case 1}. From condition \eqref{eqn:ls-c1} we have    
$0 \le i, \, i' \le r_p-2$ and $0 \le j, \, j' \le r_q-2$, 
so that estimate $0 \le a, \, a', \, b, \, b' < 1$ follows from 
representation \eqref{eqn:c1-2}. \\[1mm]  
%%%%%%%%%%%%%%%%%%%%%%%%%%%%% case 2 %%%%%%%%%%%%%%%%%%%%%%%%%%%%%%%%%%% 
{\bf Case 2}. From condition \eqref{eqn:ls-c2} we have $i \ge 0$ 
and $j \le (r-p-q)_q$ as well as $i \le (r-p-q)_p$ and $j \ge 0$.  
Thus it follows from representation \eqref{eqn:c2-2} that      
%%%%%%%%%%%%%%
\begin{align*}
1+ a &= \dfrac{r(r-p-q)_p+(r-p)i-q j}{r(r-p-q)_p} \ge 
\dfrac{r(r-p-q)_p - q(r-p-q)_q}{r(r-p-q)_p} \\ 
&= \dfrac{r(r-p-q)_p - p (r-p-q)_p}{r(r-p-q)_p} = \dfrac{r-p}{r} > 0, \\       
1- a &= \dfrac{r(r-p-q)_p-(r-p)i+q j}{r(r-p-q)_p} \ge 
\dfrac{r(r-p-q)_p - (r-p)(r-p-q)_p}{r(r-p-q)_p} 
= \dfrac{p}{r} > 0,  
\end{align*}
%%%%%%%%%%%%% 
which shows $-1 < a < 1$. 
In similar manners we have $-1 < b, \, a', \, b' < 1$. \\[1mm]
%%%%%%%%%%%%%%%%%%%%%%%%%%%%% case 3 %%%%%%%%%%%%%%%%%%%%%%%%%%%%%%%%%%% 
{\bf Case 3}. From condition \eqref{eqn:ls-c3} we have  
$0 \le i, \, i' \le r_p-2$, so that estimate $0 \le a, \, a' < 1$ 
follows from representation \eqref{eqn:c3-2}. 
Similarly, from condition \eqref{eqn:ls-c3} we have $j \ge 0$ and 
$i \le r_p-2$ as well as $j \le (r-p-q)_q$ and $i \ge 0$.  
Thus it follows from representation \eqref{eqn:c3-2} that 
%%%%%%%%%%%%%%%
\begin{align*}
1+b &= \dfrac{r_p(r-p)_q + r_p j-i}{r_p(r-p)_q} 
\ge \dfrac{r_p(r-p)_q -(r_p -2)}{r_p(r-p)_q} = 
\dfrac{r_p(r-p-q)_q + 2}{r_p(r-p)_q} > 0, \\  
1-b &= \dfrac{r_p(r-p)_q - r_p j+i}{r_p(r-p)_q} 
\ge \dfrac{r_p(r-p)_q -r_p (r-p-q)_q}{r_p(r-p)_q} = 
\dfrac{1}{(r-p)_q} > 0,   
\end{align*}
%%%%%%%%%%%%%%%
which shows $-1 < b < 1$. 
In a similar manner we have $-1 < b' < 1$.  
\\[1mm]
%%%%%%%%%%%%%%%%%%%%%%%%%%%%% case 4 %%%%%%%%%%%%%%%%%%%%%%%%%%%%%%%%%%% 
{\bf Case 4}. From condition \eqref{eqn:ls-c4} we have 
$0 \le i, \, i' \le r_p-2$, so that estimate $0 \le a, \, a' < 1$ 
follows from representation \eqref{eqn:c4-2}. 
Similarly, from condition \eqref{eqn:ls-c4} we have 
$0 \le (r_p -1) j \le (r_p(r-p-q))_q$, that is, 
$0 \le q j \le r(r-p-q)/(r-p)$. 
Thus representation \eqref{eqn:c4-2} yields   
%%%%
\[
0 \le b = \dfrac{q j}{r} \le \dfrac{r-p-q}{r-p} < 1. 
\]
%%%% 
From condition \eqref{eqn:ls-c4} we have $j' \ge 0$ and $i' \le r_p-2$ 
as well as $r_p j' \le (r_p(r-p-q))_q$ and $i' \ge 0$.  
Thus it follows from representation \eqref{eqn:c4-2} that  
%%%%%%%%%%%%%%%%
\begin{align*}
1+b' &= \dfrac{(r_p(r-p))_q + r_p j' - i'}{(r_p(r-p))_q} 
\ge \dfrac{(r_p(r-p))_q -(r_p -2)}{(r_p(r-p))_q} = 
 \dfrac{(r_p(r-p-q))_q + 2}{(r_p(r-p))_q} > 0, \\
1-b' &= \dfrac{(r_p(r-p))_q - r_p j' + i'}{(r_p(r-p))_q} 
\ge \dfrac{(r_p(r-p))_q -(r_p(r-p-q))_q}{(r_p(r-p))_q} = 
 \dfrac{q}{r-p} > 0.  
\end{align*}
%%%%%%%%%%%%%%%   
Therefore we have $0 \le a, \, a', \, b < 1$ and $-1 < b' < 1$. \\[1mm]
%%%%%%%%%%%%%%%%%%%%%%%%%%%%% case 5 %%%%%%%%%%%%%%%%%%%%%%%%%%%%%%%%%%% 
{\bf Case 5}. From condition \eqref{eqn:ls-c5b} we have  
$(r-q)_p j \le ((r-q)(r-p-q)_p)_q$, that is, $q j \le r-p-q = 
p(r-p-q)_p$. 
Thus it follows from representation \eqref{eqn:c5-3} and $i \ge 0$ that 
%%%%%
\[
1+a = \dfrac{r(r-p-q)_p + (r-p)i-q j}{r(r-p-q)_p} 
\ge \dfrac{r(r-p-q)_p - p(r-p-q)_p}{r(r-p-q)_p} = \frac{r-p}{r} > 0.  
\]
%%%%%
From condition \eqref{eqn:ls-c5a} we have $i \le (r-p-q)_p$.  
It then follows from \eqref{eqn:c5-3} and $j \ge 0$ that 
%%%%%%%%%%%%%%%
\begin{align*}
1-a &= \dfrac{r(r-p-q)_p - (r-p)i+ q j}{r(r-p-q)_p} 
\ge \dfrac{r(r-p-q)_p - (r-p)(r-p-q)_p}{r(r-p-q)_p} = \frac{p}{r} > 0, \\  
1+b &= \dfrac{(r(r-p-q)_p)_q + (r-q)_p j -i}{(r(r-p-q)_p)_q} 
\ge \dfrac{(r(r-p-q)_p)_q -(r-p-q)_p}{(r(r-p-q)_p)_q} = 
\dfrac{r-q}{r} > 0.  
\end{align*}
%%%%%%%%%%%%%%%
Since $(r-q)_p j \le ((r-q)(r-p-q)_p)_q$ and $i \ge 0$, 
representation \eqref{eqn:c5-3} yields   
%%%%
\[
1-b = \dfrac{(r(r-p-q)_p)_q -(r-q)_p j + i}{(r(r-p-q)_p)_q -(r-q)_p j + i} 
\ge \dfrac{(r(r-p-q)_p)_q - ((r-q)(r-p-q)_p)_q}{(r(r-p-q)_p)_q -(r-q)_p j + i} 
= \dfrac{q}{r} > 0. 
\]
%%%%  
From condition \eqref{eqn:ls-c5b} we have $(r-p-q)_p j' \le ((r-q)(r-p-q)_p)_q$, 
that is, $q j' \le r-q = p(r-q)_p$. 
Thus it follows from formula \eqref{eqn:c5-3} and $i' \ge 0$ that  
%%%%
\[
1+a' = \dfrac{r(r-q)_p + r i'-q j'}{r(r-q)_p} 
\ge  \dfrac{r(r-q)_p - p(r-q)_p}{r(r-q)_p} = \dfrac{r-p}{r} > 0.  
\]
%%%%
From condition \eqref{eqn:ls-c5a} we have $i' \le (r-p-q)_p$.  
Then it follows from \eqref{eqn:c5-3} and $j' \ge 0$ that  
%%%%
\[
1-a' = \dfrac{r(r-q)_p - r i'+ q j'}{r(r-q)_p} 
\ge  \dfrac{r(r-q)_p - r(r-p-q)_p}{r(r-q)_p} = \dfrac{1}{(r-q)_p} > 0.  
\]
%%%%
Finally, since $0 \le q j' \le r-q$, 
we have $0 \le b' = q j'/r \le (r-q)/r < 1$.  
\\[1mm]
%%%%%%%%%%%%%%%%%%%%%%%%%%%%% case 6 %%%%%%%%%%%%%%%%%%%%%%%%%%%%%%%%%%% 
{\bf Case 6}. From condition \eqref{eqn:ls-c6a} we have $(r-p-q) i \le 
((r-p)(r-p-q))_p$, that is, $p i \le r-p$. 
Thus it follows from representation \eqref{eqn:c6-2} and $i \ge 0$ that 
$0 \le a = p i/r \le (r-p)/r < 1$. 
Similarly, we have $0 \le b' < 1$. 
From condition \eqref{eqn:ls-c6b} we have $(r-p-q)j' \le ((r-q)(r-p-q))_q$, 
that is, $q j' \le r-q$. 
Thus it follows from representation \eqref{eqn:c6-2} and $i' \ge 0$ that
%%%%%%
\[
1+ a' = \dfrac{(r(r-q))_p + r i'-q j'}{(r(r-q))_p} 
\ge \dfrac{(r(r-q))_p - (r-q)}{(r(r-q))_p} = \dfrac{r-p}{r} > 0. 
\]
%%%%%% 
On the other hand, from condition \eqref{eqn:ls-c6a} we have 
$(r-p)i' \le ((r-p)(r-p-q))_p$, that is, $p i' \le r-p-q$.  
Thus it follows from representation \eqref{eqn:c6-2} and $j' \ge 0$ that
%%%%%%
\[
1-a' = \dfrac{(r(r-q))_p - r i' + q j'}{(r(r-q))_p} 
\ge \dfrac{(r(r-q))_p - r i'}{(r(r-q))_p}
= \dfrac{(r-q) - p i'}{r-q} \ge \dfrac{p}{r-q} > 0. 
\]
%%%%%% 
In a similar manner we have $-1 < b < 1$.  \hfill $\Box$
%%%%%%%%%%%%%%%%%%%%%%%%%%%%% end proof %%%%%%%%%%%%%%%%%%%%%%%%%%%%%%%%
%%%%%%%%%%%%%%%%%%%%%%%%%%%%% lem:incl %%%%%%%%%%%%%%%%%%%%%%%%%%%%%%%%%
\begin{lemma} \label{lem:incl} 
We have $0 \le a, \, b < 1$ and there exist inclusions of multi-sets:   
%%%%%%%%%%%%%%%%%%%%%%%% eqn:incl %%%%%%%%%%%%%%%%%%%%%%%%%%%%%%%%%%%%%%
\begin{subequations} \label{eqn:incl}
\begin{align}
\left\{ \frac{i+a}{p}\right\}_{i\in \bar{I}_p} \bigcup \,\,  
\left\{\frac{i+b}{q}\right\}_{i \in \bar{I}_q}
& \subset \, 
\left\{\frac{j-a}{r-p}\right\}_{j \in J_p} \bigcup \,\, 
\left\{\frac{j-b}{r-q}\right\}_{j \in J_q}, 
\label{eqn:incl-a} \\
\left\{ \frac{i+a}{p}\right\}_{i\in I_p} \bigcup \,\,  
\left\{\frac{i+b}{q}\right\}_{i \in I_q}
& \subset \, 
\left\{\frac{j-a}{r-p}\right\}_{j \in \bar{J}_p} \bigcup \,\, 
\left\{\frac{j-b}{r-q}\right\}_{j \in \bar{J}_q},   
\label{eqn:incl-b}
\end{align}    
\end{subequations}
%%%%%%%%%%%%%%%%%%%%%%%%%%%%%%%%%%%%%%%%%%%%%%%%%%%%%%%%%%%%%%%%%%%%%%%%
where $\bar{I}_p := [p] \setminus I_p$, $\bar{I}_q := [q] \setminus I_q$, 
$\bar{J}_p := [r-p] \setminus J_p$ and $\bar{J}_q := [r-q] \setminus J_q$.  
\end{lemma}
%%%%%%%%%%%%%%%%%%%%%%%%%%%%%%%%%%%%%%%%%%%%%%%%%%%%%%%%%%%%%%%%%%%%%%%%
%%%%%%%%%%%%%%%%%%%%%%%%%%%% begin proof %%%%%%%%%%%%%%%%%%%%%%%%%%%%%%%
{\it Proof}. Recall that $J_p \subset [r-p] = \{0, \dots, r-p-1\}$.  
By Lemma \ref{lem:ab-bound} we have $-a < 1$ so that 
$(j-a)/(r-p) < (r-p-1+1)/(r-p) = 1$ for any $j \in J_p$.  
If $a \le 0$ then obviously we have $(j-a)/(r-p) \ge 0$ for any 
$j \in J_p$.  
If $a > 0$ then we have $J_p \subset \{1, \dots, r-p-1\}$ by 
Lemma \ref{lem:ab-J} and $-a > -1$ by Lemma \ref{lem:ab-bound}, so that  
$(j-a)/(r-p) > (1-1)/(r-p) = 0$ for any $j \in J_p$. 
In either case we have $0 \le (j-a)/(r-p) < 1$ for any $j \in J_p$. 
In a similar manner we have $0 \le (j-b)/(r-q) < 1$ for any $j \in J_q$. 
So the multi-set on the right-hand side of \eqref{eqn:succ1} lies  
in the interval $[0, \, 1)$. 
Since $0 \in \bar{I}_p$ and $0 \in \bar{I}_q$, the numbers 
$a/p$ and $b/q$ belong to the multi-set on the left so that they 
must be nonnegative by the binary relation \eqref{eqn:succ1}. 
Hence we have $a, \, b \ge 0$, which together with 
Lemma \ref{lem:ab-bound} yields the estimate $0 \le a, \, b < 1$. 
%%%%%
\par
%%%%%
Since $0 \le a < 1$ and $\bar{I}_p \subset \{0, \dots, p-1\}$ we 
have $0 = (0+0)/p \le (i+a)/p < (p-1+1)/p = 1$ for any $i \in \bar{I}_p$. 
In a similar manner we have $0 \le (i+b)/q < 1$ for any $i \in \bar{I}_q$. 
Thus the multi-set on the left-hand side of \eqref{eqn:succ1} also lies 
in the interval $[0, \, 1)$, as does the multi-set on the right. 
Therefore binary relation \eqref{eqn:succ1} must be the inclusion 
\eqref{eqn:incl-a}.  
%%%%%
\par
%%%%%
To prove inclusion \eqref{eqn:incl-b}, we use the dual version of 
inclusion \eqref{eqn:incl-a}: 
%%%%%%%%%%%%%%%%%%%%%%%%%%% eqn:incl-a-p %%%%%%%%%%%%%%%%%%%%%%%%%%%%%
\begin{equation} \label{eqn:incl-a-p}
\left\{ \frac{i'+a'}{p}\right\}_{i' \in \bar{I}'_p} \bigcup \,\,  
\left\{\frac{i'+b'}{q}\right\}_{i' \in \bar{I}'_q}
\subset \, 
\left\{\frac{j'-a'}{r-p}\right\}_{j' \in J'_p} \bigcup \,\, 
\left\{\frac{j'-b'}{r-q}\right\}_{j' \in J'_q},  
\end{equation}
%%%%%%%%%%%%%%%%%%%%%%%%%%%%%%%%%%%%%%%%%%%%%%%%%%%%%%%%%%%%%%%%%%%%%% 
where $I'_p$, $I'_q$, $J'_p$, $J'_q$ are the dual counterparts 
of $I_p$, $I_q$, $J_p$, $J_q$, with $\bar{I}'_p := [p] \setminus 
I'_p$ and $\bar{I}'_q := [r-q] \setminus I'_q$. 
Dividing equation \eqref{eqn:vi-s} by \eqref{eqn:IJ}, we have     
%%%%%%
\[
\prod_{i=1}^r (w+v_i^*) = 
\prod_{i \in \bar{I}_p} \left(w+{\ts\frac{i+a}{p}}\right) 
\prod_{i \in \bar{I}_q} \left(w+{\ts\frac{i+b}{q}}\right) 
\prod_{j \in \bar{J}_p} \left(w+{\ts\frac{j-a}{r-p}}\right) 
\prod_{j \in \bar{J}_q} \left(w+{\ts\frac{j-b}{r-q}}\right).     
\]
%%%%% 
Taking the reflection $w \mapsto w'$ as in \eqref{eqn:w-p} and 
using definitions \eqref{eqn:duality} and \eqref{eqn:vi-p} there 
yield 
%%%%%%%%%%%%%%%
\begin{align*}
\prod_{i=1}^r (w+v_i') &= 
\prod_{i \in \bar{I}_p} \left(w+{\ts\frac{(p-1-i)+a'}{p}}\right) 
\prod_{i \in \bar{I}_q} \left(w+{\ts\frac{(q-1-i)+b'}{q}}\right) \\
&\phantom{=} \times 
\prod_{j \in \bar{J}_p} \left(w+{\ts\frac{(r-p-1-j)-a'}{r-p}}\right) 
\prod_{j \in \bar{J}_q} \left(w+{\ts\frac{(r-q-1-j)-b'}{r-q}}\right),    
\end{align*}
%%%%%%%%%%%%
which together with the definitions of $I'_s$, $\bar{I}'_s$ and 
$J'_s$, $s \in \{p, \, q\}$, implies 
%%%%%%%%%%%%%%%%%%%%%%%%% eqn:IJ-p %%%%%%%%%%%%%%%%%%%%%%%%%%%%%%% 
\begin{equation} \label{eqn:IJ-p}
\begin{split}
I'_s &= \{ i' = s-1-i \,:\, i \in \bar{I}_s \}, \qquad 
\bar{I}'_s = \{ i' = s -1-i \,:\, i \in I_s \}, \\   
J'_s &= \{ j' = r-s-1-j \,:\, j \in \bar{J}_s \}  
\qquad \qquad (s \in \{p, \, q\}).   
\end{split}
\end{equation}
%%%%%%%%%%%%%%%%%%%%%%%%%%%%%%%%%%%%%%%%%%%%%%%%%%%%%%%%%%%%%%%%%%% 
Under the correspondences $i' \leftrightarrow i$ and 
$j' \leftrightarrow j$ in \eqref{eqn:IJ-p}, it follows from 
definition \eqref{eqn:duality} that  
%%%%%% 
\[
\dfrac{i'+c'}{s} = 1 - \frac{2}{r} - \dfrac{i+c}{s}, \qquad 
\dfrac{j'-c'}{r-s} = 1 - \frac{2}{r} - \dfrac{j-c}{r-s}, 
\qquad ((s, c) = (p, a), \, (q, b)).
\]
%%%%%% 
Thus inclusion \eqref{eqn:incl-a-p} and relation \eqref{eqn:IJ-p} 
lead to inclusion \eqref{eqn:incl-b}. \hfill $\Box$   
%%%%%%%%%%%%%%%%%%%%%%%%%%%% end proof %%%%%%%%%%%%%%%%%%%%%%%%%%%%%%%%%
%%%%%%%%%%%%%%%%%%%%%%%%%%%%% thm:ab %%%%%%%%%%%%%%%%%%%%%%%%%%%%%%%%%%%%
\begin{theorem} \label{thm:ab} 
Any $(\rA)$-solution $\lambda = (p,q,r;a,b;x) \in \cD^-$ together with its dual 
solution $\lambda' = (p,q,r;a',b';x) \in \cD^-$ must be subject to the four 
conditions $(1)$--$(4)$ in Table $\ref{tab:ab}$ as well as to  
%%%%%%%%%%%%%%%%%%%%%%%%%%%% eqn:estimate2 %%%%%%%%%%%%%%%%%%%%%%%%%%%%%
\begin{equation} \label{eqn:estimate2} 
a, b, \, v_1, \dots, v_r \in \Q; \qquad 
0 \le a, \, b < 1, \qquad 0 \le v_1, \dots, v_r < 1, 
\end{equation} 
%%%%%%%%%%%%%%%%%%%%%%%%%%%%%%%%%%%%%%%%%%%%%%%%%%%%%%%%%%%%%%%%%%%%%%%%
along with the same condition for $a', b'$ and $v_1', \dots, v_r'$.         
\end{theorem}
%%%%%%%%%%%%%%%%%%%%%%%%%%%%%%%%%%%%%%%%%%%%%%%%%%%%%%%%%%%%%%%%%%%%%%%%%%
%%%%%%%%%%%%%%%%%%%%%%%%%%%% begin proof %%%%%%%%%%%%%%%%%%%%%%%%%%%%%%%
{\it Proof}. 
The first assertion concerning Table \ref{tab:ab} is an immediate 
consequence of Lemmas \ref{lem:case1-6} and \ref{lem:case7}. 
To prove the second assertion \eqref{eqn:estimate2}, notice that the 
rationality of $a$ and $b$ is already proved in Lemma \ref{lem:ab}, 
while $v_1, \dots, v_r \in \Q$ follows from $a, b \in \Q$ and 
formula \eqref{eqn:IJ}.  
The estimate $0 \le a, \, b < 1$ is already proved in 
Lemma \ref{lem:incl}, while $0 \le v_1, \dots, v_r < 1$ follows from 
a combination of $0 \le a, \, b < 1$, formula \eqref{eqn:IJ} and 
Lemma \ref{lem:ab-J}.  
Thus we have condition \eqref{eqn:estimate2}. \hfill $\Box$ 
%%%%%%%%%%%%%%%%%%%%%%%%%%%% end proof %%%%%%%%%%%%%%%%%%%%%%%%%%%%%%%%%%%  
%%%%%%%%%%%%%%%%%%%%%%%%%% rem:ab %%%%%%%%%%%%%%%%%%%%%%%%%%%%%%%%%%%%%%%%
\begin{remark} \label{rem:ab} 
Two remarks about Theorem \ref{thm:ab} should be in order at this stage. 
\begin{enumerate}
\item In each case of Table \ref{tab:ab} the $\Z$-linear 
equations for $(i,j;i',j') \in \Z_{\ge0}^4$ are of the form 
%%%%%%%%%%%%%%%%%%%%%%%%%%%%%% eqn:z-linear %%%%%%%%%%%%%%%%%%%%%%%%%%%%%%%
\begin{equation} \label{eqn:z-linear}
\mu_1 \, i + \mu_2 \, j + \mu_3 \, i' + \mu_4 \, j' = \mu_5, \qquad 
\nu_1 \, i + \nu_2 \, j + \nu_3 \, i' + \nu_4 \, j' = \nu_5, 
\end{equation}
%%%%%%%%%%%%%%%%%%%%%%%%%%%%%%%%%%%%%%%%%%%%%%%%%%%%%%%%%%%%%%%%%%%%%%%%%%%
where $\mu_k$, $\nu_k \in \Z_{\ge0}$, $\mu_k + \nu_k \ge 1$ for 
$k =1,2,3,4$, so that system \eqref{eqn:z-linear} cannot have 
infinitely many solutions. 
Thus the finite cardinality of $(\rA)$-solutions with a 
given principal part is an immediate corollary to Theorem \ref{thm:ab} 
and assertion (1) of \cite[Theorem 2.3]{Iwasaki}.  
\item A data $\lambda = (p,q,r;a,b;x) \in \cD^-$ with $\bp = (p,q;r) 
\in D_{\rA}^-$ is said to be a {\sl candidate} for an $(\rA)$-solution, 
if  $\lambda$ is subject to one of the six patterns in Table \ref{tab:ab} 
and if $\bp$ satisfies condition \eqref{eqn:dr} to be established in 
Proposition \ref{prop:dr}.  
Note that any data cannot be an $(\rA)$-solution unless it is a candidate, 
but it may (perhaps quite often) happen that a candidate is 
not an actual $(\rA)$-solution.   
\end{enumerate}
\end{remark}
%%%%%%%%%%%%%%%%%%%%%%%%%%%%%%%%%%%%%%%%%%%%%%%%%%%%%%%%%%%%%%%%%%%%%%%%%%%
\par
%%%%% 
Theorem \ref{thm:D-r} is contained in the following.   
%%%%%%%%%%%%%%%%%%%%%%%%%%%% prop:division %%%%%%%%%%%%%%%%%%%%%%%%%%%%% 
\begin{proposition} \label{prop:division} 
For any $(\rA)$-solution $\lambda = (p,q,r;a,b;x) \in \cD^-$ the 
division relation \eqref{eqn:division2} holds true along with 
yet another division relation            
%%%%%%%%%%%%%%%%%%%%%%%%%%%% eqn:division3 %%%%%%%%%%%%%%%%%%%%%%%%%%%%%
\begin{equation} \label{eqn:division3} 
\prod_{i=0}^{p-1} \left(w+{\ts\frac{i+a}{p}}\right) 
\prod_{i=0}^{q-1} \left(w+{\ts\frac{i+b}{q}}\right) 
\,\, \Big| \,\, 
\prod_{j=0}^{r-p-1} \left(w+{\ts\frac{j-a}{r-p}}\right)
\prod_{j=0}^{r-q-1} \left(w+{\ts\frac{j-b}{r-q}}\right) 
\qquad \mbox{in \,\, $\Q[w]$},  
\end{equation} 
%%%%%%%%%%%%%%%%%%%%%%%%%%%%%%%%%%%%%%%%%%%%%%%%%%%%%%%%%%%%%%%%%%%%%%%
and hence one can arrange $v_1, \dots, v_r$ so that equation 
\eqref{eqn:vi-pq} holds. 
With this convention the reciprocal $\check{\lambda}$ of $\lambda$ 
becomes a solution to Problem $\rI$ in $\cF^-$ with gamma product 
formula    
%%%%%%%%%%%%%%%%%%%%%%%%%%%%% eqn:gpf-c2 %%%%%%%%%%%%%%%%%%%%%%%%%%%%%% 
\begin{equation} \label{eqn:gpf-c2} 
f(w; \check{\lambda}) = \check{C} \cdot \check{d}^w \cdot 
\dfrac{\prod_{i=0}^{r-p-q-1}\varGamma\left(w+\frac{i}{r-p-q} \right)}{
\prod_{i=1}^{r-p-q} \varGamma\left(w+ \check{v}_i \right)}, 
\end{equation}
%%%%%%%%%%%%%%%%%%%%%%%%%%%%%%%%%%%%%%%%%%%%%%%%%%%%%%%%%%%%%%%%%%%%%%%
where $\check{C}$ is a positive constant, $\check{d}$ and 
$\check{v}_1, \dots, \check{v}_{r-p-q}$ are given by formulas 
\eqref{eqn:d-c} and \eqref{eqn:vi-c}. 
\end{proposition}
%%%%%%%%%%%%%%%%%%%%%%%%%%%%%%%%%%%%%%%%%%%%%%%%%%%%%%%%%%%%%%%%%%%%%%%
%%%%%%%%%%%%%%%%%%%%%%%%%%%%% begin proof %%%%%%%%%%%%%%%%%%%%%%%%%%%%%
{\it Proof}. 
The inclusions of multi-sets \eqref{eqn:incl-a} and \eqref{eqn:incl-b} 
are equivalent to the division relations  
%%%%%%%%%%%%%%%%%%%%%%%% eqn:division4 %%%%%%%%%%%%%%%%%%%%%%%%%%%%%%%%%
\begin{subequations} \label{eqn:division4}
\begin{align}
\prod_{i\in \bar{I}_p} \left(w+ \ts \frac{i+a}{p} \right)   
\prod_{i \in \bar{I}_q} \left(w+ \ts \frac{i+b}{q} \right) \,\,
& \Big | \, 
\prod_{j \in J_p} \left(w+ \ts \frac{j-a}{r-p} \right)   
\prod_{j \in J_q} \left(w+ \ts \frac{j-b}{r-q} \right), 
\label{eqn:division4a} \\
\prod_{i\in I_p} \left(w+ \ts \frac{i+a}{p} \right)    
\prod_{i \in I_q} \left(w+ \ts \frac{i+b}{q} \right) \,\,
& \Big | \, 
\prod_{j \in \bar{J}_p} \left(w+ \ts \frac{j-a}{r-p} \right)  
\prod_{j \in \bar{J}_q} \left(w+ \ts \frac{j-b}{r-q} \right),    
\label{eqn:division4b}
\end{align}    
\end{subequations}
%%%%%%%%%%%%%%%%%%%%%%%%%%%%%%%%%%%%%%%%%%%%%%%%%%%%%%%%%%%%%%%%%%%%%%%%
in $\C[w]$, respectively. 
Then division relation \eqref{eqn:division2} follows from formula 
\eqref{eqn:IJ} and relation \eqref{eqn:division4a}, while division 
relation \eqref{eqn:division3} is derived by multiplying relations 
\eqref{eqn:division4a} and \eqref{eqn:division4b} together. 
%%%%%%
\par
%%%%%%
Now that division relation \eqref{eqn:division2} is proved, the 
convention \eqref{eqn:vi-pq} recasts formula \eqref{eqn:h3} to   
%%%%%%%%%%%%%%%%%%%%%%%%%%%%% eqn:h5 %%%%%%%%%%%%%%%%%%%%%%%%%%%%%%%%%%%
\begin{equation} \label{eqn:h5} 
h(w;\lambda) = C_3 \cdot \tilde{d}^w \cdot 
\dfrac{\varGamma((r-p-q)w+1-a-b)}{\prod_{i=1}^{r-p-q} 
\varGamma\left(w+v_i\right)}.  
\end{equation}
%%%%%%%%%%%%%%%%%%%%%%%%%%%%%%%%%%%%%%%%%%%%%%%%%%%%%%%%%%%%%%%%%%%%%%%%
Replacing $w$ by $\hat{w}$ in formula \eqref{eqn:f-h}, where $\hat{w}$ 
is defined by formula \eqref{eqn:w-h}, we have 
%%%%%
\[
\begin{split}
f(w; \check{\lambda}) 
&= x^{r \check{w}-1} (1-x)^{(p+q-r)\check{w}+a+b} \, 
h(\hat{w}; \lambda) \\[2mm]
&= C_3 \cdot x^{r \check{w}-1} (1-x)^{(p+q-r)\check{w}+a+b} \,  
\tilde{d}^{\hat{w}} \cdot 
\dfrac{\varGamma((r-p-q) \hat{w} +1-a-b)}{\prod_{i=1}^{r-p-q} 
\varGamma\left(\hat{w} + v_i\right)} \\[2mm]
&= C_4 \cdot x^{r w} (1-x)^{(p+q-r) w} \, \tilde{d}^w \cdot 
\dfrac{\varGamma\left( (r-p-q)w \right)}{\prod_{i=1}^{r-p-q} 
\varGamma\left(w + \check{v}_i\right)},   
\end{split}
\]
%%%%%
where \eqref{eqn:h5} is used in the second equality, 
definitions \eqref{eqn:vi-c} and \eqref{eqn:w-h} are used in 
the third equality and $C_4$ is a real constant. 
Applying the multiplication formula \eqref{eqn:mult} to   
$\varGamma\left((r-p-q)w \right)$ yields 
%%%%%
\[
f(w; \check{\lambda}) 
= \check{C} \cdot \{ x^r (1-x)^{p+q-r} \, \tilde{d} \, (r-p-q)^{r-p-q} \}^w 
\cdot \dfrac{\prod_{i=0}^{r-p-q-1}
\varGamma\left(w + \frac{i}{r-p-q}\right)}{\prod_{i=1}^{r-p-q} 
\varGamma\left(w + \check{v}_i\right)}, 
\]
%%%%%
where $\check{C}$ is a real constant. 
This gives GPF \eqref{eqn:gpf-c2}, since 
$x^r (1-x)^{p+q-r} \, \tilde{d} \, (r-p-q)^{r-p-q} = \check{d}$ 
by definition \eqref{eqn:d-c}.  
To see $\check{C} > 0$, we look at formula \eqref{eqn:gpf-c2} 
for a large positive value of $w$. 
The right-hand side of it without constant factor $\check{C}$ is positive, 
while the left-hand side $f(w;\check{\lambda})$ is also positive since 
$\check{\lambda} \in \cF^-$.  
Here we used the fact that if $\mu \in \cF^-$ then $f(w; \mu) > 0$ 
for every large $w > 0$, which will be shown in the proof of 
Lemma \ref{lem:zeros}; see claim \eqref{eqn:f>0}. 
Thus $\check{C} > 0$. \hfill $\Box$
%%%%%%%%%%%%%%%%%%%%%%%%%%%%% end proof %%%%%%%%%%%%%%%%%%%%%%%%%%%%%%%%
%%%%%%%%%%%%%%%%%%%%%%%%%%%%% rem:division %%%%%%%%%%%%%%%%%%%%%%%%%%%%%
\begin{remark} \label{rem:division} 
In the situation of Proposition \ref{prop:division} we have 
$\check{a}$, $\check{b} \in \Q$, $\check{x}$ algebraic, and 
%%%%%%%%%%%%%% eqn:sum2 %%%%% eqn:estimate3 %%%%%%%%%%%%%%%%%%%%%%%%%%%% 
\begin{gather}
\check{v}_1 + \cdots + \check{v}_{r-p-q} = \frac{r-p-q-1}{2} = 
\frac{\check{r}-1}{2}, \label{eqn:sum2} \\
\check{v}_1, \dots, \check{v}_{r-p-q} \in \Q, \qquad 
-c \le \check{v}_1, \dots, \check{v}_{r-p-q} < 1-c, \qquad 
c := c(\lambda) = \frac{1-a-b}{r-p-q}. \label{eqn:estimate3}   
\end{gather}
%%%%%%%%%%%%%%%%%%%%%%%%%%%%%%%%%%%%%%%%%%%%%%%%%%%%%%%%%%%%%%%%%%%%%%%% 
Indeed, since $a$ and $b$ are rational by Lemma \ref{lem:ab} while 
$x$ is algebraic by assertion (1) of \cite[Theorem 2.3]{Iwasaki}, 
so are $\check{a}$, $\check{b}$ and $\check{x}$ via the definition 
of reciprocity \eqref{eqn:recip}.    
Moreover summation \eqref{eqn:sum2} comes from \eqref{eqn:sum} together 
with convention \eqref{eqn:vi-pq} and definition \eqref{eqn:vi-c},    
while condition \eqref{eqn:estimate3} follows from \eqref{eqn:estimate2} 
via definition \eqref{eqn:vi-c}.   
These remarks will play an important part in \S \ref{subsec:gpf-F}. 
\end{remark}
%%%%%%%%%%%%%%%%%%%%%%%%%%%%%%%%%%%%%%%%%%%%%%%%%%%%%%%%%%%%%%%%%%%%%%%%
%%%%%%%%%%%%%%%%%%%%%%%%%%%%% sec:dr %%%%%%%%%%%%%%%%%%%%%%%%%%%%%%%%%%%
\section{Division Relations} \label{sec:dr} 
%%%%%%%%%%%%%%%%%%%%%%%%%%%%%%%%%%%%%%%%%%%%%%%%%%%%%%%%%%%%%%%%%%%%%%%%
In Lemma \ref{lem:case1-6} we derived some arithmetical constraints 
on $\bp = (p, q; r)$ for $(\rA)$-solutions in $\cD^-$. 
In this section we exploit division relation \eqref{eqn:division3} 
to amplify this kind of study.    
In what follows let $\lambda = (p,q,r;a,b;x) \in \cD^-$ be an 
arbitrary $(\rA)$-solution.    
%%%%%%%%%%%%%%%%%%%%%%%%%%%%% lem:dr1 %%%%%%%%%%%%%%%%%%%%%%%%%%%%%%%%%%
\begin{lemma} \label{lem:dr1} 
If $p \not| r$ and $p \not| (r-p-q)$, then $p | 2r$ and $p | 2(r-p-q)$ 
with $p \ge 4$.  
%%%%%%%%%%%%%%%%%%%%%%%%%%%%%%%%%%%%%%%%%%%%%%%%%%%%%%%%%%%%%%%%%%%%%%%%
\end{lemma} 
%%%%%%%%%%%%%%%%%%%%%%%%%%%%% begin proof %%%%%%%%%%%%%%%%%%%%%%%%%%%%%%
{\it Proof}. 
If $p = 1$ we obviously have $p|r$. 
If $p = 2$ we have $p|(r-p-q)$ since $r-p-q$ must be even. 
Hereafter we suppose $p \ge 3$.  
Since $\prod_{i=0}^{p-1}\left( w+(i+a)/p \right)$ divides the 
right-hand side of \eqref{eqn:division3}, there exists a mapping 
$\phi : [p] \to [r-p] \sqcup [r-q]$ such that for each $i \in [p]$,  
%%%%%%%%%%%%%%%%%%%%%%%%%%%%% eqn:hm-ht %%%%%%%%%%%%%%%%%%%%%%%%%%%%%%%%
\begin{equation} \label{eqn:hm-ht}
\dfrac{i+a}{p} = 
\begin{cases}
\dfrac{\phi(i) - a}{r-p} \quad & \mbox{if} \quad \phi(i) \in [r-p] 
\quad (\mbox{$i$ is homogeneous, or of type $\rI$}), \\[4mm]
\dfrac{\phi(i) - b}{r-q} \quad & \mbox{if} \quad \phi(i) \in [r-q] 
\quad (\mbox{$i$ is heterogeneous, or of type $\rII$}). 
\end{cases}
\end{equation}
%%%%%%%%%%%%%%%%%%%%%%%%%%%%%%%%%%%%%%%%%%%%%%%%%%%%%%%%%%%%%%%%%%%%%%%%
Then a sequence $\tau = (\tau_0, \tau_1, \dots, \tau_{p-1})$ of 
symbols $\rI$ and $\rII$ is defined by assigning $\tau_i = \rI$ or 
$\tau_i = \rII$ if $i$ is homogeneous or heterogeneous respectively. 
There is a dichotomy: 
\begin{enumerate}
\item There exists at least one index $i \in [p-1]$ such that 
$\tau_i = \tau_{i+1}$. 
\item The sequence $\tau$ is interlacing, that is, 
$\tau_i \neq \tau_{i+1}$ for any index $i \in [p-1]$. 
\end{enumerate}
%%%%%
\par
%%%%%
We begin with case (1). 
In this case we have either $\tau_i = \tau_{i+1} = \rI$ or 
$\tau_i = \tau_{i+1} = \rII$. 
In the former subcase, it follows from formula \eqref{eqn:hm-ht} that  
%%%%%%
\[
\dfrac{i+a}{p} = \dfrac{\phi(i) - a}{r-p} \quad \mbox{and} \quad  
\dfrac{(i+1)+a}{p} = \dfrac{\phi(i+1) - a}{r-p}, \quad \mbox{so that} \quad 
\dfrac{1}{p} = \dfrac{\phi(i+1) - \phi(i)}{r-p}, 
\]
%%%%%% 
that is, $r-p = p \{\phi(i+1)-\phi(i)\}$, which implies $p|r$. 
In the latter subcase, similarly we have    
%%%%%%
\[
\dfrac{i+a}{p} = \dfrac{\phi(i) - b}{r-q} \quad \mbox{and} \quad  
\dfrac{(i+1)+a}{p} = \dfrac{\phi(i+1) - b}{r-q}, \quad \mbox{so that} \quad 
\dfrac{1}{p} = \dfrac{\phi(i+1) - \phi(i)}{r-q}, 
\]
%%%%%% 
that is, $r-q = p \{\phi(i+1)-\phi(i)\}$, which implies $p|(r-p-q)$.    
%%%%%
\par
%%%%%
We proceed to case (2) with $p = 3$. 
In this case we have either $\tau = (\rI, \rII, \rI)$ or 
$\tau = (\rII, \rI, \rII)$. 
In the former subcase, it follows from \eqref{eqn:hm-ht} 
and $\tau_0 = \tau_2 = \rI$ that 
%%%%%%
\[
\dfrac{a}{p} = \dfrac{\phi(0) - a}{r-p} \quad \mbox{and} \quad  
\dfrac{2+a}{p} = \dfrac{\phi(2) - a}{r-p}, \quad \mbox{so that} \quad 
\dfrac{2}{p} = \dfrac{\phi(2) - \phi(0)}{r-p}, 
\]
%%%%%% 
that is, $2(r-p) = p \{\phi(2)-\phi(0)\}$, which implies $p|2 r$, 
but since $p = 3$ is odd, we must have $p|r$. 
In the latter subcase, a similar reasoning with $\tau_0 = \tau_2 = 
\rII$ leads to $2(r-q) = p \{\phi(2)-\phi(0)\}$, which implies 
$p|2(r-p-q)$, but since $p = 3$ is odd, we must have $p|(r-p-q)$.   
%%%%%%
\par
%%%%%%
Finally we consider case (2) with $p \ge 4$. 
In this case we have either $(\tau_0, \tau_1, \tau_2, \tau_3) = 
(\rI, \rII, \rI, \rII)$ or $(\tau_0, \tau_1, \tau_2, \tau_3) = 
(\rII, \rI, \rII, \rI)$.  
In the former subcase, it follows from \eqref{eqn:hm-ht} that 
%%%%%%%%%%%%%%
\begin{align*}
\dfrac{a}{p} &= \dfrac{\phi(0) - a}{r-p} \quad \mbox{and} \quad  
\dfrac{2+a}{p} = \dfrac{\phi(2) - a}{r-p}, \quad \mbox{so that} \quad 
\dfrac{2}{p} = \dfrac{\phi(2) - \phi(0)}{r-p}, \\[4mm]
\dfrac{1+a}{p} &= \dfrac{\phi(1)-b}{r-q} \quad \mbox{and} \quad  
\dfrac{3+a}{p} = \dfrac{\phi(3)-b}{r-q}, \quad \mbox{so that} \quad 
\dfrac{2}{p} = \dfrac{\phi(3) - \phi(1)}{r-q},
\end{align*}
%%%%%%%%%%%%%
that is, $2(r-p) = p \{\phi(2)-\phi(0)\}$ and 
$2(r-q) = p \{\phi(3)-\phi(1)\}$, which imply 
$p|2 r$ and $p|2(r-p-q)$. 
In the latter subcase, a similar reasoning leads to 
$2(r-p) = p \{\phi(3)-\phi(1)\}$ and $2(r-q) = p \{\phi(2)-\phi(0)\}$, 
which again imply $p|2 r$ and $p|2(r-p-q)$. 
Thus if $p \not| r$ and $p \not| (r-p-q)$, then we must be in case (2) 
with $p \ge 4$, which forces $p | 2r$ and $p | 2(r-p-q)$. \hfill $\Box$ 
%%%%%%%%%%%%%%%%%%%%%%%%%%%%% end proof %%%%%%%%%%%%%%%%%%%%%%%%%%%%%%%%
%%%%%%%%%%%%%%%%%%%%%%%%%%%%% lem:dr2 %%%%%%%%%%%%%%%%%%%%%%%%%%%%%%%%%%
\begin{lemma} \label{lem:dr2} 
If $p \not| r$ and $p \not|(r-p-q)$, then 
\begin{enumerate} 
\item there exist positive integers $s$, $t$ and $k$ such that 
%%%%%%%%%%%%%%%%%%%%%%%%%%%%% eqn:stk %%%%%%%%%%%%%%%%%%%%%%%%%%%%%%%%%
\begin{equation} \label{eqn:stk}
p = 2^{k+1} \, s; \qquad p|q; \qquad r = 2^k \, t; 
\qquad s, \, t : \mbox{odd}; \qquad s|t,   
\end{equation}
%%%%%%%%%%%%%%%%%%%%%%%%%%%%%%%%%%%%%%%%%%%%%%%%%%%%%%%%%%%%%%%%%%%%%%%%
\item we must also have $q \not| r$ and $q \not|(r-p-q)$. 
\end{enumerate}
\end{lemma}
%%%%%%%%%%%%%%%%%%%%%%%%%%%%% begin proof %%%%%%%%%%%%%%%%%%%%%%%%%%%%%%
{\it Proof}. 
We use Lemma \ref{lem:dr1}, upon writing $p = 2^i \, s$, $r = 2^k \, t$, 
$r-q = 2^j \, u$ with $i,j,k \in \Z_{\ge0}$ and odd $s, t, u \in \N$.  
Division relation $p| 2 r$ implies $i \le k+1$ and $s|t$, while 
$p \not| r$ yields $i \ge k+1$ and hence $i = k+1$. 
In a similar manner the division relation $p| 2 (r-q)$ implies 
$i \le j+1$ and $s|u$, while $p \not| (r-q)$ yields $i \ge j+1$ 
and hence $i = j+1$. 
In summary we have 
%%%%%
\[
p = 2^{k+1} \, s; \quad r = 2^k \, t; \quad 
r-q = 2^k \, u; \quad k \in \Z_{\ge0}; \quad 
s, \,t, \,u : \mbox{odd}; \quad s|t, \quad s|u.  
\]
%%%%%  
Since $t$ and $u$ are odd, we can write $t = 2 t'+1$ and $u = 2 u'+1$ 
with $t', u' \in \Z_{\ge0}$, so that $s|t$ and $s|u$ imply 
$s|(t-u) = s|2(t'-u')$, which in turn yields $s|(t'-u')$ since 
$s$ is also odd.  
It then implies $p|q$ because $p = 2^{k+1} \, s$ and $q = r-(r-q) = 
2^k (t-u) = 2^{k+1} (t'-u')$. 
Note that $p$ is even by $p = 2^{k+1} \, s$ with $k \ge 0$, 
so $q$ is even too by $p|q$ and $r$ is also even, since $r-p-q$ is even.  
Thus we have $k \ge 1$ and all the conditions in \eqref{eqn:stk} 
have been proved. 
%%%%%
\par
%%%%%
Next we show that $q \not| r$ and $q \not|(r-p-q)$ by contradiction. 
Indeed, if $q|r$ then the division relation $p|q$ in condition 
\eqref{eqn:stk} yields $p|r$ contrary to the assumption $p \not|r$, 
while if $q|(r-p-q)$ then $p|q$ gives $p|(r-p-q)$ contrary to the 
assumption $p \not|(r-p-q)$. \hfill $\Box$ \par\medskip 
%%%%%%%%%%%%%%%%%%%%%%%%%%%%% end proof %%%%%%%%%%%%%%%%%%%%%%%%%%%%%%%%
Note that Lemmas \ref{lem:dr1} and \ref{lem:dr2} remain true 
if the roles of $p$ and $q$ are exchanged. 
%%%%%%%%%%%%%%%%%%%%%%%%%%%%% prop:dr %%%%%%%%%%%%%%%%%%%%%%%%%%%%%%%%%%
\begin{proposition} \label{prop:dr} 
For any $(\rA)$-solution $\lambda = (p,q,r;a,b;x) \in \cD^-$ the integer 
vector $\bp = (p,q;r) \in D^-_{\rA}$ must satisfy division relations 
\eqref{eqn:dr}.   
\end{proposition}
%%%%%%%%%%%%%%%%%%%%%%%%%%%%%%%%%%%%%%%%%%%%%%%%%%%%%%%%%%%%%%%%%%%%%%%%
%%%%%%%%%%%%%%%%%%%%%%%%%%%%% begin proof %%%%%%%%%%%%%%%%%%%%%%%%%%%%%%
{\it Proof}. 
To prove the lemma by contradiction, suppose the contrary that 
%%%%%%%%%%%%%
\[
(\, p \not| r \quad \mbox{and} \quad p \not| (r-p-q) \,) 
\quad \mbox{or} \quad 
(\, q \not| r \quad \mbox{and} \quad q \not| (r-p-q) \,).
\]
%%%%%%%%%%%%% 
By symmetry we may take the former condition in the ``or" sentence 
above, but the latter condition also follows from assertion (2) of 
Lemma \ref{lem:dr2}, 
so that we are led to the ``and" sentence:   
%%%%%%%%%%%%%
\[
(\, p \not| r \quad \mbox{and} \quad p \not| (r-p-q) \,) 
\quad \mbox{and} \quad 
(\, q \not| r \quad \mbox{and} \quad q \not| (r-p-q) \,).
\]
%%%%%%%%%%%%%
By a part of condition \eqref{eqn:stk} we have $p|q$ and likewise  
$q|p$ upon exchanging the role of $p$ and $q$. 
Hence $p=q$ and condition \eqref{eqn:stk} becomes 
%%%%%%%%%%%%%%%%%%%%%%%%%%%%% eqn:stk2 %%%%%%%%%%%%%%%%%%%%%%%%%%%%%%%%%
\begin{equation} \label{eqn:stk2}
p = q = 2^{k+1} \, s; \qquad r = 2^k \, t; 
\qquad k \in \N, \quad s, \, t : \mbox{odd}; \qquad s|t,   
\end{equation}
%%%%%%%%%%%%%%%%%%%%%%%%%%%%%%%%%%%%%%%%%%%%%%%%%%%%%%%%%%%%%%%%%%%%%%%%
\par
%%%%%
A look at Table \ref{tab:ab} shows that if $p = q$ then we must have 
$p|r$ in cases (1)--(5), while currently we have the contrary 
$p \not| r$. 
Thus we must be in case 6 of Table \ref{tab:ab}, so that $\Z$-linear 
equation \eqref{eqn:ls-c6a} must be satisfied. 
It follows from formula \eqref{eqn:stk2} that $r-p-q = 2^k s(t_s-4)$, 
$r-p = 2^k s(t_s-2)$, $q = 2^{k+1} s$ and 
$((r-p)(r-p-q))_p = 2^{k-1} s(t_s-2)(t_s-4)$, so equation 
\eqref{eqn:ls-c6a} is equivalent to 
%%%%%%
\[
(t_s-4) i + (t_s-2) i' + 2 j = \frac{1}{2} (t_s-2)(t_s-4), \qquad 
i, i', j \in \Z_{\ge0},  
\]
%%%%%%
where we used the notation introduced in item (1) of Remark \ref{rem:ab}.  
Notice that the left-hand side above is an integer, while the 
right-hand side is a half-integer (not an integer), because $s$ and 
$t$ are odd integers with $s|t$ so $t_s = t/s$ is also an odd integer. 
This contradiction shows that our initial assumption is false 
and condition \eqref{eqn:dr} must be true. 
\hfill $\Box$ \par\medskip 
%%%%%%%%%%%%%%%%%%%%%%%%%%%%% end proof %%%%%%%%%%%%%%%%%%%%%%%%%%%%%%%%%
Division relation \eqref{eqn:dr} in Theorem \ref{thm:ab-rf} is now 
established by Proposition \ref{prop:dr}. 
If $\bp = (p,q;r) \in D^-_{\rA}$ then $\check{r} := r-p-q \in 2 \N$.  
Conversely we have the following.  
%%%%%%%%%%%%%%%%%%%%%%%%%%%%% lem:dr-bound %%%%%%%%%%%%%%%%%%%%%%%%%%%%%%
\begin{lemma} \label{lem:dr-bound}  
Given any $\check{r} \in 2 \N$, there are only a finite number of triples 
$\bp = (p,q;r) \in D^-_{\rA}$ that satisfy $r-p-q = \check{r}$ and the 
division relations \eqref{eqn:dr}; any such $\bp$ must be bounded by 
%%%%%%%%%%%%%%%%%%%%%%%%%%%%% eqn:dr-bound %%%%%%%%%%%%%%%%%%%%%%%%%%%%%%
\begin{equation} \label{eqn:dr-bound}
1 \le p, \, q \le 3 \check{r}, 
\qquad 2 \le p+q \le 5 \check{r}, 
\qquad 4 \le r \le 6 \check{r}.  
\end{equation}
%%%%%%%%%%%%%%%%%%%%%%%%%%%%%%%%%%%%%%%%%%%%%%%%%%%%%%%%%%%%%%%%%%%%%%%%%  
\end{lemma}
%%%%%%%%%%%%%%%%%%%%%%%%%%%%%%%%%%%%%%%%%%%%%%%%%%%%%%%%%%%%%%%%%%%%%%%%%
%%%%%%%%%%%%%%%%%%%%%%%%%%%%% begin proof %%%%%%%%%%%%%%%%%%%%%%%%%%%%%%%
{\it Proof}. 
Since $p, \, q \in \N$ and $\check{r} \in 2 \N$, it is evident that 
$p, \, q \ge 1$ and $r = \check{r}+p+q \ge 4$. 
Division relation \eqref{eqn:dr} is equivalent to the condition that 
$( p|\check{r} \,\, \mbox{or} \,\, p|(\check{r}+q) )$ and 
$( q|\check{r} \,\, \mbox{or} \,\, q|(\check{r}+p) )$, which is divided 
into four cases: (i) $p|\check{r}$ and $q|\check{r}$; 
(ii) $p|\check{r}$ and $q|(\check{r}+p)$;   
(iii) $p|(\check{r}+q)$ and $q|\check{r}$; 
(iv) $p|(\check{r}+q)$ and $q|(\check{r}+p)$. 
In case (i) we have $p, \, q \le \check{r}$ and $r = \check{r}+p+q \le 3 
\check{r}$. 
In case (ii) we have $p \le \check{r}$, $q \le \check{r}+p \le 2 \check{r}$ 
and $r = \check{r}+p+q \le 4 \check{r}$. 
Case (iii) is similar to case (ii). 
In case (iv) there exist $i$, $j \in \N$ such that 
$\check{r}+q = i p$ and $\check{r}+p = j q$. 
Note that $i j \ge 2$, for otherwise $i = j = 1$ would imply 
$\check{r} = p-q = - \check{r} = 0$, a contradiction to $\check{r} \ge 4$. 
Thus the two equations above for $(p, q)$ are uniquely settled as 
$p = l \, \check{r}$, $q = m \, \check{r}$ and hence  
$r = \check{r}+p+q = n \, \check{r}$, where 
%%%%%% 
\[
l = \dfrac{j+1}{i j-1}, \qquad m = \dfrac{i+1}{i j-1}, \qquad 
n = 1+l+m = \dfrac{(i+1)(j+1)}{i j-1} \qquad (i, j \in \N; \, i j \ge 2). 
\] 
%%%%%%
To estimate these numbers we may assume $i \ge j$ and so 
$i \ge 2$ and $j \ge 1$ by symmetry.  
Then  
%%%%%%
\[
l \le \frac{j+1}{2j-1} = \frac{1}{2} + \frac{3}{2(2j-1)} \le 
\frac{1}{2} + \frac{3}{2} = 2, \qquad  
m \le \frac{i+1}{i-1} = 1 + \frac{2}{i-1} \le 1+2=3,  
\]
%%%%%% 
and $n = 1+l+m \le 1+ 2+3 =6$. 
This establishes the bound \eqref{eqn:dr-bound}. 
\hfill $\Box$ \par\medskip
%%%%%%%%%%%%%%%%%%%%%%%%%%%%% end proof %%%%%%%%%%%%%%%%%%%%%%%%%%%%%%%%%
For example, those $\bp = (p,q;r) \in D_{\rA}^-$ with 
$\check{r} = 2$ and $p \ge q$ are exactly $\bp = (1,1;4)$, $(2,1;5)$, 
$(2,2;6)$, $(3,1;6)$, $(4,2;8)$, among which only $(2,1;5)$ leads to 
no solutions (see \cite[Table 1]{Iwasaki}). 
%%%%%%%%%%%%%%%%%%%%%%%%%%%%%% sec:SW %%%%%%%%%%%%%%%%%%%%%%%%%%%%%%%%%%%
\section{The South-West Domain} \label{sec:SW}
%%%%%%%%%%%%%%%%%%%%%%%%%%%%%%%%%%%%%%%%%%%%%%%%%%%%%%%%%%%%%%%%%%%%%%%%%
%%%%%%%%%%%%%%%%%%%%%%%%%%%%%% subsec:cfcr-F %%%%%%%%%%%%%%%%%%%%%%%%%%%%
\subsection{Coming from Contiguous Relations} \label{subsec:cfcr-F}
%%%%%%%%%%%%%%%%%%%%%%%%%%%%%%%%%%%%%%%%%%%%%%%%%%%%%%%%%%%%%%%%%%%%%%%%
It is important to think of the linear independence of $f(w; \lambda)$ 
and $\tilde{f}(w; \lambda)$ over the rational function field $\C(w)$, 
where $\tilde{f}(w; \lambda)$ is defined right after formula \eqref{eqn:3trf}. 
This issue was already discussed in \cite[\S 7]{Iwasaki} on the 
domain $\cD \cup \cI \cup \cE$, where the equivalence of Problems 
$\rI$ and $\rII$ (see \cite[Theorem 2.1]{Iwasaki}) made it relatively tractable.   
Without such an advantage, the discussion on $\cF^-$ should be more 
elaborate and require some function-theoretic preliminaries.    
%%%%%%%%%%%%%%%%%%%%%%%%%%%%%% lem:sin %%%%%%%%%%%%%%%%%%%%%%%%%%%%%%%%%
\begin{lemma} \label{lem:sin}
If $\alpha_1, \dots, \alpha_k$ and $\beta_1, \dots, \beta_k$ are 
real numbers with $\alpha_1 \ge \cdots \ge \alpha_k > 0$, then   
for any integer $m \ge \beta_1$ there exists a positive number 
$\rho_m \in I_m$ such that 
%%%%%%%%%%%%%%%%%%%%%%%%%%%%%% eqn:sin1 %%%%%%%%%%%%%%%%%%%%%%%%%%%%%%%%%
\begin{equation} \label{eqn:sin1}
\prod_{j=1}^k \dfrac{1}{|\sin \pi(\alpha_j w + \beta_j)|} \le C 
\qquad \mbox{on the circle} \quad |w| = \rho_m, 
\end{equation}
%%%%%%%%%%%%%%%%%%%%%%%%%%%%%%%%%%%%%%%%%%%%%%%%%%%%%%%%%%%%%%%%%%%%%%%%%
where $C$ is a positive constant independent of $m$ and $I_m$ is 
an open interval defined by  
%%%%%
\[
I_m := \left(\frac{m-\beta_1}{\alpha_1}, \,\, 
\frac{m+1-\beta_1}{\alpha_1}\right) \qquad (m \in \Z).  
\]
%%%%%
\end{lemma}
%%%%%%%%%%%%%%%%%%%%%%%%%%%%%%%%%%%%%%%%%%%%%%%%%%%%%%%%%%%%%%%%%%%%%%%%%
%%%%%%%%%%%%%%%%%%%%%%%%%%%%%%% begin proof %%%%%%%%%%%%%%%%%%%%%%%%%%%%%
{\it Proof}. 
Note that $I_m$ is an open interval of width $1/\alpha_1$, whose 
endpoints are a pair of consecutive zeros of $\sin \pi(\alpha_1 w+ \beta_1)$. 
For each $j = 2, \dots,k$, all the zeros of $\sin \pi(\alpha_j w + \beta_j)$ 
form an arithmetic progression 
$w_{j,m}:=\frac{m - \beta_j}{\alpha_j}$ $(m \in \Z)$ of common difference 
$1/\alpha_j \ge 1/\alpha_1$, so that $I_m$ can contain at most one zero of 
$\sin \pi(\alpha_j w + \beta_j)$. 
Thus $I_m$ can contain at most $k-1$ zeros of $\prod_{j=2}^k 
\sin \pi(\alpha_j w+\beta_j)$, which partition $I_m$ into at most 
$k$ subintervals. 
Among them let $J_m$ be a subinterval of the largest width.   
By pigeon hole principle the width of $J_m$ cannot be smaller than 
$\frac{1}{k \alpha_1}$, so that the midpoint $\rho_m$ of $J_m$ is 
at least $\delta := \frac{1}{2 k \alpha_1}$ distant from all the 
zeros of $\prod_{j=1}^k \sin \pi(\alpha_j w + \beta_j)$. 
For $j = 1, \dots,k$ and $m \in \Z$ let $D_{j,m}$ be the open disk 
of radius $\delta$ with center at $w_{j,m}$. 
Then it is not hard to see that there exists a positive constant 
$C_j$ such that 
%%%%%
\[
\dfrac{1}{|\sin \pi(\alpha_j w + \beta_j)|} \le C_j \qquad 
\mbox{for all} \quad w \in \C \setminus D_j, 
\quad \mbox{with} \quad D_j := \bigcup_{m \in \Z} D_{j,m}. 
\]
%%%%%
Accordingly, if we set $C := C_1 \cdots C_k$ then we have  
%%%%%%%%%%%%%%%%%%%%%%%%%% eqn:sin2 %%%%%%%%%%%%%%%%%%%%%%%%%%%%%%%%%%%%
\begin{equation} \label{eqn:sin2}
\prod_{j=1}^k \dfrac{1}{|\sin \pi(\alpha_j w + \beta_j)|} \le C 
\qquad \mbox{for all} \quad w \in K := \C \setminus \bigcup_{j=1}^k D_j. 
\end{equation}
%%%%%%%%%%%%%%%%%%%%%%%%%%%%%%%%%%%%%%%%%%%%%%%%%%%%%%%%%%%%%%%%%%%%%%%%  
Notice that $\rho_m \in K$ for any $m \in \Z$. 
If $m \ge \beta_1$ then $\rho_m > 0$ and the circle $|w| = \rho_m$ 
is contained in $K$, hence estimate \eqref{eqn:sin2} leads to 
estimate \eqref{eqn:sin1}. \hfill $\Box$
%%%%%%%%%%%%%%%%%%%%%%%% end proof %%%%%%%%%%%%%%%%%%%%%%%%%%%%%%%%%%%%%%
%%%%%%%%%%%%%%%%%%%%%%%% lem:hg-estimate %%%%%%%%%%%%%%%%%%%%%%%%%%%%%%%%
\begin{lemma} \label{lem:hg-estimate}
If $\lambda = (p,q,r;a,b;x)$ lies in domain \eqref{eqn:real2} $($with 
$b \not\in \Z$ when $q = 0$, $r$$)$, then there exists an infinite 
sequence of positive numbers $\{ \rho_m \}_{m \ge m_0}$ such that 
for any integer $m \ge m_0$,  
%%%%%%%%%%%%%%%%%%%%%%%% eqn:hg-estimate %%%%%%%%%%%%%%%%%%%%%%%%%%%%%%%%
\begin{equation} \label{eqn:hg-estimate}
\sigma m + \tau < \rho_m < \sigma(m+1) + \tau, \qquad 
|f(w; \lambda)| \le M_1 e^{c_1 \rho_m} \quad \mbox{on the circle} 
\quad |w| = \rho_m,   
\end{equation}
%%%%%%%%%%%%%%%%%%%%%%%%%%%%%%%%%%%%%%%%%%%%%%%%%%%%%%%%%%%%%%%%%%%%%%%%%
where $M_1$, $c_1$, $\sigma > 0$ and $\tau \in \R$ are independent of 
$m$ and $m_0$ is an integer with $\sigma m_0 + \tau \ge 0$. 
\end{lemma}
%%%%%%%%%%%%%%%%%%%%%%%%%%%%%%%%%%%%%%%%%%%%%%%%%%%%%%%%%%%%%%%%%%%%%%%%%
%%%%%%%%%%%%%%%%%%%%%%%% begin proof %%%%%%%%%%%%%%%%%%%%%%%%%%%%%%%%%%%%
{\it Proof}.  
The hypergeometric function admits Euler's contour integral representation 
%%%%%
\[
\hgF(\alpha,\beta;\gamma;z) = - \dfrac{e^{-\pi i \gamma}}{4 \sin \pi \beta 
\cdot \sin \pi(\gamma-\beta) \cdot B(\beta, \, \gamma-\beta)} 
\int_{\wp} t^{\beta-1} (1-t)^{\gamma-\beta-1} (1-z t)^{-\alpha} d t,     
\]
%%%%%
along a Pochhammer loop $\wp$ around $t=0$ and $t=1$. 
Notice that  
%%%%% 
\[
\begin{split}
\dfrac{1}{B(\beta, \gamma-\beta)} 
&= (1-\gamma) \cdot 
\dfrac{\sin \pi \beta \cdot \sin \pi(\gamma-\beta)}{\pi \sin \pi \gamma} \cdot 
B(1-b, \, 1-(\gamma-\beta)) \\
&= - \dfrac{(1-\gamma) e^{\pi i \gamma}}{4 \pi \cdot \sin \pi \gamma} \int_{\wp} 
t^{-\beta} (1-t)^{-(\gamma-\beta)} d t,  
\end{split}
\]
%%%%%
where the first equality follows from the reflection formula 
for the beta function: 
%%%%%
\[
B(\alpha, \beta) B(1-\alpha, 1-\beta) = 
\frac{\pi \sin \pi(\alpha+\beta)}{(1-\alpha-\beta) \cdot 
\sin \pi \alpha \cdot \sin \pi \beta}, 
\]
%%%%%
while the second equality stems from Euler's integral representation 
of the beta function along the Pochhammer loop $\wp$.  
Putting these together we have  
%%%%%
\[
\begin{split}
\hgF(\alpha, \beta; \gamma; z) &= \dfrac{1-\gamma}{16 \pi \cdot \sin \pi \beta 
\cdot \sin \pi \gamma \cdot \sin \pi(\gamma-\beta)} \\
&\phantom{=} \times 
\left(\int_{\wp} t^{-\beta}(1-t)^{-(\gamma-\beta)} d t\right) 
\left(\int_{\wp} t^{\beta-1} (1-t)^{\gamma-\beta-1} (1-z t)^{-\alpha} d t \right),  
\end{split}
\]
%%%%%
which evaluated at $(\alpha, \beta; \gamma; z) = (p w + a, q w + b; r w; x)$ 
yields $f(w; \lambda) = \psi_1(w) \psi_2(w)$ with 
%%%%%
\[
\begin{split}
\psi_1(w) &= \dfrac{1- r w}{16 \pi \cdot \sin \pi (q w+b) 
\cdot \sin(\pi r w) \cdot \sin \pi((r-q)w-b)}, \\[1mm]
\psi_2(w) &= \left(\int_{\wp} t^{- q w-b}(1-t)^{-(r-q)w+b} d t\right) 
\left(\int_{\wp} t^{q w+b-1} (1-t)^{(r-q)w-b-1} (1-x t)^{-p w-a} d t \right). 
\end{split}
\]
%%%%%
\par
%%%%%
We can apply Lemma \ref{lem:sin} to the first factor $\psi_1(w)$. 
For some constants $\sigma > 0$ and $\tau \in \R$ there exists an infinite 
sequence $\rho_m \in (\sigma m + \tau, \, \sigma(m+1) + \tau)$ such that 
$|\psi_1(w)| = O(\rho_m)$ on the circle $|w| = \rho_m$ as $m \to + \infty$. 
For the second factor $\psi_2(w)$ it is easy to see that there exists 
a constant $c_2$ such that $|\psi_2(w)| = O(e^{c_2|w|})$ as $|w| \to +\infty$, 
because the integrands in $\psi_2(w)$ admit similar exponential estimates,  
uniform in $t \in \wp$, as $|w| \to +\infty$.  
Now \eqref{eqn:hg-estimate} follows readily. \hfill $\Box$ 
%%%%%%%%%%%%%%%%%%%%%%%% end proof %%%%%%%%%%%%%%%%%%%%%%%%%%%%%%%%%%%%%%
%%%%%%%%%%%%%%%%%%%%%%%% lem:zeros %%%%%%%%%%%%%%%%%%%%%%%%%%%%%%%%%%%%%% 
\begin{lemma} \label{lem:zeros}
For any data $\lambda = (p,q,r;a,b;x) \in \cF^-$ the following hold.       
\begin{enumerate}
\item Any pole of $f(w;\lambda)$ is simple and lies in the arithmetic 
progression $\{w_j := -j/r\}_{j=0}^{\infty}$. 
Conversely, $f(w; \lambda)$ actually has a pole at $w = w_j$ for every 
sufficiently large integer $j$, in particular it has infinitely many poles.        
\item $f(w;\lambda)$ has infinitely many zeros.   
\item $f(w;\lambda)$ and $\tilde{f}(w;\lambda)$ are linearly independent 
over the rational function field $\C(w)$. 
\end{enumerate} 
\end{lemma}
%%%%%%%%%%%%%%%%%%%%%%%%%%%%%%%%%%%%%%%%%%%%%%%%%%%%%%%%%%%%%%%%%%%%%%%%%
%%%%%%%%%%%%%%%%%%%%%%%% begin proof %%%%%%%%%%%%%%%%%%%%%%%%%%%%%%%%%%%%
{\it Proof}. Assertion (1). 
It is evident from definition \eqref{eqn:f} that every pole of 
$f(w;\lambda)$ is simple and lies in the sequence $\{w_j\}_{j=0}^{\infty}$. 
By \cite[Lemma 4.1]{Iwasaki} we have   
%%%%%%%%%%%%%%%%%%%%%%%%%%%%%% eqn:residue %%%%%%%%%%%%%%%%%%%%%%%%%%%%
\[ 
\underset{\scriptstyle w = w_j}{\mathrm{Res}} \, f(w; \lambda) = C_j \cdot 
\hgF(a_j, \, b_j; \, j+2; \, x) \qquad (j = 0,1,2,\dots),  
\]
%%%%%%%%%%%%%%%%%%%%%%%%%%%%%%%%%%%%%%%%%%%%%%%%%%%%%%%%%%%%%%%%%%%%%%%
where $a_j := pw_j+j+a+1$, $b_j := qw_j+j+b+1$ and 
%%%%%%%%%%%%%%%%%%%%%%%%%%%%%% eqn:Cj %%%%%%%%%%%%%%%%%%%%%%%%%%%%%%%%%
\[
C_j := \frac{(-1)^j}{r} \cdot
\frac{(p w_j+a)_{j+1} \, (q w_j+b)_{j+1}}{j! \, (j+1)!} \, x^{j+1}. 
\]
%%%%%%%%%%%%%%%%%%%%%%%%%%%%%%%%%%%%%%%%%%%%%%%%%%%%%%%%%%%%%%%%%%%%%%%
Since $\lambda \in \cF^-$, we have $p < 0$, $q < 0$, $r > 0$ and 
$0 < x < 1$, so there exists an integer $j_0$ such that 
$p w_j + a > 0$, $q w_j + b > 0$ and hence $(-1)^j C_j > 0$ 
for every $j \ge j_0$. 
Notice also that if $j \ge j_0$ then $a_j > p w_j + a > 0$, 
$b_j > q w_j + b > 0$ and thus $\hgF(a_j, b_j; j+2; x) > 0$.  
Therefore $f(w)$ actually has a simple pole with a non-vanishing 
residue at $w = w_j$ for every $j \ge j_0$.  
%%%%%
\par
%%%%%
Assertion (2). 
Suppose the contrary that $f(w;\lambda)$ has at most finitely many zeros. 
Then it follows from assertion (1) that $u(w) := f(w;\lambda)/\varGamma(r w)$ 
is an entire holomorphic function with at most finitely many zeros. 
Here we have a uniform estimate $|1/\varGamma(r w)| = 
O\left(e^{r |w| (\log|w| + c_3)}\right)$ as $|w| \to \infty$, 
which follows from Stirling's formula: for any $\varepsilon \in 
(0, \, \pi)$ one has uniformly 
%%%%%%%%%%%%%%%%%%%%%%%%%% eqn:stirling %%%%%%%%%%%%%%%%%%%%%%%%%%%%%%%%%
\begin{equation} \label{eqn:stirling}
\dfrac{1}{\varGamma(w)} \sim 
\dfrac{1}{\sqrt{2\pi}} \, w^{1/2} \, e^{w(1-\log w)} \qquad 
\mbox{as} \quad w \to \infty \quad \mbox{in} \quad |\arg w| < 
\pi - \varepsilon, 
\end{equation}
%%%%%%%%%%%%%%%%%%%%%%%%%%%%%%%%%%%%%%%%%%%%%%%%%%%%%%%%%%%%%%%%%%%%%%%%%
combined with the reflection formula \eqref{eqn:refl} for 
the gamma function.  
By this estimate and Lemma \ref{lem:hg-estimate} we have  
$|u(w)| \le M_2 \, e^{r \rho_m (\log \rho_m + c_4)}$ on the circle 
$|w| = \rho_m$ for every $m \ge m_0$. 
Given any $z \in \C$ with sufficiently large $|z|$, take an $m \ge m_0$ 
so that $\rho_{m-1} < |z| \le \rho_m$. 
Since $u(w)$ is entire and $\rho_m < \rho_{m-1} + 2 \sigma < |z| + 2 \sigma$ 
by estimate \eqref{eqn:hg-estimate}, the maximum principle yields 
%%%
\[
|u(w)| \le M_2 \, e^{r \rho_m (\log \rho_m + c_4)} 
\le M_2 \, e^{r (|z|+ 2 \sigma) \{\log (|z| + 2 \sigma) + c_4\}} 
\le M_3 \, e^{r |z| (\log |z| + c_5)},  
\]
%%% 
which means that $u(w)$ is an entire function of order at most $1$.  
Since $u(w)$ has at most finitely many zeros, Hadamard's factorization 
theorem allows us to write $u(w) = u_0(w) \, e^{c_6 w}$, where $u_0(w)$ 
is a nonzero polynomial (or constant) and $c_6 \in \C$.  
As a polynomial, $u_0(w) \sim M_4 \, w^k$ as $w \to \infty$ 
for some $M_4 \in \C^{\times}$ and $k \in \Z_{\ge0}$. 
It follows from estimates \eqref{eqn:hg-estimate} and 
\eqref{eqn:stirling} that 
%%%%
\[
\begin{split}
1 &= \left|\dfrac{u(\rho_m)}{u_0(\rho_m) \, e^{c_6 \rho_m}} \right| 
= \left|\dfrac{f(\rho_m; \lambda)}{u_0(\rho_m) \, e^{c_6 \rho_m}} 
\cdot \frac{1}{\varGamma(r \rho_m)} \right| \lesssim  
\dfrac{M_1 \, e^{c_1 \rho_m} \cdot (r \rho_m)^{1/2} 
e^{r \rho_m \{1-\log(r \rho_m) \}}}{|M_4| \, \rho_m^k \cdot e^{\rRe(c_6) \rho_m} \cdot 
\sqrt{2\pi}} \\
&= M_5 \cdot \rho_m^{\frac{1}{2}-k} \cdot e^{\rho_m (c_7 - r \log \rho_m)} 
\to 0 \qquad \mbox{as} \quad m \to \infty, 
\end{split}
\]
%%%% 
with $M_5 := \frac{M_1}{|M_4|}\sqrt{\frac{r}{2\pi}}$ and 
$c_7 := c_1 - \rRe(c_6) + r - r \log r$, where $s_m \lesssim t_m$ means  
$\displaystyle \limsup_{m \to \infty} \frac{s_m}{t_m} \le 1$.  
This contradiction shows that $f(w;\lambda)$ must have infinitely many zeros. 
%%%%
\par
%%%%
Assertion (3) can be established by modifying the proof of 
\cite[Proposition 7.1]{Iwasaki}. 
Suppose the contrary that $f(w; \lambda)$ and $\tilde{f}(w; \lambda)$ are 
linearly dependent over $\C(w)$, so that there exists a rational function 
$T(w) \in \C(w)$ such that $\tilde{f}(w; \lambda) = T(w) f(w;\lambda)$, 
because $f(w; \lambda)$ does not vanish identically. 
Then as in the proof of \cite[Proposition 7.1]{Iwasaki} we have 
%%%%%%%%%%%%%%%%%%%%%%%%% eqn:f-f1 %%%%%%%%%%%%%%%%%%%%%%%%%%%%%%%%%%%%%%
\begin{equation} \label{eqn:f-f1}
f(w; \lambda) f_1(w; \lambda) = (1-x)^{(r-p-q)w-a-b-1},  
\end{equation}
%%%%%%%%%%%%%%%%%%%%%%%%%%%%%%%%%%%%%%%%%%%%%%%%%%%%%%%%%%%%%%%%%%%%%%%%%
where $f_1(w; \lambda)$ is defined by    
%%%%%%
\[
\begin{split}
f_1(w; \lambda) &:= \hgF( (p-r)w+a+1, (q-r)w+b+1; 1-r w; x) \\
&\phantom{:=} + \dfrac{(p w+a)(q w+b) x}{r w(r w-1)} \cdot T(w) \cdot 
\hgF( (p-r)w+a+1, (q-r)w+b+1; 2-r w; x).    
\end{split}
\]
%%%%%%
Observe that any pole of $f_1(w; \lambda)$ is either a pole of 
$T(w)$ or in the discrete set $\frac{1}{r} \Z_{\ge0}$. 
Thus $f_1(w; \lambda)$ cannot have infinitely many 
poles off the positive real axis $\R_+$.  
On the other hand, applying Euler's transformation \cite[(7b)]{Iwasaki}  
to definition \eqref{eqn:f}, we have   
%%%%%%
\[
f(w;\lambda) = (1-x)^{(r-p-q)w-a-b} \hgF( (r-p)w-a, (r-q)w-b; r w; x), 
\]
%%%%%% 
where $r-p > 0$, $r-q > 0$ and $r > 0$ while $0 < x < 1$ 
by the assumption $\lambda \in \cF^-$.  
Thus 
%%%%%%%%%%%%%%%%%%%%%%%%% eqn:f>0 %%%%%%%%%%%%%%%%%%%%%%%%%%%%%%%%%%%%%%
\begin{equation} \label{eqn:f>0}
f(w;\lambda) > 0 \qquad \mbox{for every real} \quad  
w > \max \{a/(r-p), \, b/(r-q), \, 0\}, 
\end{equation}
%%%%%%%%%%%%%%%%%%%%%%%%%%%%%%%%%%%%%%%%%%%%%%%%%%%%%%%%%%%%%%%%%%%%%%%% 
so $f(w;\lambda)$ cannot have infinitely many zeros on $\R_+$.   
Assertion (2) then implies that $f(w; \lambda)$ must have  
infinitely many zeros off $\R_+$.  
Therefore $f(w; \lambda)$ admits a zero $w_0 \in \C \setminus \R_+$ 
that is not a pole of $f_1(w; \lambda)$. 
Substituting $w = w_0$ into equation \eqref{eqn:f-f1} yields a 
contradiction $0 = f(w_0; \lambda) f_1(w_0; \lambda) = 
(1-x)^{(r-p-q)w_0-a-b} \neq 0$, which shows that 
$f(w; \lambda)$ and $\tilde{f}(w; \lambda)$ are linearly independent 
over the rational function field $\C(w)$. \hfill $\Box$ 
%%%%%%%%%%%%%%%%%%%%%%%% end proof %%%%%%%%%%%%%%%%%%%%%%%%%%%%%%%%%%%%%%
%%%%%%%%%%%%%%%%%%%%%%%% prop:cfcr-F %%%%%%%%%%%%%%%%%%%%%%%%%%%%%%%%%%%%
\begin{proposition} \label{prop:cfcr-F}
Let $\lambda = (p,q,r;a,b;x) \in \cF^-$ in what follows.   
\begin{enumerate}
\item Any solution to Problem $\rI$ or $\rII$ in $\cF^-$ is non-elementary,  
that is, $f(w; \lambda)$ has infinitely many poles in $\C_w$.  
\item Any integral solution in $\cF^-$ to Problem $\rII$ comes from 
contiguous relations. 
\item For any integral solution $\lambda \in \cF^-$ to Problem $\rII$, 
its reciprocal $\check{\lambda}$ is an $(\rA)$-solution in $\cD^-$ 
so that $r \equiv 0 \mod 2$ and $\lambda = \lambda^{\vee\vee}$ 
is the reciprocal of an $(\rA)$-solution $\check{\lambda} \in \cD^-$.     
\item Any rational solution in $\cF^-$ to Problem $\rII$ 
essentially comes from contiguous relations. 
\end{enumerate} 
\end{proposition}
%%%%%%%%%%%%%%%%%%%%%%%%%%%%%%%%%%%%%%%%%%%%%%%%%%%%%%%%%%%%%%%%%%%%%%%%%
%%%%%%%%%%%%%%%%%%%%%%%% begin proof %%%%%%%%%%%%%%%%%%%%%%%%%%%%%%%%%%%%
{\it Proof}. 
Assertion (1) follows directly from assertion (1) of Lemma \ref{lem:zeros}. 
With assertion (3) of Lemma \ref{lem:zeros} the proof of assertion (2) 
is exactly the same as that of \cite[Theorem 7.2]{Iwasaki}.  
The proof of assertion (3) proceeds as follows.  
Since the $\lambda$ in assertion (3) comes from contiguous relations by 
assertion (2), Lemma \ref{lem:r-cf1} can be used to infer that its 
reciprocal $\check{\lambda}$ is a solution to Problem $\rII$ and hence 
becomes an $(\rA)$-solution in $\cD^-$ by \cite[Theorems 2.1 and 2.2]{Iwasaki}. 
In other words, $\lambda = \lambda^{\vee\vee}$ is the reciprocal of an 
$(\rA)$-solution $\check{\lambda} \in \cD^-$. 
Applying assertion (2) of \cite[Theorem 2.2]{Iwasaki} to 
$\check{\lambda}$ we have 
$r = \check{r} - \check{p} - \check{q} \equiv 0 \mod 2$. 
Assertion (4) is an immediate consequence of assertion (2) applied to 
a multiplication $k \lambda$, where $k \in \N$ is chosen so that 
$k \lambda$ becomes integral. \hfill $\Box$
%%%%%%%%%%%%%%%%%%%%%%%% end proof %%%%%%%%%%%%%%%%%%%%%%%%%%%%%%%%%%%%%%
%%%%%%%%%%%%%%%%%%%%%%%% subsec:gpf-F %%%%%%%%%%%%%%%%%%%%%%%%%%%%%%%%%%%
\subsection{Gamma Product Formulas} \label{subsec:gpf-F} 
%%%%%%%%%%%%%%%%%%%%%%%%%%%%%%%%%%%%%%%%%%%%%%%%%%%%%%%%%%%%%%%%%%%%%%%%%
Now that Proposition \ref{prop:cfcr-F} is established, we can use 
reciprocity and multiplication to consider  whether any {\sl rational} 
solution to Problem $\rII$ in $\cF^-$ leads back to a solution to Problem $\rI$.    
%%%%%%%%%%%%%%%%%%%%%%%% lem:gpf-F %%%%%%%%%%%%%%%%%%%%%%%%%%%%%%%%%%%%%%
\begin{lemma} \label{lem:gpf-F} 
If $\lambda = (p,q,r;a,b;x) \in \cF^-$ is a rational solution to 
Problem $\rII$ with 
%%%%%%%%%%%%%%%%%%%%%%% eqn:ocf-F %%%%%%%%%%%%%%%%%%%%%%%%%%%%%%%%%%%%%%%
\begin{equation} \label{eqn:ocf-F}
\dfrac{f(w+1; \lambda)}{f(w; \lambda)} = d 
\dfrac{(w+u_1)\cdots(w+u_m)}{(w+v_1)\cdots(w+v_n)},   
\end{equation}
%%%%%%%%%%%%%%%%%%%%%%%%%%%%%%%%%%%%%%%%%%%%%%%%%%%%%%%%%%%%%%%%%%%%%%%%%
for some $d \in \C^{\times}$, $u_1, \dots, u_m \in \C$ and 
$v_1, \dots, v_m \in \C$, then $m = n$, $d$ is a positive number given 
by formula \eqref{eqn:d-F}, and $\lambda$ becomes a solution to Problem 
$\rI$ with gamma product formula     
%%%%%%%%%%%%%%%%%%%%%%% eqn:gpf-F2 %%%%%%%%%%%%%%%%%%%%%%%%%%%%%%%%%%%%%%
\begin{equation} \label{eqn:gpf-F2}
f(w; \lambda) = C \cdot d^w \cdot 
\dfrac{\varGamma(w+u_1) \cdots \varGamma(w+u_m)}{\varGamma(w+v_1) 
\cdots \varGamma(w+v_m) },
\end{equation}
%%%%%%%%%%%%%%%%%%%%%%%%%%%%%%%%%%%%%%%%%%%%%%%%%%%%%%%%%%%%%%%%%%%%%%%%
where $C$ is a positive constant. 
\end{lemma}
%%%%%%%%%%%%%%%%%%%%%%%%%%%%%%%%%%%%%%%%%%%%%%%%%%%%%%%%%%%%%%%%%%%%%%%%%
%%%%%%%%%%%%%%%%%%%%%%%% begin proof %%%%%%%%%%%%%%%%%%%%%%%%%%%%%%%%%%%%
{\it Proof}.  
Let $k \in \N$ be the least common denominator of $p$, $q$, $r \in \Q$. 
Then the multiplication $k \lambda = (k p, k q, kr; a, b; x) \in \cF^-$ 
is an integral solution to Problem $\rII$.   
It follows from assertion (3) of Proposition \ref{prop:cfcr-F} that 
$(k \lambda)^{\vee}$ is an $(\rA)$-solution in $\cD^-$ and $k \lambda 
= (k \lambda)^{\vee\vee}$ is its reciprocal.   
Accordingly we can apply Proposition \ref{prop:division} to the 
$(\rA)$-solution $(k \lambda)^{\vee} = (k|p|, k|q|, k(r-p-q); \, \check{a}, 
\check{b}; \, \check{x}) \in \cD^-$. 
Adapting GPF \eqref{eqn:gpf-c2} to the current situation, we have 
%%%%%%%%%%%%%%%%%%%%%%%%%%%%%%% eqn:gpf-F3 %%%%%%%%%%%%%%%%%%%%%%%%%%%%%%
\begin{equation} \label{eqn:gpf-F3}
f(w; k \lambda) = f\left(w; (k \lambda)^{\vee\vee} \right) 
= \check{C} \cdot \delta^{k w} \cdot \dfrac{\prod_{i=0}^{k r-1} 
\varGamma\left(w + \frac{i}{k r} \right)}{\prod_{i=1}^{k r} 
\varGamma(w+ \xi_i)}, 
\end{equation}
%%%%%%%%%%%%%%%%%%%%%%%%%%%%%%%%%%%%%%%%%%%%%%%%%%%%%%%%%%%%%%%%%%%%%%%%%
where $r-p-q$ in the general formula \eqref{eqn:gpf-c2} is now $k r$, 
the constant $\check{d}$ there is now $\delta^k$ with 
%%%%%%%%%%%%%%%%%%%%%%%%%%%%%% eqn:d-F2 %%%%%%%%%%%%%%%%%%%%%%%%%%%%%%%%%
\begin{equation} \label{eqn:d-F2}
\delta := r^r \sqrt{\dfrac{|p|^{|p|} \, |q|^{|q|} \, 
(1-x)^{r-p-q}}{(r-p)^{r-p} \, (r-q)^{r-q} \, x^r}},   
\end{equation}
%%%%%%%%%%%%%%%%%%%%%%%%%%%%%%%%%%%%%%%%%%%%%%%%%%%%%%%%%%%%%%%%%%%%%%%%%
(which is different from $\delta$ in \eqref{eqn:D-d}), while the numbers 
$\check{v}_i$ there are now written $\xi_i$.  
Under these circumstances the conditions \eqref{eqn:sum2} and 
\eqref{eqn:estimate3} with $\lambda$ replaced by $(k \lambda)^{\vee}$ 
are represented as  
%%%%%%%%%%%%%%%%%% eqn:sum3 %%%%%% eqn:estimate4 %%%%%%%%%%%%%%%%%%%%%%%%
\begin{gather} 
\xi_1 + \cdots + \xi_{k r} = \frac{k r -1}{2}, \label{eqn:sum3} \\
\xi_1, \dots, \xi_{k r} \in \Q; \qquad 
\frac{c}{k} \le \xi_1, \dots, \xi_{k r} < \frac{c}{k} + 1, \qquad 
c = c(\lambda) := \dfrac{1-a-b}{r-p-q},  \label{eqn:estimate4} 
\end{gather}
%%%%%%%%%%%%%%%%%%%%%%%%%%%%%%%%%%%%%%%%%%%%%%%%%%%%%%%%%%%%%%%%%%%%%%%%% 
where we used $c\left((k \lambda)^{\vee}\right) = - c(k \lambda) = 
- c(\lambda)/k = -c/k$ to derive condition \eqref{eqn:estimate4}.    
%%%%%%
\par
%%%%%%
On the one hand formula \eqref{eqn:gpf-F3} leads to a closed-form 
expression 
%%%%%%%%%%%%%%%%%%%%%%%%%%%%%%% eqn:ocf-F2 %%%%%%%%%%%%%%%%%%%%%%%%%%%%%%
\begin{equation} \label{eqn:ocf-F2}
\frac{f(w+1; k \lambda)}{f(w; k \lambda)} 
= \delta^k \cdot \dfrac{\prod_{i=0}^{k r-1} 
\left( w + \frac{i}{k r} \right)}{\prod_{i=1}^{k r} (w+ \xi_i) },  
\end{equation}
%%%%%%%%%%%%%%%%%%%%%%%%%%%%%%%%%%%%%%%%%%%%%%%%%%%%%%%%%%%%%%%%%%%%%%%%%
whereas on the other hand applying formula \eqref{eqn:R-mult} to 
assumption \eqref{eqn:ocf-F} yields 
%%%%%%%%%%%%%%%%%%%%%%%%%%%%% eqn:ocf-F3 %%%%%%%%%%%%%%%%%%%%%%%%%%%%%%%%
\begin{equation} \label{eqn:ocf-F3}
\frac{f(w+1; k \lambda)}{f(w; k \lambda)} 
= d^k \cdot k^{k(m-n)} \cdot 
\prod_{i=0}^{k-1} 
\dfrac{\left(w+ \frac{u_1 + i}{k}\right) \cdots 
\left(w+ \frac{u_m +i}{k}\right)}{\left(w+ \frac{v_1 + i}{k}\right) 
\cdots \left(w+ \frac{v_n +i}{k}\right)}. 
\end{equation}
%%%%%%%%%%%%%%%%%%%%%%%%%%%%%%%%%%%%%%%%%%%%%%%%%%%%%%%%%%%%%%%%%%%%%%%%%%
Since \eqref{eqn:ocf-F2} and \eqref{eqn:ocf-F3} must be exactly the same,   
they must be so asymptotically, that is, 
$\delta^k \sim d^k \cdot k^{k(m-n)} \cdot w^{k(m-n)}$ as $w \to \infty$,  
which forces $m = n$ and $\delta^k = d^k$. 
Formula \eqref{eqn:ocf-F} together with \eqref{eqn:f>0} and $m = n$ 
yields $0 < f(w+1; \lambda)/f(w; \lambda) \to d$ as $\R \ni w \to +\infty$ 
and so $d \ge 0$.   
Since $\delta$ is positive by definition \eqref{eqn:d-F2}, 
equation $\delta^k = d^k$ gives $d = \delta$.  
Thus in view of definition \eqref{eqn:d-F2}, $d$ must be given 
by formula \eqref{eqn:d-F}. 
The coincidence of \eqref{eqn:ocf-F2} and \eqref{eqn:ocf-F3} yields  
%%%%%%
\[
\dfrac{\prod_{i=0}^{k r-1} 
\left( w + \frac{i}{k r} \right)}{\prod_{i=1}^{k r} (w+ \xi_i) }
= 
\prod_{i=0}^{k-1} 
\dfrac{\left(w+ \frac{u_1 + i}{k}\right) \cdots 
\left(w+ \frac{u_m +i}{k}\right)}{\left(w+ \frac{v_1 + i}{k}\right) 
\cdots \left(w+ \frac{v_m +i}{k}\right)},  
\]
%%%%%% 
so that gamma product formula \eqref{eqn:gpf-F3} can be recast to  
%%%%%%%%%%%%%%%%%%%%%%%%%%%%%%% eqn:gpf-F4 %%%%%%%%%%%%%%%%%%%%%%%%%%%%%%
\begin{equation} \label{eqn:gpf-F4}
f(w; k \lambda) = \check{C} \cdot d^{k w} \cdot \prod_{i=0}^{k-1} 
\dfrac{\varGamma\left(w+ \frac{u_1 + i}{k}\right) \cdots 
\varGamma\left(w+ \frac{u_m +i}{k}\right)}{\varGamma 
\left(w+ \frac{v_1 + i}{k}\right) 
\cdots \varGamma \left(w+ \frac{v_m +i}{k}\right)}, 
\end{equation}
%%%%%%%%%%%%%%%%%%%%%%%%%%%%%%%%%%%%%%%%%%%%%%%%%%%%%%%%%%%%%%%%%%%%%%%%%
where $m = n$ and $d = \delta$ are also incorporated. 
Replacing $w$ by $w/k$ in \eqref{eqn:gpf-F4} and using the multiplication 
formula \eqref{eqn:mult} for the gamma function, we obtain formula 
\eqref{eqn:gpf-F2} with $C = \check{C} \cdot k^{v-u}$, where 
$u := u_1+\cdots+u_m$ and $v := v_1+\cdots+v_m$ are real numbers 
as $\lambda$ is a real data.   
Since $\check{C}$ is positive, so is the constant $C$.  \hfill $\Box$
%%%%%%%%%%%%%%%%%%%%%%%% end proof %%%%%%%%%%%%%%%%%%%%%%%%%%%%%%%%%%%%%%
%%%%%%%%%%%%%%%%%%%%%%%% prop:gpf-F %%%%%%%%%%%%%%%%%%%%%%%%%%%%%%%%%%%%%
\begin{proposition} \label{prop:gpf-F} 
If $\lambda = (p,q,r;a,b;x) \in \cF^-$ is a rational solution to Problem 
$\rII$, then $r \in \N$, $a, b \in \Q$, $x$ algebraic, and $\lambda$ is 
a solution to Problem $\rI$ with GPF 
\eqref{eqn:gpf-F} where $C$ is a positive constant, $d$ is given by 
formula \eqref{eqn:d-F} and $v_1, \dots, v_r$ satisfy condition 
\eqref{eqn:vi-F}.  
\end{proposition}
%%%%%%%%%%%%%%%%%%%%%%%%%%%%%%%%%%%%%%%%%%%%%%%%%%%%%%%%%%%%%%%%%%%%%%%%% 
%%%%%%%%%%%%%%%%%%%%%%%% begin proof %%%%%%%%%%%%%%%%%%%%%%%%%%%%%%%%%%%%
{\it Proof}. 
By assertion (1) of Lemma \ref{lem:zeros} the function $f(w; \lambda)$ 
has infinitely many poles so we must have $m \ge 1$ in formula 
\eqref{eqn:gpf-F2}. 
There exists an integer $s$ with $0 \le s \le m$ such that 
%%%%%%%%%%%%%%%%%
\begin{center}
(i) \,\, $u_i - v_j \not\in \Z$ \, for any \, $i, j = 1, \dots, s$, 
\qquad   
(ii) \,\, $u_i - v_i \in \Z$ \, for any \, $i = s+1, \dots, m$, 
\end{center}
%%%%%%%%%%%%%%%%% 
after suitable rearrangements of $\{u_1, \dots, u_m\}$ and 
$\{v_1, \dots, v_m \}$, where condition (i) resp. (ii) should be 
ignored if $s = 0$ resp. $s = m$. 
In view of property (ii) a repeated use of the recursion formula 
$\varGamma(w+1) = w \varGamma(w)$ allows us to rewrite formula 
\eqref{eqn:gpf-F2} as 
%%%%%%%%%%%%%%%%%%%%%%%% eqn:gpf-F5 %%%%%%%%%%%%%%%%%%%%%%%%%%%%%%%%%%%%%
\begin{equation} \label{eqn:gpf-F5}
f(w; \lambda) = S(w) \cdot d^w \cdot 
\dfrac{\varGamma(w+u_1) \cdots \varGamma(w+u_s)}{\varGamma(w+v_1) 
\cdots \varGamma(w+v_s) },
\end{equation}
%%%%%%%%%%%%%%%%%%%%%%%%%%%%%%%%%%%%%%%%%%%%%%%%%%%%%%%%%%%%%%%%%%%%%%%%%
where $S(w)$ is a nontrivial rational function of $w$. 
We must have $s \ge 1$, for otherwise $f(w; \lambda) = S(w) \cdot 
d^w$ could not have infinitely many poles, contradicting assertion (1) 
of Lemma \ref{lem:zeros}. 
%%%%%
\par
%%%%%
Take a nonnegative integer $i$ sufficiently large so that neither 
$\omega_0 := -u_1-i$ nor $\omega_1 := -u_1-i-1$ are zeros of $S(w)$. 
Then by property (i), $w = \omega_0$ and $w = \omega_1$ are actually 
poles of $f(w; \lambda)$, so by assertion (1) of Lemma \ref{lem:zeros} 
there exist nonnegative integers $j_0$ and $j_1$ with $j_0 < j_1$ 
such that $\omega_0 = -j_0/r$ and $\omega_1 = -j_1/r$. 
Thus we have $1 = \omega_0 - \omega_1 = (j_1 - j_0)/r$, namely,  
$r = j_1 - j_0 \in \N$ and hence the number $k$ in the proof of 
Lemma \ref{lem:gpf-F} is the least common denominator of $p$, $q \in \Q$.   
The assertion that $a$, $b \in \Q$ and $x$ is algebraic follows 
from the first part of Remark \ref{rem:division} applied to 
$(k \lambda)^{\vee}$ in place of $\lambda$. 
%%%%%
\par
%%%%%
Again by assertion (1) of Lemma \ref{lem:zeros}, for a sufficiently 
large $j_2$, the set  
%%%%%%
\[ 
\left\{-\frac{j}{r} \right\}_{j \ge j_2 \, r} = 
\coprod_{i=1}^r \left\{ -j- \frac{i-1}{r} \right\}_{j \ge j_2}  
\]
%%%%%%
constitutes all but a finite number of poles of $f(w; \lambda)$. 
The same is true with the multi-set
%%%%%%
\[ 
\bigcup_{i=1}^s \left\{ -j- u_i \right\}_{j \ge j_3} \qquad 
(\mbox{union as multi-sets}),    
\]
%%%%%% 
due to formula \eqref{eqn:gpf-F5} and property (i). 
Since all poles are simple, we have $u_i - u_j \not\in \Z$ for every 
distinct $i,j = 1, \dots, s$, so the union of multi-sets above is 
just a disjoint union of ordinary sets. 
As the two sets above can differ only by a finite number of elements, 
we must have 
%%%%%%
\[
r = s, \qquad u_i - (i-1)/r \in \Z \qquad (i = 1, \dots,r), 
\]
%%%%%% 
after taking a suitable rearrangement of $u_1, \dots, u_r$. 
Thus property (i) is equivalent to 
%%%%%%%%%%%%%%%%%%%%%%%% eqn:vi-F2 %%%%%%%%%%%%%%%%%%%%%%%%%%%%%%%%%%%%%%
\begin{equation} \label{eqn:vi-F2}
v_1, \dots, v_r \not\in \frac{1}{r} \, \Z,  
\end{equation}
%%%%%%%%%%%%%%%%%%%%%%%%%%%%%%%%%%%%%%%%%%%%%%%%%%%%%%%%%%%%%%%%%%%%%%%%%  
and a further repeated use of the recursion formula 
$\varGamma(w+1) = w \varGamma(w)$ converts \eqref{eqn:gpf-F5} to 
%%%%%%%%%%%%%%%%%%%%%%%% eqn:gpf-F6 %%%%%%%%%%%%%%%%%%%%%%%%%%%%%%%%%%%%%
\begin{equation} \label{eqn:gpf-F6}
f(w; \lambda) = S(w) \cdot d^w \cdot \dfrac{ \prod_{i=0}^{r-1} 
\varGamma\left(w+\frac{i}{r} \right)}{ \prod_{i=1}^r 
\varGamma(w+v_i)},
\end{equation}
%%%%%%%%%%%%%%%%%%%%%%%%%%%%%%%%%%%%%%%%%%%%%%%%%%%%%%%%%%%%%%%%%%%%%%%%%
with a possibly different rational function $S(w)$ of $w$, where we may 
assume   
%%%%%%%%%%%%%%%%%%%%%%% eqn:vi-F3 %%%%%%%%%%%%%%%%%%%%%%%%%%%%%%%%%%%%%%%
\begin{equation} \label{eqn:vi-F3}
c \le \mathrm{Re} \, v_1, \dots, \mathrm{Re}\, v_r < c + 1, 
\qquad c := \frac{1-a-b}{r-p-q},  
\end{equation}
%%%%%%%%%%%%%%%%%%%%%%%%%%%%%%%%%%%%%%%%%%%%%%%%%%%%%%%%%%%%%%%%%%%%%%%%% 
after translating $v_1, \dots, v_r$ by suitable integers and making yet 
another use of the recursion formula $\varGamma(w+1) = w \varGamma(w)$ 
with an ensuing modification of the rational function $S(w)$. 
%%%%%
\par
%%%%%
Replacing $w$ by $k w$ in formula \eqref{eqn:gpf-F6} and using 
formulas \eqref{eqn:f-mult} and \eqref{eqn:mult}, we have 
%%%%%%%%%%%%%%%%%%%%%%%% eqn:gpf-F7 %%%%%%%%%%%%%%%%%%%%%%%%%%%%%%%%%%%%%
\begin{equation} \label{eqn:gpf-F7}
f(w; k \lambda) = k^{(r-1)/2-v} \cdot S(k w) \cdot d^{k w} \cdot 
\dfrac{ \prod_{i=0}^{k r-1} \varGamma\left(w+\frac{i}{k r} \right)}{ 
\prod_{i=1}^{k r} \varGamma(w+\eta_i)}, 
\end{equation}
%%%%%%%%%%%%%%%%%%%%%%%%%%%%%%%%%%%%%%%%%%%%%%%%%%%%%%%%%%%%%%%%%%%%%%%%% 
where $v := v_1 + \cdots + v_r$ and $\eta_1, \dots, \eta_{k r}$ are the 
numbers defined by the multi-set 
%%%%%%%%%%%%%%%%%%%%%%%% eqn:eta %%%%%%%%%%%%%%%%%%%%%%%%%%%%%%%%%%%%%%%%
\begin{equation} \label{eqn:eta} 
\{\eta_1, \dots, \eta_{k r}\} := 
\{v_{i,j} \,:\, i = 1, \dots, r, \, j = 0, \dots, k-1\}, \qquad  
v_{i,j} := \frac{v_i + j}{k}. 
\end{equation} 
%%%%%%%%%%%%%%%%%%%%%%%%%%%%%%%%%%%%%%%%%%%%%%%%%%%%%%%%%%%%%%%%%%%%%%%%%
In view of this definition the estimate \eqref{eqn:vi-F3} leads to 
%%%%%%%%%%%%%%%%%%%%%%% eqn:vi-F4 %%%%%%%%%%%%%%%%%%%%%%%%%%%%%%%%%%%%%%%
\begin{equation} \label{eqn:vi-F4}
\frac{c}{k} \le \mathrm{Re} \, \eta_1, \dots, 
\mathrm{Re} \, \eta_{k r} < \frac{c}{k} + 1.       
\end{equation}
%%%%%%%%%%%%%%%%%%%%%%%%%%%%%%%%%%%%%%%%%%%%%%%%%%%%%%%%%%%%%%%%%%%%%%%%% 
Comparing formulas \eqref{eqn:gpf-F3} and \eqref{eqn:gpf-F7} with 
$\delta = d$ taken into account, we find 
that  
%%%%%%
\[ 
\dfrac{ \prod_{i=1}^{k r} \varGamma(w+ \eta_i)}{\prod_{i=1}^{k r} 
\varGamma(w+\xi_i)} = \left(\check{C} \right)^{-1} k^{(r-1)/2-v} 
\cdot S(k w),    
\]
%%%%%%%
which must be a rational function of $w$, having only at most finitely 
many poles and zeros. 
This forces $\eta_i - \xi_i\in \Z$ for every $i = 1, \dots, k r$, 
after a suitable rearrangement of $\xi_1, \dots, \xi_{k r}$. 
But in view of conditions \eqref{eqn:estimate4} and \eqref{eqn:vi-F4},  
this coincidence modulo $\Z$ must be an exact coincidence 
%%%%%%%%%%%%%%%%%%%%%%%%%% eqn:eta-xi %%%%%%%%%%%%%%%%%%%%%%%%%%%%%%%%%
\begin{equation} \label{eqn:eta-xi} 
\eta_i = \xi_i \qquad (i = 1, \dots, k r).  
\end{equation} 
%%%%%%%%%%%%%%%%%%%%%%%%%%%%%%%%%%%%%%%%%%%%%%%%%%%%%%%%%%%%%%%%%%%%%%%
So $\left(\check{C} \right)^{-1} k^{(r-1)/2-v} \cdot S(k w) = 1$  
and thus $S(w) = \check{C} k^{v-(r-1)/2} =: C$ must be a positive 
constant and formula \eqref{eqn:gpf-F6} reduces to GPF 
\eqref{eqn:gpf-F}. 
We can arrange $\eta_1, \dots, \eta_{k r}$ so that 
$\eta_i = v_{i,0} = v_i/k$ for each $i = 1, \dots, r$. 
Then coincidence \eqref{eqn:eta-xi} and the rationality in 
\eqref{eqn:estimate4} imply $v_i = k \eta_i = k \xi_i \in \Q$ for 
every $i = 1, \dots, r$. 
This together with \eqref{eqn:vi-F2} and \eqref{eqn:vi-F3} gives 
$v_1, \dots, v_r \in \left(\Q \setminus \frac{1}{r} \, \Z\right)  
\cap \left[c, \, c+1\right)$ in conditions \eqref{eqn:vi-F}. 
To prove the remaining condition in \eqref{eqn:vi-F} we observe that 
%%%%%%%%%%%%%%%%%%%%%%%%% eqn:xi-eta %%%%%%%%%%%%%%%%%%%%%%%%%%%%%%%%%%%%
\begin{equation} \label{eqn:xi-eta}
\sum_{i=1}^{k r} \xi_i = \sum_{i=1}^{k r} \eta_i = 
\sum_{i=1}^r \sum_{j=0}^{k-1} v_{i,j} \\
= \sum_{i=1}^r \sum_{j=0}^{k-1} \frac{v_i + j}{k} = 
\sum_{i=1}^r v_i + \frac{(k-1) r}{2}, 
\end{equation}
%%%%%%%%%%%%%%%%%%%%%%%%%%%%%%%%%%%%%%%%%%%%%%%%%%%%%%%%%%%%%%%%%%%%%%%%% 
where we used \eqref{eqn:eta-xi} in the first equality and 
\eqref{eqn:eta} in the second and third equalities. 
Comparing \eqref{eqn:xi-eta} with \eqref{eqn:sum3} 
we have $v_1 + \cdots + v_r = (r-1)/2$ and hence proves conditions 
\eqref{eqn:vi-F}. \hfill $\Box$ \par\medskip 
%%%%%%%%%%%%%%%%%%%%%%%% end proof %%%%%%%%%%%%%%%%%%%%%%%%%%%%%%%%%%%%%%
It is evident that Theorem \ref{thm:F} follows from Propositions 
\ref{prop:cfcr-F} and \ref{prop:gpf-F}. 
The proof of Theorem \ref{thm:sol-B} is virtually contained in the proof  
of Proposition \ref{prop:gpf-F} as a special case $k = 2$; the only 
necessary modification is to deduce the condition 
$0 \le v_1, \dots, v_r < 1$ for a $(\rB)$-solution from a similar 
condition for an $(\rA)$-solution.  
%%%%%%%%%%%%%%%%%%%%%%%%% subsec:sol-F %%%%%%%%%%%%%%%%%%%%%%%%%%%%%%%%%%%%%%%%%
\subsection{Examples of Reciprocal Solutions} \label{subsec:sol-F}
%%%%%%%%%%%%%%%%%%%%%%%%%%%%%%%%%%%%%%%%%%%%%%%%%%%%%%%%%%%%%%%%%%%%%%%%%%%%%
%%%%%%%%%%%%%%%%%%%%%%%% tab:sol-F %%%%%%%%%%%%%%%%%%%%%%%%%%%%%%%%%%%%%%%%%%%
\begin{table}[t]
\[
\begin{array}{c|c|c|c|c|c|c|cccc|l|c}
\hline
  &    &    &             &                  &      &     &      &       &       &      &               &            \\[-4mm]
r & p  & q  & x           & d                & a    & b   & v_1  & v_2   & v_3   & v_4  & \mbox{remark} & \mbox{No.} \\[1mm]
\cline{1-13}
  &    &    &             &                 &              &             &              &               &  &  &                            &   \\[-3mm]  
2 & -1 & -1 & \frac{1}{9} & \frac{2^8}{3^5} & \frac{11}{8} & \frac{9}{8} & \frac{5}{24} & \frac{7}{24}  &  &  & \mbox{GS \cite[(6.6)]{GS}} & 1 \\[2mm]
\cline{6-13}
  &    &    &             &                 &              &             &              &               &  &  &                            &   \\[-3mm]   
  &    &    &             &                 & \frac{5}{8}  & \frac{7}{8} & \frac{1}{24} & \frac{11}{24} &  &  & \mbox{GS \cite[(6.5)]{GS}} & 2 \\[2mm]
\cline{6-13}
  &    &    &             &                 &              &             &              &               &  &  &                  &   \\[-3mm]   
  &    &    &             &                 & \frac{5}{4}  & \frac{3}{4} & \frac{1}{12} & \frac{5}{12}  &  &  & \mbox{self-dual} & 3 \\[2mm]  
\cline{1-13}
  &    &    &             &                    &             &             &              &              &               &               &                         &   \\[-3mm]   
4 & -1 & -1 & \frac{1}{5} & \frac{2^{14}}{5^6} & \frac{9}{8} & \frac{5}{8} & \frac{3}{40} & \frac{7}{40} & \frac{23}{40} & \frac{27}{40} & \mbox{divisible by $2$} & 4 \\[2mm]
\cline{6-13}
  &    &    &             &                    &             &             &              &              &               &               &                         &   \\[-3mm]   
  &    &    &             &                    & \frac{3}{8} & \frac{7}{8} & \frac{1}{40} & \frac{9}{40} & \frac{21}{40} & \frac{29}{40} & \mbox{divisible by $2$} & 5 \\[2mm]
\cline{1-13}
  &    &    &                           &                          &             &             &              &              &     & &                  &   \\[-3mm]   
2 & -2 & -2 & \frac{1}{4}(3 \sqrt{3}-5) & \frac{3^4}{2^7} \sqrt{3} & \frac{5}{3} & \frac{4}{3} & \frac{1}{12} & \frac{5}{12} &     & & \mbox{self-dual} & 6 \\[2mm]  
\cline{2-13}
  &    &    &                &                              &             &             &              &              &   &  &    &   \\[-3mm]   
  & -3 & -1 & 9 -4 \sqrt{5}  & \frac{2^8}{5^3}(5-2\sqrt{5}) & \frac{9}{4} & \frac{5}{4} & \frac{3}{20} & \frac{7}{20} &   &  &    & 7 \\[2mm]
\cline{6-13}
  &    &    &                &                              &             &             &              &              &   &  &    &   \\[-3mm]   
  &    &    &                &                              & \frac{7}{4} & \frac{3}{4} & \frac{1}{20} & \frac{9}{20} &   &  &    & 8 \\[2mm]
\cline{2-13}
  &    &    &                &                                   &             &             &              &              & & &                   &   \\[-3mm]   
  & -4 & -2 & 17-12 \sqrt{2} & \frac{2^{10}}{3^3}(17-12\sqrt{2}) & \frac{5}{2} & \frac{3}{2} & \frac{1}{12} & \frac{5}{12} & & &  \mbox{self-dual} & 9 \\[2mm]
\hline  
\end{array}
\]
\caption{Nine integral solutions in $\cF^-$ reciprocal to the $(\rA)$-solutions in \cite[Table 1]{Iwasaki}.} 
\label{tab:sol-F}
\end{table} 
%%%%%%%%%%%%%%%%%%%%%%%%%%%%%%%%%%%%%%%%%%%%%%%%%%%%%%%%%%%%%%%%%%%%%%%%%%%%%%%%%%%%%%%
Some examples of $(\rA)$-solutions in $\cD^-$ were presented in the previous article \cite[Table 1]{Iwasaki}. 
Their reciprocal solutions in $\cF^-$ are given in the corresponding places of 
Table \ref{tab:sol-F}, which exhibits the data $\lambda = (p,q,r;a,b;x) \in \cF^-$ 
itself as well as the values of $d$ and $v_1, \dots, v_r$ in GPF \eqref{eqn:gpf-F} 
with a brief remark, if any. 
Up to classical symmetries, terminating versions of solutions 1 and 2 are Gosper's 
conjectural identities in Gessel and Stanton \cite[formulas (6.6) and (6.5)]{GS} 
proved later by Karlsson \cite[formulas (1.5) and (1.4)]{Karlsson} (and also by 
Koepf \cite[Table 5]{Koepf} using an extension of WZ-method \cite{WZ}), 
while solution 3 can be found in 
Erd\'elyi et al. \cite[Chap. I\!I, \S 2.8, formula (54)]{Erdelyi}.  
Self-dual solutions are indicated so in the remark column, otherwise any pair of 
consecutive solutions is a dual pair. 
Solutions 4 and 5 with $(p,q;r) = (-1,-1;4)$ are divisible by $2$, 
halves of which are tabulated in Table \ref{tab:sol-F-half}. 
These two solutions in $\cF^-$ are ``dual" to each other as well as 
``reciprocal" to the $(\rB)$-solutions in \cite[Table 2]{Iwasaki}, 
although we have to be careful in applying these terminologies to 
non-integral solutions (as mentioned in Remark \ref{rem:kummer}). 
%%%%%%%%%%%%%%%%%%%%%%%% tab:sol-F %%%%%%%%%%%%%%%%%%%%%%%%%%%%%%%%%%%%%%%%%%%%%%
\begin{table}[t]
\[
\begin{array}{c|c|c|c|c|c|c|cc|c}
\hline
  &              &              &             &                    &             &             &              &               &            \\[-4mm]
r & p            & q            & x           & d                  &  a          & b           & v_1          & v_2           & \mbox{No.} \\[1mm]
\cline{1-10}
  &              &              &             &                    &             &             &              &               &   \\[-3mm]   
2 & -\frac{1}{2} & -\frac{1}{2} & \frac{1}{5}  & \frac{2^{7}}{5^3}  & \frac{9}{8} & \frac{5}{8} & \frac{3}{20} & \frac{7}{20} & 4 \\[2mm]
\cline{6-10}
  &              &              &             &                    &             &             &              &               &   \\[-3mm]   
  &              &              &             &                    & \frac{3}{8} & \frac{7}{8} & \frac{1}{20} & \frac{9}{20}  & 5 \\[2mm]
\hline  
\end{array}
\]
\caption{Two non-integral rational solutions in $\cF^-$.}  
\label{tab:sol-F-half}
\end{table} 
%%%%%%%%%%%%%%%%%%%%%%%%%%%%%%%%%%%%%%%%%%%%%%%%%%%%%%%%%%%%%%%%%%%%%%%%%%%%%%%%%%%%%%%  
%%%%%%%%%%%%%%%%%%%%%%%%%%%%% sec:algorithm %%%%%%%%%%%%%%%%%%%%%%%%%%%%%
\section{Algorithmic Point of View} \label{sec:algorithm}
%%%%%%%%%%%%%%%%%%%%%%%%%%%%%%%%%%%%%%%%%%%%%%%%%%%%%%%%%%%%%%%%%%%%%%%%% 
We put all the results on $\cD^- \cup \cF^-$ (in \S \ref{sec:results})  
into context from an algorithmic point of view and discuss how to 
enumerate all (rational) solutions $\lambda = (p,q,r;a,b,x)$ with a 
prescribed value of $\bp = (p,q;r)$, or perhaps how to prove the 
nonexistence of such solutions, both in finite steps.      
%%%%%%%%%%%%%%%%%%%%%%%%%%%%% subsec:integral %%%%%%%%%%%%%%%%%%%%%%%%%%%
\subsection{Integral Solutions} \label{subsec:integral} 
%%%%%%%%%%%%%%%%%%%%%%%%%%%%%%%%%%%%%%%%%%%%%%%%%%%%%%%%%%%%%%%%%%%%%%%%%  
If $\lambda = (p,q,r;a,b;x) \in \cD^-$ is an $(\rA)$-solution then 
$\bp = (p,q;r) \in D^-_{\rA}$ by assertion (2) of \cite[Theorem 2.2]{Iwasaki}. 
Putting \cite[Theorems 2.1 and 2.3]{Iwasaki} and Theorems \ref{thm:ab-rf} 
and \ref{thm:ab} together along with Theorem \ref{thm:truncate} below, 
we are able to develop an algorithm to find all $(\rA)$-solutions 
$\lambda \in \cD^-$ with any prescribed $\bp \in D^-_{\rA}$. 
To describe it we introduce some notation.  
Let $\langle \varphi(z) \rangle_k := \sum_{j=0}^k c_j \, z^j$ 
be the {\sl truncation} at degree $k$ of a power series $\varphi(z) = 
\sum_{j=0}^{\infty} c_j \, z^j$. 
For a data $\lambda \in \cD^-$ with $\bp \in D^-_{\rA}$ we consider 
two ``truncated hypergeometric products":       
%%%%%%%%%%%%%%%%%%
\begin{align*} 
V(w;\lambda) 
&:= (r w)_{r-1} \, \langle \, \hgF(\bal^*(w); z) \cdot 
\hgF(\bv-\bal^*(w+1);z) \, \rangle_{k} \big|_{z = x},  
\\[2mm]
P(w;\lambda) 
&:= (r w)_r \, \langle \, \hgF(\bal^*(w); z) \cdot 
\hgF(\1-\bal^*(w+1);z) \, \rangle_{k}\big|_{z = x},          
\end{align*}
%%%%%%%%%%%%%%%%%%
where $(s)_n := \varGamma(s+n)/\varGamma(s)$, $\bal^*(w) := 
((r-p)w-a, \, (r-q)w-b; \, r w)$, $\bv := (1,1;2)$, 
$\1 := (1,1;1)$ and $k := \max\{r-p-1, \, r-q-1\}$. 
It is shown in \cite[Lemma 10.2]{Iwasaki} that $V(w; \lambda)$ is a 
polynomial of degree at most $r-1$ in $w$ and so admits an expansion 
%%%%%
\[
V(w; \lambda) = \sum_{\nu=0}^{r-1} V_{\nu}(\lambda) \, w^{\nu} 
= \sum_{\nu=0}^{r-1} V_{\nu}(a,b;x) \, w^{\nu} ,  
\]
%%%%% 
where $V_{\nu}(\lambda)$ is written $V_{\nu}(a,b; x)$ when 
$\bp$ is understood to be given a priori.      
It is also known that $P(w; \lambda)$ is a polynomial of degree 
at most $r$. 
Our algorithm is based on the following.   
%%%%%%%%%%%%%%%%%%%%%%%%%%% thm:truncate %%%%%%%%%%%%%%%%%%%%%%%%%%%%%%
\begin{theorem} \label{thm:truncate} 
$\lambda \in \cD^-$ is an $(\rA)$-solution if and only if 
$\bp \in D^-_{\rA}$ and $V(w; \lambda)$ vanishes identically as a 
polynomial of $w$, that is, $(a,b;x)$ is a simultaneous root of 
algebraic equations    
%%%%%%%%%%%%%%%%%%%%%%%%%%% eqn:cfcr2 %%%%%%%%%%%%%%%%%%%%%%%%%%%%%%%%%
\begin{equation} \label{eqn:cfcr2}
V_{\nu}(a, b; x) = 0 \qquad (\nu = 0, \dots, r-1),  
\end{equation} 
%%%%%%%%%%%%%%%%%%%%%%%%%%%%%%%%%%%%%%%%%%%%%%%%%%%%%%%%%%%%%%%%%%%%%%%%
in which case $P(w; \lambda)$ is exactly of degree $r$ in $w$  
and $R(w; \lambda)$ in formula \eqref{eqn:ocf} is given by 
%%%%%%%%%%%%%%%%%%%%%%%%%%% eqn:R-l %%%%%%%%%%%%%%%%%%%%%%%%%%%%%%%
\begin{equation} \label{eqn:R-l}
R(w;\lambda) = (1-x)^{r-p-q-1} \cdot \dfrac{(r w)_r}{P(w;\lambda)}. 
\end{equation}
%%%%%%%%%%%%%%%%%%%%%%%%%%%%%%%%%%%%%%%%%%%%%%%%%%%%%%%%%%%%%%%%%%%%%%%%% 
\end{theorem} 
%%%%%%%%%%%%%%%%%%%%%%%%%%%%%%%%%%%%%%%%%%%%%%%%%%%%%%%%%%%%%%%%%%%%%%%%%
\par
%%%%%%
This theorem follows from \cite[Theorem 10.3]{Iwasaki}. 
Now we have the following.   
%%%%%%%%%%%%%%%%%%%%%%%%% algorithm %%%%%%%%%%%%%%%%%%%%%%%%%%%%%%%%%%%%%
\begin{algorithm} \label{algorithm}
To enumerate all $(\rA)$-solutions $\lambda \in \cD^-$ with any  
prescribed $\bp \in D^-_{\rA}$:     
\begin{enumerate}
\item Check if $\bp = (p,q;r)$ satisfies division relation \eqref{eqn:dr} 
in Theorem \ref{thm:ab-rf}.  
\item If it is alright, then following Theorem \ref{thm:ab} find a 
candidate for $\ba =(a,b)$ explicitly in terms of $\bp$, which must be 
in one of Cases 1--6 of Table \ref{tab:ab}.      
\item Substitute the ensuing candidate $\ba$ into \eqref{eqn:cfcr2} 
to derive algebraic equations for $x$ over $\Q$:       
%%%%%%%%%%%%%%%%%%%%%%%%% eqn:V %%%%%%%%%%%%%%%%%%%%%%%%%%%%%%%%%%%%%%%%%
\begin{equation} \label{eqn:V}
V_{\nu}(x) = 0 \qquad (\nu = 0, \dots, r-1),   
\end{equation}
%%%%%%%%%%%%%%%%%%%%%%%%%%%%%%%%%%%%%%%%%%%%%%%%%%%%%%%%%%%%%%%%%%%%%%%%%
where $V_{r-1}(x) = 0$ is equivalent to the algebraic equation $Y(x; \bp) = 0$ in 
\cite[Theorem 2.3]{Iwasaki}. 
\item Check if equations \eqref{eqn:V} admits a simultaneous root $x$ 
with $0 < x < 1$.  
If so, we actually get an $(\rA)$-solution $\lambda = (\bp; \ba; x) \in 
\cD^-$; otherwise the current candidate gives no solution.    
\item If we have a solution $\lambda$ in step (4), put it into 
formula \eqref{eqn:R-l} to find $R(w;\lambda)$ explicitly. 
If $R(w;\lambda)$ is of the form \eqref{eqn:R} (inevitably with $m = n$),  
then we have a GPF \eqref{eqn:gpf} by \cite[Theorem 2.1]{Iwasaki}, 
where the constant $C$ can be evaluated by putting   
$w = -a/p$ or $w = -b/q$ into GPF \eqref{eqn:gpf}.       
\item Repeat the procedures (2)--(5) until all candidates are exhausted. 
By Theorem \ref{thm:ab-rf} there are only a finite number of candidates   
so that the algorithm terminates in finite steps.   
\end{enumerate}
\end{algorithm}
%%%%%%%%%%%%%%%%%%%%%%%%%%%%%%%%%%%%%%%%%%%%%%%%%%%%%%%%%%%%%%%%%%%%%%%%%
\par
%%%%%
Empirically, it almost surely occurs that algebraic equations 
\eqref{eqn:V} have no roots in common for all candidates, 
in which case Algorithm \ref{algorithm} serves as a rigorous proof of 
the nonexistence of solutions in $\cD^-$ with given 
$\bp \in D^-_{\rA}$.     
Observe also that for any $\bp = (p,q;r) \in D_{\rA}^-$ with  
division relation \eqref{eqn:dr} there is an estimate 
$1 \le p, \, q \le r/2$.  
Thus by the use of Algorithm \ref{algorithm} the enumeration of all 
$(\rA)$-solutions in $\cD^-$ with any  prescribed bound for $r$ 
terminates in finite steps after producing at most finitely many solutions.    
In view of estimate \eqref{eqn:dr-bound} in Lemma \ref{lem:dr-bound},  
this is true even with any prescribed bound for $\check{r} := r-p-q$. 
For example, there are exactly seven $(\rA)$-solutions in $\cD^-$ with 
$\check{r} = 2$, all of which are contained in \cite[Table 1]{Iwasaki}.   
%%%%%
\par
%%%%% 
We proceed to the treatment of integral solutions in $\cF^-$, which 
is based on the following.  
%%%%%%%%%%%%%%%%%%%%%%%%%%%% thm:D-F %%%%%%%%%%%%%%%%%%%%%%%%%%%%%%%%%%%%
\begin{theorem} \label{thm:D-F} 
Reciprocity \eqref{eqn:recip} induces a bijection between the set of all 
$(\rA)$-solutions in $\cD^-$ and the set of all integral solutions in 
$\cF^-$. 
\end{theorem}
%%%%%%%%%%%%%%%%%%%%%%%%%%%%%%%%%%%%%%%%%%%%%%%%%%%%%%%%%%%%%%%%%%%%%%%%%
\par
%%%%
This theorem is an immediate consequence of Theorems \ref{thm:D-r} and 
\ref{thm:F}. 
It implies that all integral solutions $\lambda = (p,q,r;a,b;x) \in \cF^-$ 
with a given $\bp = (p,q;r)$ is in one-to-one correspondence with all 
$(\rA)$-solutions $\check{\lambda} = 
(\check{p}, \check{q}, \check{r}; \check{a}, \check{b}; \check{x}) \in \cD^-$ 
with a given $\check{\bp} = (\check{p}, \check{q}; \check{r}) = 
(-p,-q;r-p-q) = (|p|,|q|;r+|p|+|q|)$, so that the enumeration of former 
solutions is accomplished through the enumeration of latter solutions 
by the use of Algorithm \ref{algorithm}.  
%%%%%%%%%%%%%%%%%%%%%%%%%%%%%% cor:D-F %%%%%%%%%%%%%%%%%%%%%%%%%%%%%%%%%
\begin{corollary} \label{cor:D-F} 
There are at most finitely many integral solutions $\lambda = 
(p,q,r;a,b;x) \in \cF^-$ with any prescribed bound for $r$. 
Moreover any integral solution $\lambda \in \cF^-$ must satisfy  
%%%%%%%%%%%%%%% 
\[
-3 r \le p, \, q \le -1, \qquad -5 r \le p + q \le -2.  
\]
%%%%%%%%%%%%%%%
\end{corollary} 
%%%%%%%%%%%%%%%%%%%%%%%%%%%%%%%%%%%%%%%%%%%%%%%%%%%%%%%%%%%%%%%%%%%%%%%%
\par
%%%%%
These inequalities follow from estimates \eqref{eqn:dr-bound} in 
Lemma \ref{lem:dr-bound} by replacing $\lambda$ with $\check{\lambda}$.      
%%%%%%%%%%%%%%%%%%%%%%%%%%%% subsec:rational %%%%%%%%%%%%%%%%%%%%%%%%%%%%
\subsection{Rational Solutions} \label{subsec:rational}
%%%%%%%%%%%%%%%%%%%%%%%%%%%%%%%%%%%%%%%%%%%%%%%%%%%%%%%%%%%%%%%%%%%%%%%%%
We turn our attention to dealing with rational solutions in 
$\cD^- \cup \cF^-$ that are non-integral. 
In $\cD^-$ this amounts to considering solutions of type 
$(\rB)$, since all solutions in $\cD^-$ are rational and the 
dichotomy of types $(\rA)$ and $(\rB)$ is exactly that of 
`integral' and `non-integral'.     
%%%%% 
\par
%%%%%   
A solution $\lambda = (p, q, r; a, b; x)$ to Problem $\rI$ or $\rII$ 
is said to be {\sl divisible} by an integer $k \ge 2$, if the data 
$\lambda/k := (p/k, q/k, r/k; a, b; x)$ is also a solution to the 
same problem, in which case $\lambda/k$ is referred to as the 
{\sl division} of $\lambda$ by $k$, in particular, $\lambda/2$ 
is a {\sl half} of $\lambda$. 
%%%%%%%%%%%%%%%%%%%%%%%%%%%% lem:divisible %%%%%%%%%%%%%%%%%%%%%%%%%%%%%%
\begin{lemma} \label{lem:divisible}
For an integer $k \ge 2$, an integral solution $\lambda = (p,q,r;a,b;x) 
\in \cD^- \cup \cF^-$ is divisible by $k$ if and only if $k|r$ and there 
exist $d > 0$ and $v_1, \dots, v_s \in \Q$ with $s := r/k$ such that  
%%%%%%%%%%%%%%%%%%%%%%%%%%%% eqn:divisible %%%%%%%%%%%%%%%%%%%%%%%%%%%%%%
\begin{equation} \label{eqn:divisible}
R(w; \lambda) = d \cdot 
\dfrac{\prod_{i=0}^{r-1} \left(w+ \frac{i}{r} \right)}{\prod_{i=1}^{s} 
\prod_{j=0}^{k-1} \left(w+ \frac{v_i + j}{k} \right)}.  
\end{equation}
%%%%%%%%%%%%%%%%%%%%%%%%%%%%%%%%%%%%%%%%%%%%%%%%%%%%%%%%%%%%%%%%%%%%%%%%%
\end{lemma}
%%%%%%%%%%%%%%%%%%%%%%%%%%%%%%%%%%%%%%%%%%%%%%%%%%%%%%%%%%%%%%%%%%%%%%%%%
%%%%%%%%%%%%%%%%%%%%%%%%%%%% begin proof %%%%%%%%%%%%%%%%%%%%%%%%%%%%%%%%
{\it Proof}. 
Essential is the ``if" part. 
Since Problems $\rI$ and $\rII$ are equivalent for integral (or more 
generally rational) data in $\cD^- \cup \cF^-$ by \cite[Theorem 2.1]{Iwasaki} 
and assertion (3) of Theorem \ref{thm:F}, the closed-form condition 
\eqref{eqn:divisible} lifts to a ``multiplied" gamma product formula  
%%%%%%%%%
\[ 
f(w; \lambda) = C \cdot d^w \cdot 
\dfrac{\prod_{i=0}^{r-1} 
\varGamma \left(w+ \frac{i}{r} \right)}{\prod_{i=1}^{s} 
\prod_{j=0}^{k-1} \varGamma \left(w+ \frac{v_i + j}{k} \right)} 
= C \cdot d^w \cdot \prod_{j=0}^{k-1} 
\dfrac{\prod_{i=0}^{s-1}  
\varGamma \left(w+ \frac{(i/s) + j}{k} \right)}{\prod_{i=1}^{s} 
\varGamma \left(w+ \frac{v_i + j}{k} \right)},  
\]
%%%%%%%%%
which in turn descends through the multiplication formula 
\eqref{eqn:mult} to a GPF  
%%%%%%%%%
\[ 
f(w; \lambda/k) = f(w/k; \lambda) = C \cdot d^{w/k} \cdot 
\dfrac{ \prod_{i=0}^{s-1} 
\varGamma \left(w+ \frac{i}{s} \right)}{\prod_{i=1}^{s} 
\varGamma \left(w+ v_i \right)},   
\]
%%%%%%%%%
with a possibly different $C$. 
Thus $\lambda/k$ is a solution and hence $\lambda$ is divisible by $k$. 
The proof of ``only if" part is a simple reversing of the 
argument so far, with the positivity $d > 0$ coming from 
formulas \eqref{eqn:d} and \eqref{eqn:d-F}, while the 
rationality $v_1, \dots, v_s \in \Q$ from properties 
\eqref{eqn:estimate} and \eqref{eqn:vi-F}. 
\hfill $\Box$ \par\medskip 
%%%%%%%%%%%%%%%%%%%%%%%%%%%% end proof %%%%%%%%%%%%%%%%%%%%%%%%%%%%%%%%%%
By assertions (2) and (4) of \cite[Theorem 2.2]{Iwasaki} and 
assertion (3) of Theorem \ref{thm:F},  
%%%%%%%%%%%%%%%%%%
\begin{enumerate}
\item any $(\rB)$-solution in $\cD^-$ is a {\sl half} of an 
$(\rA)$-solution $\lambda_* = (p_*, q_*, r_*; a_*, b_*; x_*) 
\in \cD^-$ such that $p_*$ and $q_*$ are odd positive integers 
and $r_*$ is an even positive integer,   
\item any rational solution in $\cF^-$ is the division 
$\lambda_*/k$ of an integral solution $\lambda_* \in \cF^-$ 
by a {\sl divisor} $k$ of $r_*$, where $r_*$ is necessarily 
an even positive integer.     
\end{enumerate}
%%%%%%%%%%%%%%%%%  
Thus finding a rational solution consists of finding an integral   
solution $\lambda_*$ by using the recipe in \S \ref{subsec:integral} 
and verifying if $\lambda_*$ is divisible by a divisor of $r_*$ 
based on Lemma \ref{lem:divisible}. 
All rational solutions are obtained in this manner. 
Note that any seed solution $\lambda_*$ produces only 
a finite number of division solutions. 
All $(\rB)$-solutions $\lambda = (p,q,r;a,b;x) \in \cD^-$ with any 
prescribed bound for $r$ are only of a finite cardinality, since  
only division by $2$ is involved in this case.  
It is not known whether this is also the case for rational solutions 
in $\cF^-$, because $r_*$ and $k$ may be arbitrary large while 
$r = r_*/k$ is kept bounded.  
%%%%%%%%%%%%%%%%%%%%%%%%%% subsec:primitive %%%%%%%%%%%%%%%%%%%%%%%%%%%%%%%
\subsection{A Problem for Primitive Solutions} \label{subsec:primitive}
%%%%%%%%%%%%%%%%%%%%%%%%%%%%%%%%%%%%%%%%%%%%%%%%%%%%%%%%%%%%%%%%%%%%%%%%%
The main orientation of the field has been toward searching for as 
many solutions as possible, but the converse orientation of confining 
solutions into as slim a region as possible or perhaps of proving 
the nonexistence of solutions other than those already known is equally 
important toward the ultimate goal of complete enumeration. 
With this in mind we close this article by posing an important problem.  
A solution to Problem $\rI$ or $\rII$ is said to be {\sl primitive} 
if it is a multiplication of no other solution. 
Since any solution is a multiplication of a primitive solution, the 
enumeration of all solutions boils down to that of primitive ones.      
%%%%%%%%%%%%%%%%%%%%%%% prob:primitive %%%%%%%%%%%%%%%%%%%%%%%%%%%%%%%%%%
\begin{problem} \label{prob:primitive} 
Are there infinitely many {\sl primitive} solutions in domain $\cD^-$  
or only a finite number of them?   
Ask the same question for domain $\cE^{*-}$ or region $\cI^{*-}$  
or anywhere else.       
\end{problem}
%%%%%%%%%%%%%%%%%%%%%%%%%%%%%%%%%%%%%%%%%%%%%%%%%%%%%%%%%%%%%%%%%%%%%%
\par
%%%%%%  
As $r$ becomes larger, it is increasingly more difficult that the $r$ 
algebraic equations \eqref{eqn:V} admit a root in common, but there 
is no logical reasoning that prohibits this miracle.  
Solving Problem \ref{prob:primitive} is still very ambitious and 
seems to require a completely new and amazing idea.     
%%%%%%%%%%%%%%%%%%%%%%%%%%
\par\vspace{5mm} \noindent
%%%%%%%%%%%%%%%%%%%%%%%%%% 
{\bf Acknowledgments}. 
This work is supported by Grant-in-Aid for 
Scientific Research, JSPS, 16K05165 (C).        
%%%%%%%%%%%%%%%%%%%%%%%% References %%%%%%%%%%%%%%%%%%%%%%%%%%%%%%%%%%%%%

%%%%%%%%%%%%%%%%%%%%%%%%%%%%%%%%%%%%%%%%%%%%%%%%%%%%%%%%%%%%%%%%%%%%%%%
% \bibliography{indag-d-15-00307ver2}
%%%%%%%%%%%%%%%%%%%%%%%%%%%%%%%%%%%%%%%%%%%%%%%%%%%%%%%%%%%%%%%%%%%%%%%
\end{document}